\RequirePackage{fix-cm}
\documentclass[smallextended]{svjour3}       
\smartqed  

\usepackage{graphicx}

\usepackage{amsfonts, amssymb}
\usepackage{amsmath}

\numberwithin{equation}{section}

\begin{document}

\title{A Malmquist--Steinmetz theorem for difference equations \thanks{The first author is supported by a Project funded by China Postdoctoral Science Foundation~{(2020M680334)} and the Fundamental Research Funds for the Central Universities~{(FRF-TP-19-055A1)}. The second author is supported by a Visiting Professorship for Senior Foreign Experts by Ministry of Science and Technology of the People's Republic of China~{(G2021105019L)}. 
}
}


\author{ Yueyang Zhang  \and  Risto Korhonen        
}


\institute{ Y. Zhang \at School of Mathematics and Physics, University of Science and Technology Beijing, No.~30 Xueyuan Road, Haidian, Beijing, 100083, P.R. China \\
\email{zhangyueyang@ustb.edu.cn}           
\and
R. Korhonen \at Department of Physics and Mathematics, University of Eastern Finland, P.O. Box 111, FI-80101 Joensuu, Finland \\
\email{risto.korhonen@uef.fi}           
}

\date{Received: date / Accepted: date}

\maketitle

\begin{abstract}


It is shown that if the equation
    \begin{equation*}
    f(z+1)^n=R(z,f),
    \end{equation*}
where $R(z,f)$ is rational in both arguments and $\deg_f(R(z,f))\not=n$, has a transcendental meromorphic solution, then the equation above reduces into one out of several types of difference equations where the rational term $R(z,f)$ takes particular forms. Solutions of these equations are presented in terms of Weierstrass or Jacobian elliptic functions, exponential type functions or functions which are solutions to a certain autonomous first-order difference equation having meromorphic solutions with preassigned asymptotic behavior. These results complement our previous work on the case $\deg_f(R(z,f))=n$ of the equation above and thus provide a complete difference analogue of Steinmetz' generalization of Malmquist's theorem.



\end{abstract}

\section{Introduction}\label{intro} 

Nevanlinna theory (see, e.g., \cite{Hayman1964Meromorphic}) is a powerful tool when analyzing meromorphic solutions of complex differential equations. For example, by utilizing Nevanlinna theory, Yosida \cite{yosida:33} and Laine \cite{laine:71} provided elegant alternate proofs of the classical Malmquist theorem on first-order differential equations \cite{malmquist1913fonctions}, while Steinmetz \cite{steinmetz:78}, and Bank and Kaufman \cite{bank1980growth} gave a precise classification of the differential equation
\begin{equation}\label{first_order_di_n}
(f')^n=R(z,f),
\end{equation}
where $n\in \mathbb{N}$ and $R(z,f)$ is rational in both arguments. See also \cite[Chapter~10]{Laine1993} for Malmquist--Yosida--Steinmetz type theorems. In \cite{Korhonenzhang2020}, we studied a natural difference analogue of equation \eqref{first_order_di_n}, i.e., the first-order difference equation
\begin{equation}\label{first_order_de_n}
f(z+1)^n=R(z,f),
\end{equation}
where $n\in \mathbb{N}$ and $R(z,f)$ is rational in $f$ with small functions of $f$ as coefficients. Assuming that \eqref{first_order_de_n} has a meromorphic solution and $\deg_f(R(z,f))=n$, we showed that equation \eqref{first_order_de_n} can, by a bilinear transformation in $f$, be transformed into one in a list of twelve equations. In particular, we considered meromorphic solutions of hyper-order less than~1 of \eqref{first_order_de_n} and showed that if such a solution exists, then equation \eqref{first_order_de_n} with rational coefficients has to reduce into the difference linear or Riccati equation, or one in a list of five equations including Fermat type difference equations and a special case of the symmetric Quispel--Roberts--Thompson (QRT) map~\cite{QuispelRobertsThompson1988,QuispelRobertsThompson1989}. We also showed that these five equations are explicitly solved in terms of Weierstrass or Jacobian elliptic functions, or of functions which are solutions of certain difference Riccati equations. These results provide a natural difference analogue of Steinmetz' generalization of Malmquist's theorem in the sense of Ablowitz, Halburd and Herbst \cite{AblowitzHalburdHerbst2000}, who proposed that the existence of sufficiently many finite-order meromorphic solutions of a difference equation is a good difference analogue of the Painlev\'e property. Their idea was refined and successfully implemented by Halburd and the second author \cite{halburdrk:LMS2006} on the second-order difference equation
\begin{equation}\label{ahh_eq}
f(z+1)+f(z-1)=R(z,f),
\end{equation}
where $R(z,f)$ is rational in $f$ with small functions of $f$ as coefficients, reducing equation \eqref{ahh_eq} into a short list of canonical equations, including the difference Painl\'eve~I and~II equations. The finite-order condition of the proposed difference Painl\'eve property was relaxed into hyper-order strictly less than one by Halburd, the second author and Tohge \cite{halburdkt:14TAMS}, and recently into hyper-order equal to one with minimal hyper type by Zheng and the second author~\cite{zhengR:18}.


The purpose of this paper is to find out all transcendental meromorphic solutions for the case $\deg_f(R(z,f))\not=n$ of equation \eqref{first_order_de_n} without growth conditions and provide a complete difference analogue of Steinmetz' generalization of Malmquist's theorem. As is pointed out in \cite{Korhonenzhang2020}, in this case all the transcendental meromorphic solutions of \eqref{first_order_de_n} are of hyper-order at least one. Our work is a continuation of many mathematicians' research on first-order difference equations. For example, in \cite{shimomura:81}, Shimomura showed that the difference equation
    \begin{equation}\label{first_order shi eq}
    f(z+1) = P(f(z)),
    \end{equation}
where $P(f(z))$ is a polynomial in $f(z)$ with constant coefficients, always has a non-trivial entire solution; in \cite{yanagihara:80}, Yanagihara showed that the difference equation
    \begin{equation}\label{first_order_yan eq}
    f(z+1) = R(f(z)),
    \end{equation}
where $R(f(z))$ is rational in $f(z)$ having constant coefficients, has a non-trivial meromorphic solution with preassigned asymptotic behavior in a sector for all choices of $R\not\equiv0$. Yanagihara~\cite{yanagihara:80} also showed that if \eqref{first_order_yan eq}, where $R(f(z))$ is replaced by $R(z,f)$, which is rational in both arguments, has a transcendental meromorphic solution of hyper-order less than~1, then \eqref{first_order_yan eq} must reduce into the difference Riccati equation. This can be viewed as a natural difference analogue of Malmquist's theorem on first-order differential equations. For the higher-degree equation \eqref{first_order_de_n}, the classification work in the special case where the right-hand side (RHS) of \eqref{first_order_de_n} is a polynomial in $f(z)$ with constant coefficients has been done by Nakamura and Yanagihara~\cite{NakamuraYanagihara1989difference} and Yanagihara~\cite{Yanagihara1989difference}. The results obtained in the present paper, supplemented with results of the first part of this study~\cite{Korhonenzhang2020}, can be summarized as follows.

\begin{theorem}\label{completetheorem}
If equation \eqref{first_order_de_n}, where $R(z,f)$ is rational in both arguments, has a transcendental meromorphic solution, then \eqref{first_order_de_n} can be reduced into the case $n=1$ or one out of~$27$ equations in Theorems~\ref{Theorem  1a}--\ref{Theorem  2c} below and in \cite[Theorem~2]{Korhonenzhang2020}.
\end{theorem}

Here we have to point out that the list of equations in \cite[Theorem~2]{Korhonenzhang2020} is not complete, since one equation in the case $n=2$ is left out in the proof there; see~\cite{ZhangKorhonen:2022}. Equation \eqref{first_order_de_n} with $n=1$ actually includes~$3$ equations, namely the equation \eqref{yanagiharaeq0 fd} below and two equations in \cite[Theorem~2]{Korhonenzhang2020}. Autonomous versions of them are just \eqref{first_order shi eq} and \eqref{first_order_yan eq}. In other cases, we have counted each equation in Theorems~\ref{Theorem  1a}--\ref{Theorem  2c} below once even when some of them have the same form but appear in two different theorems. Moreover, autonomous versions of all the~$27$ equations can be solved in terms of elliptic and elementary functions.

In this paper, we shall confine ourselves to considering equation \eqref{first_order_de_n} with rational coefficients. The main tools from Nevanlinna theory we use, both in our previous paper~\cite{Korhonenzhang2020} and in the present one, are the generalizations of Nevanlinna's second main theorem given by Yamanoi \cite{yamanoi:04,yamanoi:05}. We refer to \cite[pp.~42--43]{Hayman1964Meromorphic} for the standard definitions of $\delta(a,f)$, $\theta(a,f)$ and $\Theta(a,f)$, etc. Recall that a value $a\in \mathbb{C}\cup\{\infty\}$ is said to be a \emph{completely ramified value} of $f(z)$ when $f(z)-a=0$ has no simple roots. Denote the field of rational functions by $\mathcal{R}$ and set $\hat{\mathcal{R}}=\mathcal{R}\cup\{\infty\}$. Throughout the paper, we say that $c(z)\in \hat{\mathcal{R}}$ is a \emph{completely ramified rational function} of a transcendental meromorphic function $f(z)$ when the equation $f(z)=c(z)$ has at most finitely many simple roots and that $c(z)$ is a \emph{Picard exceptional rational function} of $f(z)$ when $N(r,c,f)=O(\log r)$. We also say that $c(z)$ has multiplicity at least $m$ if all the roots of $f(z)=c(z)$ have multiplicity at least $m$ with at most finitely many exceptions. As is mentioned in \cite{Korhonenzhang2020}, the main theorem in \cite{yamanoi:05} yields that the inequality
\begin{equation}\label{multiplicityinequality}
\sum_{i=1}^q\Theta(c_i,f)\leq 2
\end{equation}
holds for any collection of $c_1,\cdots,c_q\in \hat{\mathcal{R}}$ when $f$ is transcendental. Moreover, we have
\begin{theorem}\label{completelyrm}
A non-constant transcendental meromorphic function $f(z)$ can have at most four completely ramified rational functions.
\end{theorem}
As in \cite{Korhonenzhang2020}, when considering a meromorphic solution $f(z)$ of \eqref{first_order_de_n}, we will do a transformation to $f$ using some algebraic functions and end up in a situation such that the considered functions are meromorphic on a finite-sheeted Riemann surface. Such functions are called \emph{algebroid functions} (see, e.g., \cite{Katajamaki1993algebroid}). We will ignore considering the degree of these functions in this paper since in any case we are dealing with at most finitely many algebraic branch points and the inequality \eqref{multiplicityinequality} and Theorem~\ref{completelyrm} hold true for such functions.

The remainder of this paper is organized in the following way. In section~\ref{Proof 0}, we will first set up some notation and build several lemmas concerning the roots of the numerator and the denominator of $R(z,f)$ in \eqref{first_order_de_n}. Our main results, i.e., Theorems~\ref{Theorem  1a}--\ref{Theorem  2c}, and the proofs for them are distributed in sections~\ref{Proof 1} and \ref{Proof 2}, respectively, where different cases of equation \eqref{first_order_de_n} are treated. Moreover, for all the equations we find we will present explicit solutions to them in the autonomous case. These solutions are presented in the corresponding theorems, however, for the elliptic solutions, in particular, we give a detailed discussion in section~\ref{discussions} separately.



For simplicity, from now on we will use the suppressed notations: $f=f(z)$, $\underline{f}=f(z-1)$ and $\overline{f}=f(z+1)$ for a meromorphic, algebraic or algebroid, function $f(z)$.

\section{Lemmas on the roots}\label{Proof 0}
We begin to consider the transcendental meromorphic solution $f$ of the following difference equation:
\begin{equation}\label{yanagiharaeq0}
\overline{f}^n=R(z,f),
\end{equation}
where $n\in \mathbb{N}$ and $R(z,f)$ is rational in both arguments and $\deg_{f} (R(z,f))\not=n$. We first set up some notation for \eqref{yanagiharaeq0} that will be used throughout the proofs in the following sections. We denote
    \begin{equation*}
    R(z,f) = \frac{P(z,f)}{Q(z,f)},
    \end{equation*}
where
    \begin{equation}\label{yanagiharaeq0    a}
    P(z,f) = a_pf^p + a_{p-1}f^{p-1}+\cdots+a_0
    \end{equation}
and
    \begin{equation}\label{yanagiharaeq0    b}
    Q(z,f) = b_qf^q + b_{q-1}f^{q-1}+\cdots+b_0
    \end{equation}
are two polynomials in $f$ having no common factors, $p,q\in\mathbb{N}$ and the coefficients $a_p$, $\cdots$, $a_0$, $b_q$, $\cdots$, $b_0$ are rational functions. We have $\deg_f (P(z,f))=p$ and $\deg_f (Q(z,f))=q$. Denote $d=\deg_{f}(R(z,f))$. Then $d=\max\{p,q\}\not=n$. If $n=1$, then equation \eqref{yanagiharaeq0} is just
\begin{equation}\label{yanagiharaeq0 fd}
\overline{f}=R(z,f).
\end{equation}
As mentioned in the introduction, in the autonomous case, equation \eqref{yanagiharaeq0 fd} always has an entire or meromorphic solution independently of the degree of $R$. Therefore, we always assume that $n\geq 2$ in what follows. Under this assumption, if the degree of $R(z,f)$ in $f$ in \eqref{yanagiharaeq0} equals~$1$, then we have
\begin{equation*}
\underline{f}=\frac{af^n+b}{cf^n+d},
\end{equation*}
where the coefficients are rational functions such that $ad-bc\not=0$. The equation above is included in \eqref{yanagiharaeq0 fd} under the transformation $x=-z$. So we also always assume
that $d\geq 2$. We write $P(z,f)$ and $Q(z,f)$ in the algebraic factorization form
    \begin{equation}\label{yanagiharaeq0    aert}
    P(z,f) = a_p\prod_{i=1}^{M}P_i(z,f)^{\mu_i},
    \end{equation}
and
    \begin{equation}\label{yanagiharaeq0    bert}
    Q(z,f) =b_q\prod_{j=1}^{N}Q_j(z,f)^{\nu_j},
    \end{equation}
where $P_i(z,f)$ and $Q_j(z,f)$ are irreducible polynomials in $f$ of the form in \eqref{yanagiharaeq0    a} or \eqref{yanagiharaeq0    b}, $\mu_i,\nu_j\in \mathbb{N}$ and $\mu_1+\cdots+\mu_M=p$ and $\nu_1+\cdots+\nu_N=q$. In this paper, we shall use the term 'a polynomial in $f$' which means that the polynomial in $f$ has rational or algebraic coefficients. Note that the inequality \eqref{multiplicityinequality} also holds when $f$ is replaced by an algebroid function with finitely many branch points and $c_1$, $\cdots$, $c_q$ are replaced by algebraic functions. For convenience, in the following we always use the terms 'completely ramified rational function' and 'Picard exceptional rational function' of $f$ even though sometimes they actually refer to algebraic functions. We also write \eqref{yanagiharaeq0    aert} and \eqref{yanagiharaeq0    bert} as
    \begin{equation}\label{P}
    P(z,f) = P_0(z,f)^n(f-\alpha_1)^{k_1}\cdots (f-\alpha_\mu)^{k_\mu}
    \end{equation}
and
    \begin{equation}\label{Q}
    Q(z,f)  =Q_0(z,f)^n(f-\beta_1)^{l_1}\cdots (f-\beta_\nu)^{l_\nu},
    \end{equation}
where $\alpha_1$, $\cdots$, $\alpha_\mu$, $\beta_1$, $\cdots$, $\beta_\nu$ are in general algebraic functions and $k_i$ and $l_j$ denote the orders of the roots $\alpha_i$ and $\beta_j$, respectively, and satisfy $n\nmid k_i$ and $n\nmid l_j$, and
    \begin{equation*}
    P_0(z,f) =a_p^{1/n}\prod_{i=1,n\mid \mu_i}^{M}P_i(z,f)^{\mu_i/n}=a_p^{1/n}\prod_{i=1,n\mid \mu_i}^{M}(f-\alpha_{\mu_{i,1}})^{\mu_i/n}\cdots(f-\alpha_{\mu_{i,p_i}})^{\mu_i/n}
    \end{equation*}
and
    \begin{equation*}
    Q_0(z,f) =b_q^{1/n}\prod_{j=1,n\mid \nu_j}^{N}Q_j(z,f)^{\nu_j/n}=b_p^{1/n}\prod_{j=1,n\mid \nu_j}^{N}(f-\beta_{\nu_{j,1}})^{\nu_j/n}\cdots(f-\beta_{\nu_{j,q_j}})^{\nu_j/n}
    \end{equation*}
are two polynomials in $f$ of degrees $p_0$ and $q_0$, respectively, with $\alpha_{\mu_{i,1}}$, $\cdots$, $\alpha_{\mu_{i,p_i}}$, $\beta_{\nu_{j,1}}$, $\cdots$, $\beta_{\nu_{j,q_j}}$ being in general algebraic functions, and $\mu_i$ and $\nu_j$ denoting the orders of these roots of $P(z,f)$ or $Q(z,f)$, respectively. We have $p_0=\deg_f(P_0(z,f))$ and $q_0=\deg_f(Q_0(z,f))$. Note that $a_p^{1/n},b_p^{1/n}$ are in general algebraic functions. For convenience, when $q\geq 1$, we always suppose that $b_q=1$. Also note that $\alpha_i$ and $\beta_j$ are neither roots of $P_0(z,f)$ nor roots of $Q_0(z,f)$. Now, if $n\nmid \mu_i$ for all $i=1,\cdots, M$ or $n\nmid \nu_j$ for all $j=1,\cdots,N$, then $P_0(z,f)^n=a_p$ or $Q_0(z,f)^n =b_q$. On the other hand, if there are no such $\mu_i$ or $\nu_j$ that $n\nmid \mu_i$ and $n\nmid \nu_j$, then after taking the $n$-th root on both sides of \eqref{yanagiharaeq0} we get equation \eqref{yanagiharaeq0 fd}. So in the sequel we always suppose that there is at least one such $\mu_i$ or $\nu_j$. The combined number of $\alpha_i$ and $\beta_j$ in \eqref{P} and \eqref{Q} is $\mu+\nu$ and, for convenience, we always denote $N_c:=\mu+\nu$ even when there is no such $\alpha_i$ or $\beta_j$. Moreover, we may suppose that the \emph{greatest common divisor} of $n$, $k_1$, $\cdots$, $k_\mu$, $l_1$, $\cdots$, $l_\nu$, which is denoted by $k=(n,k_1,\cdots,k_\mu,l_1,\cdots,l_\nu)$, is~$1$. Otherwise, after taking the $k$-th root on both sides of \eqref{yanagiharaeq0}, we get a new equation of the same form as \eqref{yanagiharaeq0} with the power of $\overline{f}$ being $n/k$.

With the notation above, we are now ready to construct lemmas regarding the roots of $P(z,f)$ and $Q(z,f)$ through elementary multiplicity analysis on $f$. For simplicity, in what follows, when considering the zeros, poles or $\alpha_i$-points of $f$, etc., we will omit giving the corresponding Taylor or Laurent series expansions for $f$. The first lemma below provides some basic upper bounds for $N_c$.

\begin{lemma}\label{basiclemma}
Let $f$ be a transcendental meromorphic solution of equation \eqref{yanagiharaeq0}. Then $\alpha_i$ is either a Picard exceptional rational function of $f$ or a completely ramified rational function of $f$ with multiplicity $n/(n,k_i)$ and $\beta_j$ is either a Picard exceptional rational function of $f$ or a completely ramified rational function of $f$ with multiplicity $n/(n,l_j)$. Moreover, if $q=0$, then $\infty$ is a Picard exceptional rational function of $f$ and $N_c\leq 2$ and, in particular, if $n\geq 3$, then $N_c=1$; if $q\geq 1$, then $N_c\leq 4$ and, in particular, if $n\geq 3$, then $N_c\leq 3$; if $q\geq1$ and $n\nmid|p-q|$, then $\infty$ is either a Picard exceptional rational function of $f$ or a completely ramified rational function of $f$ with multiplicity $n/(n,|p-q|)$ and $N_c\leq 3$.
\end{lemma}

\renewcommand{\proofname}{Proof.}
\begin{proof}
By making use of the factorizations \eqref{P} and \eqref{Q}, it follows that $\alpha_1$, $\cdots$, $\alpha_\mu$ and $\beta_1$, $\cdots$, $\beta_\nu$ are roots of $P(z,f)$ and $Q(z,f)$, respectively. For each $\alpha_i$, if $\alpha_i$ is not a Picard exceptional rational function of $f$, then we let $z_0\in\mathbb{C}$ be such that $f(z_0)-\alpha_i(z_0)=0$ with multiplicity $m\in\mathbb{Z}^{+}$. We write $k_i=nk_{i1}+k_{i2}$, where $k_{i1},k_{i2}\in \mathbb{N}$ and $k_{i2}<n$. Note that $(n,k_{i2})=(n,k_i)$. Now, $n|mk_{i2}$ with at most finitely many exceptions since otherwise $z_0+1$ would be a branch point of $f$. Hence, $n\leq mk_{i2}$ and so $m \geq n/(n,k_{i2})$, i.e., $m \geq n/(n,k_i)$. Therefore, we have
    \begin{equation*}
    \overline{N}(r,\alpha_i,f) \leq \frac{(n,k_i)}{n} N(r,\alpha_i,f)+O(\log r).
    \end{equation*}
In particular, we have $m\geq 2$ since $n/(n,k_i)>1$. Thus $\alpha_i$ is a completely ramified rational function of $f$ with multiplicity at least $n/(n,k_i)$. The same analysis above applies for each $\beta_j$ by writing $l_j=nl_{j1}+l_{j2}$, where $l_{j1},l_{j2}\in \mathbb{N}$ and $l_{j2}<n$, as well. Therefore, if $\beta_j$ is not a Picard exceptional rational function of $f$, then we also have
    \begin{equation*}
    \overline{N}(r,\beta_j,f) \leq \frac{(n,l_j)}{n} N(r,\beta_j,f)+O(\log r),
    \end{equation*}
and thus $\beta_j$ is a completely ramified rational function of $f$ with multiplicity at least $n/(n,l_j)$. Below we consider the two cases $q=0$ and $q\geq 1$, respectively.

When $q=0$, equation \eqref{yanagiharaeq0} takes the following form:
\begin{equation}\label{ljkjh}
\overline{f}^n=P(z,f).
\end{equation}
We claim that $f$ has at most finitely many poles. Otherwise, let $z_0\in\mathbb{C}$ be a pole of $f$ with multiplicity $m\in\mathbb{Z}^{+}$. We may choose $z_0$ such that $|z_0|$ is large enough so that none of the coefficients of $P(z,f)$ has poles or zeros outside of $\{z\in \mathbb{C}: |z|<|z_0|\}$. When $n>p$, from \eqref{ljkjh} we see that $z_0+1$ is a pole of $f$ of order $pm/n$ and by iterating along the pole sequence we have $z_0+s$ is a pole of $f$ of order $p^sm/n^s$, $s\in \mathbb{N}$. By letting $s\to \infty$, it follows that there is necessarily a branch point of $f$ at some $z_0+s_0$, $s_0\in \mathbb{N}$, a contradiction to our assumption that $f$ is meromorphic. On the other hand, when $n<p$, from \eqref{ljkjh} we see that $z_0-s$ is a pole of $f$ of order $n^sm/p^s$, $s\in \mathbb{N}$, and by letting $s\to \infty$ we still get the same contradiction as above. Therefore, $f$ has at most finitely many poles, i.e., $\infty$ is a Picard exceptional rational function of $f$. Then the inequality \eqref{multiplicityinequality} implies that $N_c\leq 2$. In particular, when $n\geq 3$, $N_c=2$ is impossible; otherwise, $\alpha_i$ would have multiplicity at least~$n$ for at least one $i$ under our assumptions, a contradiction to the inequality \eqref{multiplicityinequality}. Thus we have the assertions for the case $q=0$.


When $q\geq 1$, since $\alpha_i$, as well as $\beta_j$, is either a Picard exceptional rational function of $f$ or a completely ramified rational function of $f$, then by Picard's theorem and Theorem~\ref{completelyrm} we conclude that $N_c\leq 4$. In particular, when $n\geq 3$, since $n/(n,k_i)>2$ or $n/(n,l_j)>2$ for at least one index $i$ or $j$ by our assumption, then such $\alpha_i$ or $\beta_j$ is either a Picard exceptional rational function of $f$ or a completely ramified rational function with multiplicity at least~$3$ and so by the inequality \eqref{multiplicityinequality} it follows that $N_c\leq 3$; when $n\nmid|p-q|$, by applying the above analysis to poles of $f$, it follows that $\infty$ is either a Picard exceptional rational function of $f$ or a completely ramified rational function of $f$, then by the inequality \eqref{multiplicityinequality} we also have $N_c\leq 3$. Thus the assertions of the lemma for the case $q\geq 1$ follow and this also completes the proof.
\end{proof}

By giving a more careful analysis on the roots of $P(z,f)$ and $Q(z,f)$, we have the four following Lemmas~\ref{basiclemma0}--\ref{keylemma1}, which play key roles in reducing equation \eqref{yanagiharaeq0} into certain forms in the following sections.

\begin{lemma}\label{basiclemma0}
Let $f$ be a transcendental meromorphic solution of equation \eqref{yanagiharaeq0}. Suppose that some $\alpha_i$ in \eqref{P} is~$0$. Then $0$ is a Picard exceptional rational function of $f$. Moreover, if $q=0$, or $q\geq 1$ and $n\nmid |p-q|$, then $\infty$ is also a Picard exceptional rational function of $f$.
\end{lemma}

\renewcommand{\proofname}{Proof.}
\begin{proof}
Suppose that $f$ has infinitely many zeros. Let $z_0\in \mathbb{C}$ be a zero of $f$ with multiplicity $m\in\mathbb{Z}^{+}$. We may choose $z_0$ such that $|z_0|$ is large enough so that none of the coefficients of $P(z,f)$ and $Q(z,f)$ has poles or zeros outside of $\{z\in \mathbb{C}: |z|<|z_0|\}$. Since some $\alpha_i$ is zero, from \eqref{yanagiharaeq0} we see that $z_0+1$ is a zero of $f$ of order $k_im/n$ and by iterating along the zero sequence we have $z_0+s$ is a zero of $f$ of order $k_i^sm/n^s$, $s\in \mathbb{N}$. By letting $s\to \infty$, since $n\nmid k_i$, it follows that there is necessarily a branch point of $f$ at some $z_0+s_0$, $s_0\in \mathbb{N}$, a contradiction to our assumption that $f$ is meromorphic. Therefore, $f$ has at most finitely many zeros, i.e., $0$ is a Picard exceptional rational function of $f$.

Moreover, if $q=0$, or if $q\geq 1$ and $p<q$, then it follows immediately from the proof of Lemma~\ref{basiclemma} that $\infty$ is a Picard exceptional rational function of $f$. Consider the case when $p>q\geq 1$ and $n\nmid (p-q)$. Suppose that $f$ has infinitely many poles. Let $z_0\in\mathbb{C}$ be a pole of $f$ with multiplicity $m\in\mathbb{Z}^{+}$. We may choose $z_0$ such that $|z_0|$ is large enough so that none of the coefficients of $P(z,f)$ and $Q(z,f)$ has poles or zeros outside of $\{z\in \mathbb{C}: |z|<|z_0|\}$. From equation \eqref{yanagiharaeq0} we see that $f$ has a pole of order $(p-q)m/n$ at $z=z_0+1$ and by iterating along the pole sequence $z=z_0+s$ is a pole of $f$ of order $(p-q)^{s}m/n^{s}$, $s\in\mathbb{N}$. Since $n\nmid(p-q)$, then by letting $s\to\infty$, it follows that there is necessarily a branch point of $f$ at $z_0+s_0$ for some $s_0\in\mathbb{N}$, a contradiction to our assumption that $f$ is meromorphic. Therefore, when $p>q\geq 1$ and $n\nmid (p-q)$, $f$ has at most finitely many poles, i.e., $\infty$ is a Picard exceptional rational function of $f$. We complete the proof.

\end{proof}

\begin{lemma}\label{keylemma}
Let $f$ be a transcendental meromorphic solution of equation \eqref{yanagiharaeq0} and $\gamma\in \mathcal{R}\setminus\{0\}$ be a rational function. Then $\gamma$ cannot be a Picard exceptional rational function of $f$. Moreover, if $\gamma$ is a completely ramified function of $f$ with multiplicity at least~$m$, then $\omega\gamma$ is a completely ramified function of $f$ with multiplicity at least~$m$, where $\omega$ is the $n$-th root of~1.
\end{lemma}

\renewcommand{\proofname}{Proof.}
\begin{proof}
To prove the assertions of the lemma, we divide equation \eqref{yanagiharaeq0} into the following three cases:
\begin{enumerate}
\item [(1)] at least one of $\alpha_i$ and $\beta_j$ in \eqref{P} and \eqref{Q} is non-zero and $q=0$, or $q\geq 1$ and $n\nmid |p-q|$;
\item [(2)] at least one of $\alpha_i$ and $\beta_j$ in \eqref{P} and \eqref{Q} is non-zero and $q\geq 1$ and $n\mid |p-q|$;
\item [(3)] there is only one $\alpha_i$ or $\beta_j$ in \eqref{P} and \eqref{Q} and this $\alpha_i$ or $\beta_j$ is zero.
\end{enumerate}

We first suppose that $\gamma$ is a Picard exceptional rational function of $f$. Under this assumption, below we show that each of the above three cases will lead to contradictions.

In the first case, we let $\beta$ be such that $\beta=\alpha_i$ or $\beta=\beta_j$ for some $\alpha_i$ or $\beta_j$ in \eqref{P} and \eqref{Q} and $\beta\not=0$. Denote the order of this root $\alpha_i$ or $\beta_j$ by $t_1$. Put
\begin{equation}\label{keylemma equa001}
u=\frac{\overline{f}}{f-\beta},  \quad  v=\frac{1}{f-\beta}.
\end{equation}
Then $u$ and $v$ are two algebroid functions with at most finitely many branch points and we have
\begin{equation*}
\overline{f}=\frac{u}{v},  \quad  f=\frac{1}{v}+\beta,
\end{equation*}
and it follows that \eqref{yanagiharaeq0} becomes
\begin{equation}\label{keylemma equa003}
u^n=\frac{P_1(z,v)}{Q_1(z,v)}v^{n_1},
\end{equation}
where $n_1\in \mathbb{Z}$, $P_1(z,v)$ and $Q_1(z,v)$ are two polynomials in $v$ having no common factors and none of the roots of $P_1(z,v)$ or $Q_1(z,v)$ is zero. Denote by $p_1=\deg_{v}(P_1(z,v))$ the degree of $P_1(z,v)$ in $v$ and by $q_1=\deg_{v}(Q_1(z,v))$ the degree of $Q_1(z,v)$ in $v$, respectively. Note that $q=0$, or $q\geq 1$ and $n\nmid |p-q|$. By elementary calculations, when $q=0$ we get $n_1=n-p$, $p_1=p-t_1$ and $q_1=0$; when $q\geq 1$ and $\beta=\alpha_i$ we get $n_1=n+q-p$, $p_1=p-t_1$ and $q_1=q$; when $q\geq 1$ and $\beta=\beta_j$ we get $n_1=n+q-p$, $p_1=p$ and $q_1=q-t_1$. Therefore, we always have $n_1\not=0$, $n_1\not=n$ and $p_1-q_1+n_1\not=n$. We consider
\begin{equation}\label{keylemma equa004}
\overline{f}-\overline{\gamma}=\frac{u}{v}-\overline{\gamma}=\frac{u-\overline{\gamma}v}{v}.
\end{equation}
Let $(u_0,v_0)$ be any pair of non-zero functions satisfying the following system of equations:
\begin{equation}\label{keylemma equa005}
u_0^n=\frac{P_1(z,v_0)}{Q_1(z,v_0)}v_0^{n_1}, \quad  u_0-\overline{\gamma}v_0=0.
\end{equation}
By assumption, the equation $\overline{f}-\overline{\gamma}=0$ has at most finitely many roots. Let $z_0\in \mathbb{C}$ be such that $u(z_0)=u_0(z_0)$, $v(z_0)=v_0(z_0)$ and $u_0=\overline{\gamma}v_0$ hold simultaneously. Then from \eqref{keylemma equa004} we see that the equation $u-\overline{\gamma}v=0$ can have at most finitely many such roots $z_0$. Now we consider the equation $\overline{f}-\omega\overline{\gamma}=0$, where $\omega$ is the $n$-th root of $1$. We have
\begin{equation}\label{keylemma equa006}
\overline{f}-\omega\overline{\gamma}=\frac{u}{v}-\omega\overline{\gamma}=\frac{u-\omega\overline{\gamma}v}{v}.
\end{equation}
Let $\hat{z}_0\in \mathbb{C}$ be such that $v(\hat{z}_0)=0$ or $v(\hat{z}_0)=\infty$. Since $n_1\not=0$, $n_1\not=n$ and $p_1-q_1+n_1\not=n$, we see from \eqref{keylemma equa003} that the multiplicity of $\hat{z}_0$ for $v(\hat{z}_0)=0$ or $v(\hat{z}_0)=\infty$ equals the multiplicity of $\hat{z}_0$ for $u(\hat{z}_0)=0$ or $u(\hat{z}_0)=\infty$ for at most finitely many such $\hat{z}_0$. Then it follows from \eqref{keylemma equa006} that the two equations $v=0$ and $\overline{f}-\omega\overline{\gamma}=0$, as well as the two equations $v=\infty$ and $\overline{f}-\omega\overline{\gamma}=0$, cannot have infinitely many common roots; otherwise, $\overline{\gamma}$ would be identically equal to $0$ or $\infty$, a contradiction to our assumption. Let $(u_1,v_1)$ be any pair of non-zero functions satisfying the following system of equations:
\begin{equation}\label{keylemma equa007}
u_1^n=\frac{P_1(z,v_1)}{Q_1(z,v_1)}v_1^{n_1}, \quad  u_1-\omega\overline{\gamma}v_1=0.
\end{equation}
Recall that $\omega$ is the $n$-th root of~$1$. Then the two systems of equations in \eqref{keylemma equa005} and \eqref{keylemma equa007} yield identically the same algebraic equations for $v_0$ and $v_1$. Since $u_0^n$ and $u_1^n$ equals the same rational term in $v_0$ or $v_1$, we see that the two systems of equations in \eqref{keylemma equa005} and \eqref{keylemma equa007} have the same pair of non-zero solutions, apart from permutations. By summarizing the above results, we conclude that the equation $\overline{f}-\omega\overline{\gamma}=0$ can have at most finitely many roots, i.e., $\omega\gamma$ is a Picard exceptional rational function of~$f$.

Now, since $n\geq 2$, we must have $n=2$; otherwise, $f$ would have three or more distinct Picard exceptional rational functions, a contradiction to Picard's theorem. However, when $q=0$, from Lemma~\ref{basiclemma} we know $\infty$ is also a Picard exceptional rational function of $f$, a contradiction to Picard's theorem; when $q\geq 1$, since $n\nmid |p-q|$, from Lemma~\ref{basiclemma} we know $\infty$ is either a Picard exceptional rational function of $f$ or a completely ramified rational function of $f$, a contradiction to Picard's theorem or the inequality \eqref{multiplicityinequality}.

In the second case, we must have $N_c\geq 2$ since $n\mid |p-q|$. We let $\beta$ be such that $\beta=\alpha_i$ or $\beta=\beta_j$ for some $\alpha_i$ or $\beta_j$ in \eqref{P} and \eqref{Q} and $\beta\not=0$, and $\alpha$ be such that $\alpha=\alpha_i$ or $\alpha=\beta_j$ for another $\alpha_i$ or $\beta_j$ in \eqref{P} and \eqref{Q} distinct from $\beta$. Put
\begin{equation}\label{keylemma equa008}
u=\frac{\overline{f}(f-\alpha)}{f-\beta},  \quad  v=\frac{f-\alpha}{f-\beta}.
\end{equation}
Then $u$ and $v$ are two algebroid functions with at most finitely many branch points and we have
\begin{equation*}\label{keylemma equa009}
\overline{f}=\frac{u}{v},  \quad  f=\frac{\beta v-\alpha}{v-1},
\end{equation*}
and it follows that \eqref{yanagiharaeq0} becomes
\begin{equation*}
u^n=\frac{P_2(z,v)}{Q_2(z,v)}v^{n_2},
\end{equation*}
where $n_2\in \mathbb{Z}$, $P_2(z,v)$ and $Q_2(z,v)$ are two polynomials in $v$ having no common factors and none of the roots of $P_2(z,v)$ or $Q_2(z,v)$ is zero. Denote by $p_2=\deg_{v}(P_2(z,v))$ the degree of $P_2(z,v)$ in $v$ and by $q_2=\deg_{v}(Q_2(z,v))$ the degree of $Q_2(z,v)$ in $v$, respectively. Note that $q\geq 1$ and $n\mid |p-q|$. As in the previous case we can show that $n_2\not=0$, $n_2\not=n$ and $p_2-q_2+n_2\not=n$ by elementary calculations. Then we consider the roots of the equation $\overline{f}-\overline{\gamma}=0$ and of the equation $\overline{f}-\omega\overline{\gamma}=0$, respectively, and by the same arguments as in the previous case we get the same conclusion as there.

We have $n=2$ in both cases~(1) and~(2). If some $\alpha_i$ or $\beta_j$ in \eqref{P} and \eqref{Q} is distinct from $\pm \gamma$, then from Lemma~\ref{basiclemma} it follows that this $\alpha_i$ or $\beta_j$ is either a Picard exceptional rational function of $f$ or a completely ramified rational function of $f$, a contradiction to Picard's theorem or the inequality \eqref{multiplicityinequality} since $\pm \gamma$ are both Picard exceptional rational functions of $f$. Therefore, $N_c=2$ and the two roots $\alpha_i$ or $\beta_j$ in \eqref{P} and \eqref{Q} are equal to $\pm \gamma$. We consider
\begin{equation}\label{keylemma equa01}
\overline{f}^2-\overline{\gamma}^2=\frac{P(z,f)-\overline{\gamma}^2Q(z,f)}{Q(z,f)}=\frac{a_{p_{\tau}}(f-\gamma_1)^{t_1}\cdots (f-\gamma_{\tau})^{t_{\tau}}}{Q(z,f)},
\end{equation}
where $\gamma_1$, $\cdots$, $\gamma_{\tau}$ are in general algebraic functions distinct from each other, $t_1,\cdots,t_{\tau}\in \mathbb{N}$ denote the orders of the
roots $\gamma_1$, $\cdots$, $\gamma_{\tau}$, respectively, and $t_1+\cdots+t_{\tau}=p_{\tau}\in \mathbb{N}$. Note that the equation $\overline{f}^2-\overline{\gamma}^2=0$ can have at most finitely many roots. If $p_{\tau}<q$, then we obtain from \eqref{keylemma equa01} that $\infty$ is a Picard exceptional rational function of $f$, a contradiction to Picard's theorem. Therefore, $p_{\tau}\geq q$. From the previous discussions we see that none of $\gamma_1$, $\cdots$, $\gamma_{\tau}$ is equal to $\pm \gamma$. But then we have $\gamma_1$, $\cdots$, $\gamma_{\tau}$ are all Picard exceptional rational functions of $f$ by analyzing the roots of the equations $f-\gamma_{i}=0$, $i=1,\cdots, \tau$, again a contradiction to Picard's theorem since $p_{\tau}\geq 1$.

In the third case, it follows from the assumptions that $n\nmid|p-q|$. Moreover, we have $q\geq 1$ and this $\alpha_i$ or $\beta_j$ is $\beta_1$ and $\beta_1\equiv0$; otherwise, this $\alpha_i$ or $\beta_j$ is $\alpha_1$ and $\alpha_1=0$, but it follows by Lemma~\ref{basiclemma0} that $0$ and $\infty$ are also both Picard exceptional rational functions of $f$, a contradiction to Picard's theorem. Put
\begin{equation}\label{keylemma equa0010a}
u=\frac{\overline{f}}{f},  \quad  v=\frac{1}{f}.
\end{equation}
Then $u$ and $v$ are two algebroid functions with at most finitely branch points and we have
\begin{equation*}
\overline{f}=\frac{u}{v},  \quad  f=\frac{1}{v},
\end{equation*}
and it follows that \eqref{yanagiharaeq0} becomes
\begin{equation*}
u^n=\frac{P_3(z,v)}{Q_3(z,v)}v^{n_3},
\end{equation*}
where $n_3\in \mathbb{Z}$, $P_3(z,v)$ and $Q_3(z,v)$ are two polynomials in $v$ having no common factors and none of the roots of $P_3(z,v)$ or $Q_3(z,v)$ is zero. Denote by $p_3=\deg_{v}(P_3(z,v))$ the degree of $P_3(z,v)$ in $v$ and by $q_3=\deg_{v}(Q_3(z,v))$ the degree of $Q_3(z,v)$ in $v$, respectively. As in case~{(1)} we can show that $n_3\not=0$, $n_3\not=n$ and $p_3-q_3+n_3\not=n$ by elementary calculations.
Then we consider the roots of the equation $\overline{f}-\overline{\gamma}=0$ and of the equation $\overline{f}-\omega\overline{\gamma}=0$, respectively, and by the same arguments as in case~{(1)} we get that $\omega\gamma$ is a Picard exceptional rational function of $f$. But this is impossible by the inequality \eqref{multiplicityinequality} since from Lemma~\ref{basiclemma} it follows that $\infty$ is either a Picard exceptional rational function of $f$ or a completely ramified rational function of $f$.

From the above reasoning, we conclude that $\gamma$ cannot be a Picard exceptional rational function of $f$. This gives our first assertion of the lemma.

Next, we suppose that $\gamma$ is a completely ramified rational function of $f$ with multiplicity at least~$m$. We also consider the three cases in the beginning of the proof. In the first case, we do the transformations in \eqref{keylemma equa001} and get the equation in \eqref{keylemma equa003}. By assumption, all roots of the equation $\overline{f}-\overline{\gamma}=0$ with at most finitely many exceptions have multiplicities at least~$m$.
Therefore, for any pair of non-zero functions $(u_0,v_0)$ such that the system of equations in \eqref{keylemma equa005} holds, if we let $z_0\in \mathbb{C}$ be such that $u(z_0)=u_0(z_0)$, $v(z_0)=v_0(z_0)$ and $u_0=\overline{\gamma}v_0$ hold simultaneously, then from \eqref{keylemma equa004} we see that the equation $u-\overline{\gamma}v=0$ has at most finitely many such roots $z_0$ with multiplicity less than~$m$. Moreover, letting $\hat{z}_0\in \mathbb{C}$ be such that $f(\hat{z}_0+1)-\omega\gamma(\hat{z}_0+1)=0$, from previous discussions we know that $v(\hat{z}_0)=0$ or $v(\hat{z}_0)=\infty$ for at most finitely many such $\hat{z}_0$. Also, we know that the two systems of equations in \eqref{keylemma equa005} and \eqref{keylemma equa007} have the same pair of solutions, apart from permutations. Then we conclude from the above reasoning that the equation $\overline{f}-\omega\overline{\gamma}=0$ can have at most finitely many roots with multiplicities less than~$m$. The second and the third cases can be discussed in an analogous way as above after doing the transformations in \eqref{keylemma equa008} and \eqref{keylemma equa0010a}, respectively. We omit those details. Thus we have the second assertion of the lemma and also the complete proof.

\end{proof}

\begin{lemma}\label{keylemma0}
Let $f$ be a transcendental meromorphic solution of equation \eqref{yanagiharaeq0}. Suppose that one of the following cases occurs:
\begin{enumerate}
\item [(1)] $q=0$ and $P(z,f)$ has a root $\alpha$ of order $k$ such that $2\leq k<n$ and $k\nmid n$, or $k\geq n+1$;
\item [(2)] $q\geq 1$ and $P(z,f)$ has a root $\alpha$  of order $k$ such that $2\leq k<n$ and $k\nmid n$, or $k\geq n+1$;
\item [(3)] $q\geq 1$ and $Q(z,f)$ has a root $\beta$  of order $l$ such that $2\leq l<n$ and $l\nmid n$, or $l\geq n+1$;
\item [(4)] $q\geq 1$ and $p$ and $q$ satisfy $2\leq |p-q|<n$ and $|p-q|\nmid n$, or $|p-q|\geq n+1$.
\end{enumerate}
Then if $q=0$, then $f$ cannot have~2 completely ramified rational functions in $\mathcal{R}$; if $q\geq1$, then $f$ cannot have~4 completely ramified rational functions in $\mathcal{R}$ and cannot have~3 non-zero completely ramified rational functions $\gamma_1,\gamma_2,\gamma_3\in \mathcal{R}$ such that $\sum_{i=1}^3\Theta(\gamma_i,f)=2$.

\end{lemma}

\renewcommand{\proofname}{Proof.}
\begin{proof}

First, when $q=0$, by Lemma~\ref{basiclemma} we know that $\infty$ is a Picard exceptional rational function of $f$ and thus $\Theta(\infty,f)=1$. Suppose that $\gamma_1$ and $\gamma_2$ are both completely ramified rational functions of $f$. By the inequality~\eqref{multiplicityinequality} it follows that $\gamma_1$ and $\gamma_2$ both have multiplicities~2 and that $\Theta(\gamma_1,f)+\Theta(\gamma_2,f)=1$. Further, by the Second Main Theorem of Yamanoi for rational functions as targets \cite{yamanoi:05}, letting $\gamma$ be any rational function distinct from $\gamma_i$, we must have $\overline{N}(r,\gamma,f)=T(r,f)+o(T(r,f))$, $N(r,\gamma_i,f)=T(r,f)+o(T(r,f))$ and $\overline{N}(r,\gamma_i,f)=\frac{1}{2}T(r,f)+o(T(r,f))$, where $r\to\infty$ outside an exceptional set $E$ with finite linear measure. If either~$\gamma_1=0$ or $\gamma_2=0$, say $\gamma_1=0$, then by Lemma~\ref{keylemma} it follows that $\omega \gamma_2$ is a completely ramified rational function of $f$, where $\omega$ is the $n$-th root of~1, a contradiction to the inequality~\eqref{multiplicityinequality}. Therefore,~$0$ is not a completely ramified rational function of $f$. Let $z_0\in\mathbb{C}$ be such that $f(z_0)-\alpha(z_0)=0$ with multiplicity $m\in\mathbb{Z}^{+}$. Since $\alpha$ is a root of $P(z,f)$ of order $k$ such that $2\leq k<n$ and $k\nmid n$, or $k\geq n+1$, then from \eqref{yanagiharaeq0} we see that $z_0$ is a root of $f(z+1)=0$ with multiplicity $m_0\in\mathbb{Z}^{+}$ such that $m_0=mk/n\geq 2$ with at most finitely many exceptions. By the Valiron--Mohon'ko identity \cite{Valiron1931,mohonko1971} (see also \cite{Laine1993}), we have from \eqref{yanagiharaeq0} that $nT(r,\overline{f})=dT(r,f)+O(\log r)$.
Now there are at least $T(r,f)+o(T(r,f))$ many points $z_0$ such that $f(z_0+1)=0$ with multiplicity~$m_0\geq 2$. This implies that $n<d$. Denote by $S$ the set of zeros of $\overline{f}$ in the finite disk $D=\{z\in \mathbb{C}: |z|< t\}$, where $t>0$, and by $S_1$ the set of zeros with multiplicity $\geq 2$ of $\overline{f}$ in $D$. Denote by $n_{S_1}(t,1/\overline{f})$ and $\overline{n}_{S_1}(t,1/\overline{f})$ the number of zeros of $\overline{f}$ in $S_1$, counting or ignoring multiplicities, respectively; denote by $n_{S\setminus S_1}(t,1/\overline{f})$ and $\overline{n}_{S\setminus S_1}(t,1/\overline{f})$ the number of zeros of $\overline{f}$ for the complement of $S_1$, counting or ignoring multiplicities, respectively. We may suppose that $f(1)\not=0$. By the definition of the truncated counting function $\overline{N}(r,1/\overline{f})$, we deduce that
     \begin{equation*}
     \begin{split}
      \overline{N}(r,1/\overline{f})&=\int_{0}^{r}\overline{n}_{S\setminus S_1}(t,1/\overline{f})\frac{dt}{t}+\int_{0}^{r}\overline{n}_{S_1}(t,1/\overline{f})\frac{dt}{t}\\
      &\leq \int_{0}^{r}n_{S\setminus S_1}(t,1/\overline{f})\frac{dt}{t}+\frac{1}{2}\int_{0}^{r}n_{S_1}(t,1/\overline{f})\frac{dt}{t}\\
      &=\int_{0}^{r}n_{S}(t,1/\overline{f})\frac{dt}{t}-\frac{1}{2}\int_{0}^{r}n_{S_1}(t,1/\overline{f})\frac{dt}{t}\\
      &\leq T(r,1/\overline{f})-\frac{1}{2}[T(r,f)+o(T(r,f))].
      \end{split}
     \end{equation*}
By the First Main Theorem of Nevanlinna we have $T(r,1/\overline{f})=T(r,\overline{f})+O(1)$ and it follows that $T(r,1/\overline{f})=\frac{d}{n}T(r,f)+O(\log r)$. Then by combining the above results together we get
     \begin{equation*}
     \Theta(0,\overline{f})=1-\limsup_{r\to\infty}\frac{\overline{N}(r,1/\overline{f})}{T(r,\overline{f})}\geq \frac{n}{2d}.
     \end{equation*}
In general, the quantity $\Theta(\gamma_i,f)$ may not be shift-invariant \cite{Korhonenkazhzh2020}, but under our assumptions we already have $\sum_{i=1}^2\Theta(\overline{\gamma}_i,\overline{f})=1$, and thus the inequality above is impossible. Therefore, if case~{(1)} occurs, then $f$ cannot have~2 completely ramified rational functions in $\mathcal{R}$. This is the first assertion of the lemma.

Second, we consider the case where $q\geq 1$. We suppose that $f$ has four completely ramified rational functions $\gamma_i$, $i=1,2,3,4$, in $\mathcal{R}$. By Theorem~\ref{completelyrm}, $\infty$ is not a completely ramified rational function of $f$. Moreover, none of $\gamma_1$, $\gamma_2$, $\gamma_3$ and $\gamma_4$ is zero; otherwise, say $\gamma_1=0$, by Lemma~\ref{keylemma} it follows that $\omega\gamma_2$, $\omega\gamma_3$ and $\omega\gamma_4$ are all completely ramified rational functions of $f$, where $\omega$ is the cubic root of~1, and thus by Theorem~\ref{completelyrm} we must have $n\geq 3$. However, by Lemma~\ref{basiclemma} it follows that at least one of $\alpha_i$ in \eqref{P} and $\beta_j$ in \eqref{Q} is a completely ramified rational function of $f$ with multiplicity~$\geq 3$, a contradiction to the inequality~\eqref{multiplicityinequality}. Therefore, $0$ is not a completely ramified rational function of $f$. Below we consider the three cases~{(2)}, {(3)} and {(4)}, respectively.

If case~{(2)} occurs, then $\alpha$ is not a Picard exceptional rational function of $f$. As in case~{(1)}, for the point $z_0\in\mathbb{C}$ such that $f(z_0)-\alpha(z_0)=0$ with multiplicity $m\in\mathbb{Z}^{+}$, we have that $z_0$ is a root of $f(z+1)=0$ with multiplicity $m_0\in\mathbb{Z}^{+}$ such that $m_0=mk/n\geq 2$ with at most finitely many exceptions. Since $f$ has four completely ramified rational functions $\gamma_i$, $i=1,2,3,4$, then by the inequality \eqref{multiplicityinequality} we know that $\gamma_i$, $i=1,2,3,4$, all have multiplicities~2. Further, by the Second Main Theorem of Yamanoi for rational functions as targets \cite{yamanoi:05}, letting $\gamma$ be any rational function distinct from $\gamma_i$, we must have $\overline{N}(r,\gamma,f)=T(r,f)+o(T(r,f))$, $N(r,\gamma_i,f)=T(r,f)+o(T(r,f))$ and $\overline{N}(r,\gamma_i,f)=\frac{1}{2}T(r,f)+o(T(r,f))$, where $r\to\infty$ outside an exceptional set $E$ with finite linear measure. Then, similarly as in case~{(1)}, we can obtain a contradiction by computing $\Theta(0,\overline{f})$. Therefore, if case~{(2)} occurs, then $f$ cannot have~4 completely ramified rational functions.

If case~{(3)} occurs, then for the point $z_0\in\mathbb{C}$ such that $f(z_0)-\beta(z_0)=0$ we have that $z_0$ is a root of $f(z+1)=\infty$ with multiplicity $m_0\in\mathbb{Z}^{+}$ such that $m_0=ml/n\geq 2$ with at most finitely many exceptions, and then we can obtain a contradiction by computing $\Theta(\infty,\overline{f})$. Therefore, we still have the same conclusion as in case~{(2)}.

If case~{(4)} occurs, then we let $z_0$ be a pole of $f$ with multiplicity~$m\in\mathbb{Z}^{+}$ and it follows that $z_0$ is a root of $f(z+1)=0$ or $f(z+1)=\infty$ with multiplicity $m_0\in\mathbb{Z}^{+}$ such that $m_0=m|p-q|/n\geq 2$ with at most finitely many exceptions. If $p<q$, then we get a contradiction by computing $\Theta(0,f)$; if $p>q$, then we get a contradiction by computing $\Theta(\infty,f)$. Therefore, we still have the same conclusion as in case~{(2)}.

Last, we suppose that $f$ has~3 non-zero completely ramified rational functions $\gamma_1,\gamma_2,\gamma_3\in \mathcal{R}$ such that $\sum_{i=1}^3\Theta(\gamma_i,f)=2$. From \cite[p.~46]{Hayman1964Meromorphic} we know that the possible multiplicity sets $(m_1,m_2,m_3)$ corresponding to $\gamma_1$, $\gamma_2$, $\gamma_3$ are $(2,4,4)$, $(2,3,6)$ or $(3,3,3)$, apart from permutations. Also, by the Second Main Theorem of Yamanoi for rational functions as targets \cite{yamanoi:05}, we have $\overline{N}(r,\gamma_i,f)=\frac{1}{m_i}T(r,f)+o(T(r,f))$, where $r\to\infty$ outside an exceptional set $E$ with finite linear measure. Note that $\infty$ is not a completely ramified rational function of $f$. If one of the three cases~{(2)}, {(3)} or {(4)} occurs, then we can use the same arguments as above to compute $\Theta(0,\overline{f})$ or $\Theta(\infty,\overline{f})$ and obtain similar contradictions. Thus our assertion follows.

\end{proof}

In Lemma~\ref{keylemma0}, we did not deal with equation~\eqref{yanagiharaeq0} for the case when $f$ has three non-zero completely ramified functions $\gamma_i$ such that $\sum_{i=1}^3\Theta(\gamma_i,f)=2$ and one $\gamma_i$ is $\infty$. We will exclude out this possibility in the proof of Theorem~\ref{Theorem  2b} in section~\ref{Proof 2} with applications of the analysis in the proof of Lemma~\ref{keylemma0}.

Finally, we consider the equation $\overline{f}^n-\overline{\gamma}^n=0$ further, where $\gamma\in \mathcal{R}\setminus\{0\}$ is a completely ramified rational function of $f$ with multiplicity at least~$m\geq 2$. By Lemma~\ref{keylemma}, $\omega\gamma$ is a completely ramified rational function of $f$ with multiplicity at least~$m$, where $\omega$ is the $n$-th root of~1. By \eqref{yanagiharaeq0}, when $q=0$ we have
\begin{equation}\label{keylemma equa0}
\overline{f}^n-\overline{\gamma}^n=P(z,f)-\overline{\gamma}^n=a_p(f-\gamma_1)^{t_1}\cdots (f-\gamma_{\tau})^{t_{\tau}},
\end{equation}
or, when $q\geq 1$, we have
\begin{equation}\label{keylemma equa1}
\overline{f}^n-\overline{\gamma}^n=\frac{P(z,f)-\overline{\gamma}^nQ(z,f)}{Q(z,f)}=\frac{a_{p_{\tau}}(f-\gamma_1)^{t_1}\cdots (f-\gamma_{\tau})^{t_{\tau}}}{Q(z,f)},
\end{equation}
where $\gamma_1$, $\cdots$, $\gamma_{\tau}$ are in general algebraic functions distinct from each other and $t_1,\cdots,t_{\tau}\in \mathbb{N}$ denote the orders of the roots $\gamma_1$, $\cdots$, $\gamma_{\tau}$, respectively, and $t_1+\cdots+t_{\tau}=p_{\tau}\in \mathbb{N}$. We apply the analysis in the proof of Lemma~\ref{keylemma0} to equations \eqref{keylemma equa0} and \eqref{keylemma equa1}, respectively, and get the following

\begin{lemma}\label{keylemma1}
Let $f$ be a transcendental meromorphic solution of equation \eqref{yanagiharaeq0} and $\gamma\in \mathcal{R}\setminus\{0\}$ be a completely ramified rational function of $f$ with multiplicity~$m\geq 2$. Suppose that $\zeta_i$, $\cdots$, $\zeta_t$ are Picard exceptional rational functions of $f$ or completely ramified rational functions of $f$ such that $\sum_{i=1}^t\Theta(\zeta_i,f)=2$. For each $\gamma_i$ in \eqref{keylemma equa0} or \eqref{keylemma equa1}, if $\gamma_i$ is not a completely ramified rational function of $f$, then $t_i=m$; if $\gamma_i$ is a completely ramified rational function of $f$ with multiplicity~$m_i\geq 2$, then $t_im_i=m$. In particular, for \eqref{keylemma equa1}, when $1\leq p_{\tau}<q$, if $\infty$ is not a completely ramified rational function of $f$, then $q-p_{\tau}=m$; if $\infty$ is a completely ramified rational function of $f$ with multiplicity $m_{\infty}\geq 2$, then $(q-p_{\tau})m_{\infty}=m$.
\end{lemma}

\renewcommand{\proofname}{Proof.}
\begin{proof}
By the assumption $\sum_{i=1}^t\Theta(\zeta_i,f)=2$, we know from the proof of Lemma~\ref{keylemma0} that $t=3$ or $t=4$. Moreover, for each $\gamma_i$ in \eqref{keylemma equa0} or \eqref{keylemma equa1} we have $\overline{N}(r,\gamma_i,f)=T(r,f)+o(T(r,f))$ when $\gamma_i$ is not a completely ramified rational function of $f$ and $\overline{N}(r,\gamma_i,f)=\frac{1}{m_i}T(r,f)+o(T(r,f))$ when $\gamma_i$ is a completely ramified rational function of $f$ with multiplicity~$m_i$, where $r\to\infty$ outside an exceptional set $E$ with finite linear measure. In particular, we have $N(r,\infty,f)=O(\log r)$ when $\infty$ is a Picard exceptional rational function of $f$ and, otherwise, we have $N(r,\infty,f)=T(r,f)+o(T(r,f))$ when $\infty$ is not a completely ramified rational function of $f$ and $\overline{N}(r,\infty,f)=\frac{1}{m_{\infty}}T(r,f)+o(T(r,f))$ when $\infty$ is a completely ramified rational function of $f$ with multiplicity~$m_{\infty}$, where again $r\to\infty$ outside an exceptional set $E$ with finite linear measure. Note that for equation \eqref{keylemma equa1}, when $p_{\tau}\geq 1$, from the proof of Lemma~\ref{keylemma} we have that $\infty$ cannot be a Picard exceptional rational function of $f$. We also have $\overline{N}(r,\gamma,f)=\frac{1}{m}T(r,f)+o(T(r,f))$, where $r\to\infty$ outside an exceptional set $E$ with finite linear measure. Now there are $T(r,f)+o(T(r,f))$ many points $z_0\in \mathbb{C}$ such that $f(z_0)-\gamma_i(z_0)=0$ (or $f(z_0)=\infty$ when $p_{\tau}<q$) and from \eqref{keylemma equa0} or \eqref{keylemma equa1} it follows that $f(z_0+1)^n-\gamma(z_0+1)^n=0$. For such $z_0$, by comparing the multiplicities on both sides of the equation \eqref{keylemma equa0} or \eqref{keylemma equa1}, we obtain the desired results. We omit those details.

\end{proof}

\section{Equation \eqref{yanagiharaeq0} with $q=0$}\label{Proof 1}

\subsection{Equation \eqref{yanagiharaeq0} with $q=0$ and $n>p\geq 2$}

For the case $q=0$ and $n>p\geq 2$ of equation \eqref{yanagiharaeq0}, we prove the following

\begin{theorem}\label{Theorem  1a}
Suppose that $q=0$ and $n>p\geq 2$. Let $f$ be a transcendental meromorphic solution of equation \eqref{yanagiharaeq0}. Then there exists a rational function $\alpha$ such that the linear transformation $f\rightarrow \alpha f$ reduces \eqref{yanagiharaeq0} into
\begin{equation}\label{Theorem  1a  eq0}
\overline{f}^n=cf^p,
\end{equation}
where $c$ is a non-zero constant. Moreover, solutions of equation \eqref{Theorem  1a  eq0} are represented as
\begin{equation}\label{Theorem  1a  eq1}
f=c^{\frac{1}{n-p}}\exp[\pi(z)(p/n)^z],
\end{equation}
where $\pi(z)$ is an arbitrary entire periodic function with period~1.
\end{theorem}

\renewcommand{\proofname}{Proof.}
\begin{proof}
From Lemma~\ref{basiclemma} we know that $\infty$ is a Picard exceptional rational function of $f$ and $N_c=1$. Therefore, we have
\begin{equation*}
\overline{f}^n=a_p(f-\alpha_1)^{k_1},
\end{equation*}
where $k_1=p$ and $(n,k_1)=1$. If $\alpha_1\not=0$, then by Lemmas~\ref{basiclemma} and~\ref{keylemma} it follows that $\omega\alpha_1$ is a completely ramified rational function of $f$ with multiplicity at least~$n$, where $\omega$ is the $n$-th root of~1, a contradiction to the inequality \eqref{multiplicityinequality} since $n\geq 3$. Therefore, $\alpha_1=0$ and thus we have
\begin{equation}\label{pre equation  2aaa   mmm  k0}
\overline{f}^n=a_pf^p.
\end{equation}
By Lemma~\ref{basiclemma0}, it follows that $0$ is also a Picard exceptional rational function of $f$. Then there is a non-zero rational function $\alpha$ such that $\alpha f$ is a zero-free entire function. Since $f$ satisfies \eqref{pre equation  2aaa   mmm  k0}, it follows that
\begin{equation*}
(\overline{\alpha}\overline{f})^n=\overline{\alpha}^n\overline{f}^n =\overline{\alpha}^na_pf^p=a_p\frac{\overline{\alpha}^n}{\alpha^p}(\alpha f)^p.
\end{equation*}
By redefining $\alpha f$ as $f$, we have
\begin{equation}\label{yanagiharaeq0    1}
\overline{f}^n=cf^p,
\end{equation}
where $c=a_p\frac{\overline{\alpha}^n}{\alpha^p}$ is a non-zero constant. By taking the logarithm on both sides of \eqref{yanagiharaeq0    1}, then $g=\log f$ is entire and satisfies
\begin{equation*}
n\overline{g}=\log c+pg.
\end{equation*}
Therefore, we can solve $f$ as
\begin{equation*}
f=c^{\frac{1}{n-p}}\exp[\pi(z)(p/n)^z],
\end{equation*}
where $\pi(z)$ is an arbitrary entire periodic function with period~1. This completes the proof.
\end{proof}

\subsection{Equation \eqref{yanagiharaeq0} with $q=0$ and $2\leq n<p$}

For the case $q=0$ and $2\leq n<p$ of equation \eqref{yanagiharaeq0}, we prove the following

\begin{theorem}\label{Theorem  1c}
Suppose that $q=0$ and $2\leq n< p$. Let $f$ be a transcendental meromorphic solution of equation \eqref{yanagiharaeq0}. Then there exists a rational function $\alpha$ such that by doing a linear transformation $f\to \alpha f$, we have either

\begin{enumerate}

\item [(1)] equation \eqref{yanagiharaeq0} reduces into
\begin{equation}\label{yanagiharaeq0    8}
\overline{f}^n=cf^p,
\end{equation}
where $c$ is a non-zero constant; solutions of \eqref{yanagiharaeq0    8} are represented as
\begin{equation}\label{yanagiharaeq0    8be}
f(z)=c^{\frac{1}{n-p}}\exp[\pi(z)(p/n)^z],
\end{equation}
where $\pi(z)$ is an arbitrary entire periodic function with period~1; or

\item [(2)] when $n=2$ and $p=2p_0+1$, $p_0\geq 1$, equation \eqref{yanagiharaeq0} reduces into
\begin{equation}\label{yanagiharaeq0    9}
\overline{f}^2=P_0(f)^2(f-1),
\end{equation}
and $P_0(f)$ is a polynomial in $f$ such that
\begin{equation*}
P_0(f)=\frac{\pm i}{2^{1/2}}[U_{p_0}(f)+U_{p_0-1}(f)]
\end{equation*}
with the Chebyshev polynomials $U_{p_0}$ and $U_{p_0-1}$ of the second kind, i.e.,
\begin{equation*}
U_{p_0}(\cos x)=\frac{\sin(p_0+1)x}{\sin x},
\end{equation*}
i.e.,
\begin{equation*}
U_{p_0}(f)=\sum_{t=0}^{[p_0/2]}\frac{(-1)^t(p_0-t)!}{t!(p_0-2t)!}(2f)^{p_0-2t};
\end{equation*}
therefore, if we write
\begin{equation*}
P_0(f)=\sum_{t=0}^{p_0}A_{p_0}f^{p_0},
\end{equation*}
then
\begin{equation*}
\begin{split}
A_{p_0-2t}&=(-1)^t2^{p_0-2t}\frac{(p_0-t)\cdots(p_0-2t+1)}{t!},\\
A_{p_0-2t-1}&=(-1)^t2^{p_0-2t-1}\frac{(p_0-t-1)\cdots(p_0-2t)}{t!};
\end{split}
\end{equation*}
solutions of \eqref{yanagiharaeq0    9} are represented as
\begin{equation}\label{yanagiharaeq0    21 hj}
f(z)=\frac{1}{2}(\delta^2+\delta^{-2}),
\end{equation}
where $\delta$ is given by
\begin{equation*}
\delta=(\pm i)^{\frac{1}{1-2p_0}}\exp[\pi(z)(p_0+1/2)^z],
\end{equation*}
or
\begin{equation*}
\delta=(\pm i)^{\frac{1}{3+2p_0}}\exp[\pi(z)(-p_0-1/2)^z],
\end{equation*}
where $\pi(z)$ is an arbitrary entire periodic function with period~1; or

\item [(3)] when $n=2$ and $p=2p_0+2$, $p_0\geq 1$, equation \eqref{yanagiharaeq0} reduces into
\begin{equation}\label{yanagiharaeq0    5}
\overline{f}^2=P_0(f)^2(f^2-1),
\end{equation}
and $P_0(f)$ is a polynomial in $f$ such that
\begin{equation}\label{yanagiharaeq0    6}
P_0(f)=\pm i\sum_{l=0}^{[p_0/2]}\binom{p_0+1}{2l+1}f^{p_0-2l}(f^2-1)^{l},
\end{equation}
where $[p_0/2]$ denotes the greatest integer not exceeding $p_0/2$; solutions of \eqref{yanagiharaeq0    5} are represented as
\begin{equation}\label{yanagiharaeq0    17}
f=\frac{1}{2}(\lambda+\lambda^{-1}),
\end{equation}
where $\lambda$ is given by
\begin{equation*}
\lambda=(\pm i)^{-\frac{1}{p_0}}\exp[\pi(z)(p_0+1)^z],
\end{equation*}
or
\begin{equation*}
\begin{split}
\lambda=(\pm i)^{\frac{1}{2+p_0}}\exp[\pi(z)(-p_0-1)^z],
\end{split}
\end{equation*}
where $\pi(z)$ is an arbitrary entire periodic function with period~1.

\end{enumerate}
\end{theorem}

We remark that the solutions with the form \eqref{yanagiharaeq0    21 hj} of equation \eqref{yanagiharaeq0    9} are not given in \cite{NakamuraYanagihara1989difference} where an existence theorem for entire solutions of \eqref{yanagiharaeq0    9} is stated instead; we also remark that the polynomial in \eqref{yanagiharaeq0    6} has different form from the one in \cite[Theorem~4(a)]{NakamuraYanagihara1989difference} since we have chosen different form of the solutions \eqref{yanagiharaeq0    17}.

Equations \eqref{Theorem  1a  eq0} and \eqref{yanagiharaeq0    8}, as well as their solutions \eqref{Theorem  1a  eq1} and \eqref{yanagiharaeq0    8be}, are apparently of the same form. We note that in Theorem~\ref{Theorem  1c}, when $n\geq 3$, we only have equation \eqref{yanagiharaeq0    8}. In fact, when $n\geq 3$, if some $\alpha_i$ in \eqref{P} is non-zero, then by Lemmas~\ref{basiclemma} and~\ref{keylemma} we have a contradiction to the inequality \eqref{multiplicityinequality} since $\omega\alpha_i$ is a completely ramified rational function of $f$, where $\omega$ is the $n$-th root of~1, and $\infty$ is a Picard exceptional rational function of $f$. Then \eqref{yanagiharaeq0    8} follows by applying Lemma~\ref{basiclemma0} and Picard's theorem.

\renewcommand{\proofname}{Proof of Theorem~\ref{Theorem  1c}.}
\begin{proof}
From Lemma~\ref{basiclemma} we know that $\infty$ is a Picard exceptional rational function of $f$ and $N_c\leq 2$. If some $\alpha_i$ in \eqref{P} is zero, then by Lemma~\ref{basiclemma0} it follows that $0$ is also a Picard exceptional rational function of $f$. Then by Picard's theorem we conclude that $P(z,f)$ cannot have any non-zero root and thus $P_0(z,f)^n=a_p$, i.e., we have the following equation:
\begin{equation}\label{yanagiharaeq0    3bev}
\overline{f}^n=a_pf^p.
\end{equation}
Similarly as in the proof of Theorem~\ref{Theorem  1a}, we choose a rational function $\alpha$ such that $\alpha f$ has no zeros and poles and then write equation \eqref{yanagiharaeq0    3bev} as
\begin{equation}\label{yanagiharaeq0    3}
\overline{f}^n=cf^p,
\end{equation}
where $c=a_p\frac{\overline{\alpha}^n}{\alpha^p}$ is a non-zero constant and $f$ above is a zero-free entire function. Moreover, solutions of \eqref{yanagiharaeq0    3} can be solved as
\begin{equation*}
f=c^{\frac{1}{n-p}}\exp[\pi(z)(p/n)^z],
\end{equation*}
where $\pi(z)$ is an arbitrary entire periodic function with period~1. This is the first part of Theorem~\ref{Theorem  1c}.

From now on, we suppose that none of $\alpha_i$ in \eqref{P} is zero. By Lemmas~\ref{basiclemma} and~\ref{keylemma} it follows that $\omega\alpha_i$ is a completely ramified rational function of $f$, where $\omega$ is the $n$-th root of~1. Since $\infty$ is a Picard exceptional rational function of $f$, then by the inequality \eqref{multiplicityinequality} we must have $n=2$ and $N_c\leq 2$. Below we consider the two cases where $N_c=1$ and $N_c=2$ separately.

\vskip 4pt

\noindent\textbf{Case~1:} $N_c=1$.

\vskip 4pt

In this case, $\infty$ is a Picard exceptional rational function of $f$ and $p$ is odd. Therefore, we have
\begin{equation}\label{proof4 eqre   1}
\overline{f}^2=P_0(z,f)^2(f-\alpha_1)^{k_1},
\end{equation}
where $\alpha_1\not=0$ and $k_1$ is an odd integer. By Lemmas~\ref{basiclemma} and~\ref{keylemma} it follows that $\pm\alpha_1$ are both completely ramified rational functions of $f$ with multiplicities~$2$. Then by Lemma~\ref{keylemma0}, we must have $k_1=1$ and $P_0(z,f)$ is a polynomial in $f$ with simple roots only. We may let $\alpha_1=1$ by doing a linear transformation $f\to \alpha_1f$. We consider
\begin{equation}\label{proof4 eqre   1a0}
\overline{f}^2-1=P_0(z,f)^2(f-1)-1.
\end{equation}
The RHS of \eqref{proof4 eqre   1a0} is a polynomial in $f$ with odd degree and thus has at least one root, say $\gamma_1$, of odd order. Since $f$ has no other completely ramified rational functions besides $\pm 1$, then by Lemma~\ref{keylemma1} we conclude that $\gamma_1$ must be $-1$ and there is only one such $\gamma_1$; moreover, the RHS of \eqref{proof4 eqre   1a0} is of the form $P_1(z,f)^2(f+1)$ for some polynomial $P_1(z,f)$ in $f$ with simple roots only. Now we have
\begin{equation}\label{proof4 eqre   1   pu}
\overline{f}^2=P_0(z,f)^2(f-1),
\end{equation}
and further that
\begin{equation}\label{proof4 eqre   1a0  pu}
\overline{f}^2-1=P_1(z,f)^2(f+1).
\end{equation}
From \eqref{proof4 eqre   1   pu} and \eqref{proof4 eqre   1a0   pu}, we see that the degree of $P_1(z,f)$ is $p_0$, at least~$1$. Put
\begin{equation}\label{proof4 eqre   1a0   nh0}
f=\frac{1}{2}(\lambda+\lambda^{-1}).
\end{equation}
Since both $\pm 1$ have multiplicities~$2$, we may write $f+1=g^2$ with an algebroid function $g$ and $g$ has at most finitely many algebraic branch points. It follows that the RHS of equation \eqref{proof4 eqre   1a0  pu} becomes $[P_1(z,g^2-1)g]^2$, which implies that $\lambda$ is an algebroid function with at most finitely many algebraic branch points. Moreover, $0$ and $\infty$ are both Picard exceptional rational functions of $\lambda$. Put
\begin{equation}\label{proof4 eqre   1a0   nh0ghu0}
\lambda=\delta^2.
\end{equation}
Then $\delta$ is an algebroid function with at most finitely many algebraic branch points. Now it follows from equation \eqref{proof4 eqre   1a0  pu} that
\begin{equation}\label{proof4 eqre   1aa0}
\frac{1}{2}(\overline{\delta}^2-\overline{\delta}^{-2})=P_1(z,f)\frac{\delta^2+1}{2^{1/2}\delta}.
\end{equation}
By solving equation \eqref{proof4 eqre   1aa0} together with \eqref{proof4 eqre   1   pu} and \eqref{proof4 eqre   1a0   pu}, we get
\begin{equation}\label{proof4 eqre   1aa1}
\overline{\delta}^2=\frac{P_1\left(z,\frac{\delta^4+1}{2\delta^2}\right)(\delta^2+1)\pm P_0\left(z,\frac{\delta^4+1}{2\delta^2}\right)(\delta^2-1)}{2^{1/2}\delta}:=\frac{P_{11}(z,\delta)}{2^{1/2}\delta(2\delta^2)^{p_0}},
\end{equation}
where $P_{11}(z,\delta)$ is a polynomial in $\delta$ of degree at most $4p_0+2$. Since~0 is a Picard exceptional rational function of $\lambda$, then by Picard's theorem we see from equation \eqref{proof4 eqre   1aa1} that $P_{11}(z,\delta)$ cannot have any non-zero root. By the Valiron--Mohon'ko identity \cite{Valiron1931,mohonko1971} (see also \cite{Laine1993}), we have from \eqref{proof4 eqre   1   pu}, \eqref{proof4 eqre   1a0   nh0} and \eqref{proof4 eqre   1a0   nh0ghu0} that $4T(r,\overline{\delta})=(4p_0+2)T(r,\delta)+O(\log r)$. Therefore, we have either
\begin{equation}\label{proof4 eqre   1aa2}
\begin{split}
\overline{\delta}^2=T_0\delta^{2p_0+1},
\end{split}
\end{equation}
or
\begin{equation}\label{proof4 eqre   1aa3}
\begin{split}
\overline{\delta}^2=T_0\delta^{-2p_0-1},
\end{split}
\end{equation}
where $T_0$ is an algebraic function. We write
\begin{equation}\label{proof4 eqre   1aa1  qw0}
P_1(z,f) = c_{p_0}f^{p_0} + c_{p_0-1}f^{p_0-1}+\cdots+c_0,
\end{equation}
where $c_{p_0}$, $\cdots$, $c_0$ are algebraic functions and $c_{p_0}\not=0$. If we have \eqref{proof4 eqre   1aa2}, then by substituting \eqref{proof4 eqre   1a0   nh0}, \eqref{proof4 eqre   1a0   nh0ghu0} and \eqref{proof4 eqre   1aa2} into \eqref{proof4 eqre   1aa0} and then comparing the terms on both sides of the resulting equation together with \eqref{proof4 eqre   1aa1  qw0}, we get
\begin{equation}\label{proof4 eqre   1aa4}
\overline{\delta}^2=2^{1/2}2^{-p_0}c_{p_0}\delta^{2p_0+1},
\end{equation}
and
\begin{equation}\label{proof4 eqre   1aa5}
\overline{\delta}^2=-2^{-1/2}2^{p_0}(1/c_{p_0})\delta^{2p_0+1}.
\end{equation}
On the other hand, if we have \eqref{proof4 eqre   1aa3}, then similarly as above from \eqref{proof4 eqre   1a0   nh0}, \eqref{proof4 eqre   1a0   nh0ghu0}, \eqref{proof4 eqre   1aa0} and \eqref{proof4 eqre   1aa1  qw0} we get
\begin{equation}\label{proof4 eqre   1aa6}
\overline{\delta}^2=2^{1/2}2^{-p_0}c_{p_0}\delta^{-2p_0-1},
\end{equation}
and
\begin{equation}\label{proof4 eqre   1aa7}
\overline{\delta}^2=-2^{-1/2}2^{p_0}(1/c_{p_0})\delta^{-2p_0-1}.
\end{equation}
We obtain from \eqref{proof4 eqre   1aa4} and \eqref{proof4 eqre   1aa5}, as well as from \eqref{proof4 eqre   1aa6} and \eqref{proof4 eqre   1aa7}, that
$c_{p_0}=\pm i2^{p_0-\frac{1}{2}}$
and it follows that
$T_0=\pm i$.
Thus the solution $f$ of \eqref{proof4 eqre   1   pu} is represented by \eqref{proof4 eqre   1a0   nh0} and \eqref{proof4 eqre   1a0   nh0ghu0} with $\delta$ such that
\begin{equation}\label{proof4 eqre   1aa10}
\overline{\delta}^2=\pm i\delta^{2p_0+1}.
\end{equation}
or
\begin{equation}\label{proof4 eqre   1aa11}
\overline{\delta}^2=\pm i\delta^{-2p_0-1}.
\end{equation}
Note that $\delta$ has at most finitely many zeros, poles and branch points. Then we can solve $\delta$ from \eqref{proof4 eqre   1aa10} and \eqref{proof4 eqre   1aa11} as
\begin{equation}\label{proof4 eqre   1aa12}
\delta=(\pm i)^{\frac{1}{1-2p_0}}\exp[\pi(z)(p_0+1/2)^z],
\end{equation}
or
\begin{equation}\label{proof4 eqre   1aa13}
\delta=(\pm i)^{\frac{1}{3+2p_0}}\exp[\pi(z)(-p_0-1/2)^z],
\end{equation}
respectively, where $\pi(z)$ is an arbitrary entire periodic function with period~1. We conclude that solutions of equation \eqref{proof4 eqre   1   pu} are given by \eqref{proof4 eqre   1a0   nh0}, \eqref{proof4 eqre   1a0   nh0ghu0} with \eqref{proof4 eqre   1aa12} or \eqref{proof4 eqre   1aa13}.

Now we determine the polynomial $P_0(z,f)$ in $f$ in \eqref{proof4 eqre   1   pu} by using \eqref{proof4 eqre   1a0   nh0}, \eqref{proof4 eqre   1a0   nh0ghu0}, \eqref{proof4 eqre   1aa10} and \eqref{proof4 eqre   1aa11}. From \eqref{proof4 eqre   1   pu} and \eqref{proof4 eqre   1aa10} we have
\begin{equation}\label{proof4 eqre   1aa14}
\frac{1}{2}(\pm i)(\delta^{2p_0+1}-\delta^{-2p_0-1})=P_0(z,f)\frac{\delta^2-1}{2^{1/2}\delta}.
\end{equation}
Moreover, by \eqref{proof4 eqre   1a0   nh0} and \eqref{proof4 eqre   1a0   nh0ghu0} we have
\begin{equation}\label{proof4 eqre   1aa15}
\begin{split}
\delta^2&=f\pm(f^2-1)^{1/2},\\
\delta^{-2}&=f\mp(f^2-1)^{1/2},
\end{split}
\end{equation}
and also that
\begin{equation}\label{proof4 eqre   1aa16}
\begin{split}
\delta^{2p_0+1}-\delta^{-2p_0-1}=\frac{\delta^{4p_0+2}-1}{\delta^{2p_0+1}}.
\end{split}
\end{equation}
From \eqref{proof4 eqre   1aa14}, \eqref{proof4 eqre   1aa15} and \eqref{proof4 eqre   1aa16} we see that $P_0(z,f)=P_0(f)$ is a polynomial in $f$ with constant coefficients and
\begin{equation}\label{proof4 eqre   1aa17}
\begin{split}
P_0(f)=\frac{\pm i}{2^{1/2}}\left\{\frac{\delta^{2p_0+2}-\delta^{-2p_0}}{\delta^2-1}\right\}=\frac{\pm i}{2^{1/2}}\left\{\frac{[f\pm(f^2-1)^{1/2}]^{p_0+1}-[f\mp(f^2-1)^{1/2}]^{p_0}}{f\pm(f^2-1)^{1/2}-1}\right\}.
\end{split}
\end{equation}
Now the polynomial $P_0(f)$ with constant coefficients takes the form in \cite[Theorem~4(b)]{NakamuraYanagihara1989difference}. But the proof there is rather complicated and here we give a simple one. Note that $f$ has no finite Picard exceptional values. 
It suffices to take the value $f=\cos x$, where $x\in(-\pi/2,\pi/2)$ is real. By substituting $f=\cos x$ into \eqref{proof4 eqre   1aa17}, we get
\begin{equation}\label{proof4 eqre   1aa17fg0}
\begin{split}
P_0(\cos x)&=\frac{\pm i}{2^{1/2}}\left\{\frac{[\cos x\pm i\sin x]^{p_0+1}-[\cos x\mp i\sin x]^{p_0}}{\cos x\pm i\sin x-1}\right\}\\
&=\frac{\pm i}{2^{1/2}}\left\{\frac{[\cos(x/2)\pm i\sin(x/2)]^{2(p_0+1)}-[\cos(x/2)\mp i\sin(x/2)]^{2p_0}}{[\cos(x/2)\pm i\sin(x/2)]^2-1}\right\}.
\end{split}
\end{equation}
By the well-known de Moivre's formula and the basic formula $\sin(x+y)+\sin(x-y)=2\sin x\cos y$ on real trigonometric functions, we deduce from \eqref{proof4 eqre   1aa17fg0} that
\begin{equation*}
\begin{split}
P_0(\cos x)&=\frac{\pm i}{2^{1/2}}\left\{\frac{e^{\pm i(p_0+1)x}-e^{\mp ip_0x}}{e^{\pm ix}-1}\right\}=\frac{\pm i}{2^{1/2}}\left\{\frac{e^{i(p_0+1/2)x}-e^{-i(p_0+1/2)x}}{e^{ix/2}-e^{-ix/2}}\right\}\\
&=\frac{\pm i}{2^{1/2}}\left\{\frac{\sin(p_0+1/2)x}{\sin (x/2)}\right\}=\frac{\pm i}{2^{1/2}}\left\{\frac{2\sin(p_0+1/2)x\cos(x/2)}{2\sin(x/2)\cos(x/2)}\right\}\\
&=\frac{\pm i}{2^{1/2}}\left\{\frac{\sin(p_0+1)x+\sin p_0x}{\sin x}\right\}=\frac{\pm i}{2^{1/2}}\left\{\frac{\sin(p_0+1)x}{\sin x}+\frac{\sin p_0x}{\sin x}\right\}.
\end{split}
\end{equation*}
Note that in the second step the equation holds for both choices of the signs $\pm$ in the exponential function. Denote $U_{p_0}(\cos x):=\frac{\sin(p_0+1)x}{\sin x}$. Then $U_{p_0}$ is the \emph{Chebyshev polynomial of the second kind}~\cite[p.~184]{Erdelyi1953}. Thus the coefficients of $P_0(f)$ are independent of the choice of $x$. We conclude that
\begin{equation*}
\begin{split}
P_0(f)=\frac{\pm i}{2^{1/2}}[U_{p_0}(f)+U_{p_0-1}(f)],
\end{split}
\end{equation*}
where $U_{p_0}$ and $U_{p_0-1}$ are the Chebyshev polynomials of the second kind. This corresponds to the second part of Theorem~\ref{Theorem  1c}.

\vskip 4pt

\noindent\textbf{Case~2:} $N_c=2$.

\vskip 4pt

In this case, $\infty$ is a Picard exceptional rational function of $f$ and $p$ is even. Moreover, by the inequality \eqref{multiplicityinequality} it follows that all $\alpha_i$ are completely ramified rational functions of $f$ with multiplicities~$2$. Therefore, by Lemma~\ref{keylemma0} we have
\begin{equation}\label{proof3 eqre   1}
\overline{f}^2=P_0(z,f)^2(f-\alpha_1)(f-\alpha_2),
\end{equation}
where $\alpha_1$ and $\alpha_2$ are in general both non-zero algebraic functions and distinct from each other, and $P_0(z,f)$ is a polynomial in $f$ with simple roots only.
By Lemmas~\ref{basiclemma} and~\ref{keylemma}, it follows that $\pm\alpha_1$ and $\pm\alpha_2$ are all completely ramified rational functions of $f$ and so by Theorem~\ref{completelyrm} we must have $\alpha_1+\alpha_2=0$. We may let $\alpha_1=1$ and $\alpha_2=-1$ by doing a linear transformation $f\to \alpha_1f$. We consider
\begin{equation}\label{proof3 eqre   1aa0}
\overline{f}^2-1=P_0(z,f)^2(f^2-1)-1.
\end{equation}
Since $\pm 1$ are both completely ramified rational functions of $f$ with multiplicities~2 and $f$ has no other completely ramified rational functions, then by Lemma~\ref{keylemma1} we conclude that the RHS of \eqref{proof3 eqre   1aa0} is of the form $P_1(z,f)^2$ for some polynomial $P_1(z,f)$ in $f$ with simple roots only. Now we have
\begin{equation}\label{proof3 eqre   3}
\overline{f}^2=P_0(z,f)^2(f^2-1),
\end{equation}
and further that
\begin{equation}\label{proof3 eqre   1aa1  fe}
\overline{f}^2-1=P_1(z,f)^2.
\end{equation}
Put
\begin{equation}\label{proof3 eqre   4}
f=\frac{1}{2}(\lambda+\lambda^{-1}).
\end{equation}
Then from \eqref{proof3 eqre   1aa1  fe} we see that $\lambda$ is an algebroid function with at most finitely many algebraic branch points. Moreover, $0$ and $\infty$ are both Picard exceptional rational functions of $\lambda$. It follows from \eqref{proof3 eqre   1aa1  fe} that
\begin{equation}\label{proof3 eqre   1aa1  fe  tre}
\frac{1}{2}(\overline{\lambda}-\overline{\lambda}^{-1})=P_1(z,f).
\end{equation}
From \eqref{proof3 eqre   3} and \eqref{proof3 eqre   1aa1  fe}, we see that the degree of $P_1(z,f)$ in $f$ is $p_0+1$, which is greater than or equal to $2$. By solving equation \eqref{proof3 eqre   1aa1  fe  tre} together with \eqref{proof3 eqre   3} and \eqref{proof3 eqre   1aa1  fe}, we get
\begin{equation}\label{proof3 eqre   6}
\overline{\lambda}=P_1(z,f)\pm P_0(z,f)(f^2-1)^{1/2}=\frac{P_1\left(z,\frac{\lambda^2+1}{2\lambda}\right)(2\lambda)\pm P_0\left(z,\frac{\lambda^2+1}{2\lambda}\right)(\lambda^2-1)}{2\lambda}:=\frac{P_{11}(z,\lambda)}{(2\lambda)^{p_0+2}},
\end{equation}
where $P_{11}(z,\lambda)$ is a polynomial in $\lambda$ of degree at most $2p_0+3$. Since~0 is a Picard exceptional rational function of $\lambda$, then by Picard's theorem we see from equation \eqref{proof3 eqre   6} that $P_{11}(z,\lambda)$ cannot have any non-zero root. Moreover, by the Valiron--Mohon'ko identity \cite{Valiron1931,mohonko1971} (see also \cite{Laine1993}), we have from \eqref{proof3 eqre   1aa1  fe  tre} that $2T(r,\overline{\lambda})=(p_0+1)T(r,f)+O(\log r)=2(p_0+1)T(r,\lambda)+O(\log r)$.
Therefore, similarly as in previous case, we may show that the solution $f$ of equation \eqref{proof3 eqre   3} is represented by \eqref{proof3 eqre   4} with $\lambda$ such that
\begin{equation}\label{proof3 eqre   19}
\overline{\lambda}= \pm i\lambda^{p_0+1},
\end{equation}
or
\begin{equation}\label{proof3 eqre   20}
\overline{\lambda}= \pm i\lambda^{-p_0-1}.
\end{equation}
Note that $\lambda$ has at most finitely many zeros, poles and branch points. Then we can solve $\lambda$ from equations \eqref{proof3 eqre   19} and \eqref{proof3 eqre   20} as
\begin{equation}\label{proof3 eqre   23}
\lambda=(\pm i)^{-\frac{1}{p_0}}\exp[\pi(z)(p_0+1)^z],
\end{equation}
or
\begin{equation}\label{proof3 eqre   24}
\lambda=(\pm i)^{\frac{1}{2+p_0}}\exp[\pi(z)(-p_0-1)^z],
\end{equation}
respectively, where $\pi(z)$ is an arbitrary entire periodic function with period~1. We conclude that solutions of equation \eqref{proof3 eqre   3} are given by \eqref{proof3 eqre   4} with \eqref{proof3 eqre   23} or \eqref{proof3 eqre   24}.

Now we determine the polynomial $P_0(z,f)$ in $f$ in \eqref{proof3 eqre   3} by using \eqref{proof3 eqre   4} with \eqref{proof3 eqre   19} or with \eqref{proof3 eqre   20}. From \eqref{proof3 eqre   3} and \eqref{proof3 eqre   19}, or \eqref{proof3 eqre   3} and \eqref{proof3 eqre   20}, we have
\begin{equation}\label{proof3 eqre   27}
\frac{\pm i}{2}(\lambda^{p_0+1}-\lambda^{-p_0-1})=\frac{1}{2}(\lambda-\lambda^{-1})P_0(z,f).
\end{equation}
By \eqref{proof3 eqre   4}, we have
\begin{equation}\label{proof3 eqre   28}
\begin{split}
\lambda&=f\pm(f^2-1)^{1/2},\\
\lambda^{-1}&=f\mp(f^2-1)^{1/2},
\end{split}
\end{equation}
and it follows that
\begin{equation}\label{proof3 eqre   29}
\begin{split}
\lambda^{p_0+1}-\lambda^{-p_0-1}
&=\sum_{k=0}^{p_0+1}\binom{p_0+1}{k}f^{p_0+1-k}\left\{[\pm(f^2-1)^{1/2}]^{k}-[\mp(f^2-1)^{1/2}]^{k}\right\},\\
&=2\sum_{l=0}^{[p_0/2]}\binom{p_0+1}{2l+1}f^{p_0-2l}[\pm(f^2-1)^{1/2}]^{2l+1}.
\end{split}
\end{equation}
From \eqref{proof3 eqre   27}, \eqref{proof3 eqre   28} and \eqref{proof3 eqre   29} we see that $P_0(z,f)=P_0(f)$ is a polynomial in $f$ with constant coefficients and
\begin{equation*}
P_0(f)=\pm i\sum_{l=0}^{[p_0/2]}\binom{p_0+1}{2l+1}f^{p_0-2l}(f^2-1)^{l},
\end{equation*}
where $[p_0/2]$ denotes the greatest integer not exceeding $p_0/2$. This corresponds to the third part of Theorem~\ref{Theorem  1c} and also completes the proof.

\end{proof}

\section{Equation \eqref{yanagiharaeq0} with $q\geq 1$ and $d=\max\{p,q\}\geq 2$}\label{Proof 2}

\subsection{Equation \eqref{yanagiharaeq0} with $q\geq 1$ and $n>d\geq 2$}

For the case $q\geq 1$ and $n>d\geq 2$ of equation \eqref{yanagiharaeq0}, we prove the following

\begin{theorem}\label{Theorem  2a}
Suppose that $q\geq 1$ and $n>d\geq 2$. Let $f$ be a transcendental meromorphic solution of equation \eqref{yanagiharaeq0}. Then there exists a rational function $\alpha$ such that the linear transformation $f\rightarrow \alpha f$ reduces \eqref{yanagiharaeq0} into
\begin{equation}\label{Theorem  2a  eq0}
\overline{f}^n=cf^{-q},
\end{equation}
where $c$ is a non-zero constant. Moreover, solutions of equation \eqref{Theorem  2a  eq0} are represented as
\begin{equation}\label{Theorem  2a  eq1}
f=c^{\frac{1}{n+q}}\exp[\pi(z)(-q/n)^z],
\end{equation}
where $\pi(z)$ is an arbitrary entire periodic function with period~1.
\end{theorem}

\renewcommand{\proofname}{Proof of Theorem~\ref{Theorem  2a}.}
\begin{proof}
Suppose that $P(z,f)$ has a non-zero root, say $\alpha_1$. Since $n\geq 3$, then by Lemmas~\ref{basiclemma} and~\ref{keylemma} it follows that $\omega\alpha_1$ is a completely ramified rational function of $f$, where $\omega$ is the $n$-th root of~1. Then by the inequality \eqref{multiplicityinequality} we conclude that $n=3$ or $n=4$. In particular, when $n=3$ we see that $\eta\alpha_1$ has multiplicity~$3$ since the order $k_1$ of the root $\alpha_1$ equals~1 or~2 under the assumption that $n>d$, where $\eta$ is the cubic root of~1. Thus by the inequality \eqref{multiplicityinequality}, when $n=3$ we have $\sum_{i=1}^3\Theta(\eta_i\alpha_1,f)=2$, where $\eta_i$ are the three numbers such that $\eta_i^3=1$. On the other hand, when $n=4$ we have $\sum_{i=1}^4\Theta(\omega_i\alpha_1,f)=2$, where $\omega_i$ are the four numbers such that $\omega_i^4=1$. By Lemma~\ref{basiclemma} and the inequality \eqref{multiplicityinequality} we see that~0 cannot be a root of $P(z,f)$ or $Q(z,f)$. Now, if $p\not=q$, then by Lemma~\ref{basiclemma} it follow that $\infty$ is either a Picard exceptional rational function of $f$ or a completely ramified function of $f$, a contradiction to the inequality \eqref{multiplicityinequality}. Therefore, $p=q$. From the above reasoning, when $n=3$, we have $p=q=2$ and by Lemma~\ref{keylemma0} we see that each of the roots of $P(z,f)$ and $Q(z,f)$ is simple since $0$ is not a completely ramified rational function of $f$, i.e., we have $N_c=4$, a contradiction to Lemma~\ref{basiclemma}. When $n=4$, we have $p=q=2$ or $p=q=3$, and by Lemma~\ref{keylemma0} we see that each of the roots of $P(z,f)$ and $Q(z,f)$ has double order since $0$ is not a completely ramified rational function of $f$, but it follows that the case where $p=q=3$ is impossible and when $p=q=2$ we have a contradiction to our assumption that at least one of $\alpha_i$ and $\beta_j$ in \eqref{P} and \eqref{Q} has order that is not divided by $n$. We conclude that $P(z,f)$ does not have any non-zero root. Similarly, $Q(z,f)$ cannot have any non-zero root either. Since $q\geq 1$, then by assumption we must have $P(z,f)=a_p$, i.e., we have the following equation:
\begin{equation}\label{Theorem  2a  eq0in pr}
\overline{f}^n= a_pf^{-q}.
\end{equation}
We claim that $f$ has at most finitely many poles. Suppose on the contrary that $f$ has infinitely many poles. Let $z_0\in \mathbb{C}$ be a pole of $f$ with multiplicity $m\in\mathbb{Z}^{+}$. We may choose $z_0$ such that $|z_0|$ is large enough so that $a_p$ has no poles or zeros outside of $\{z\in \mathbb{C}: |z|<|z_0|\}$. Then from \eqref{Theorem  2a  eq0in pr} we see that $z_0+1$ is a zero of $f$ of order $qm/n$ and it follows that $z_0+2$ is a pole of $f$ of order $q^2m/n^2$. By iteration we have that $f$ has a pole of order $q^{2s}m/n^{2s}$ at the point $z_0+2s$, $s\in \mathbb{N}$. Since $(n,q)=1$, then by letting $s\to\infty$, it follows that there is necessarily a branch point at some $z_0+2s_0$, $s_0\in \mathbb{N}$, a contradiction to our assumption that $f$ is meromorphic. Therefore, $f$ has at most finitely many poles. From \eqref{Theorem  2a  eq0in pr} we also have that $f$ has at most finitely many zeros.
Similarly as in the proof of Theorem~\ref{Theorem  1a}, we choose a rational function $\alpha$ such that $\alpha f$ has no zeros and poles and then write equation \eqref{Theorem  2a  eq0in pr} as
\begin{equation}\label{Theorem  2a  eq0in prewq}
\overline{f}^n= cf^{-q},
\end{equation}
where $c=a_p\overline{\alpha}^n\alpha^q$ is a non-zero constant and $f$ above is a zero-free entire function. By taking the logarithm on both sides of \eqref{Theorem  2a  eq0in prewq}, then $g=\log f$ is entire and satisfies
\begin{equation*}
n\overline{g}=\log c-qg.
\end{equation*}
Therefore, we can solve $f$ from \eqref{Theorem  2a  eq0in prewq} as
\begin{equation*}
f=c^{\frac{1}{n+q}}\exp[\pi(z)(-q/n)^z],
\end{equation*}
where $\pi(z)$ is an arbitrary entire periodic function with period~1. This completes the proof.
\end{proof}

\subsection{Equation \eqref{yanagiharaeq0} with $q\geq 1$ and $2\leq n<d$}

In this section, we consider the two cases $n\nmid |p-q|$ and $n\mid |p-q|$ of equation \eqref{yanagiharaeq0} separately. For the case $n\nmid |p-q|$, we actually have $N_c=1$; we prove the following

\begin{theorem}\label{Theorem  2b}
Suppose that $q\geq 1$ and $2\leq n<d$ and $n\nmid |p-q|$. Let $f$ be a transcendental meromorphic solution of equation \eqref{yanagiharaeq0}. Then there exists a rational function $\alpha$ such that by doing a linear transformation $f\rightarrow \alpha f$, we have either
\begin{enumerate}
\item[(1)] equation \eqref{yanagiharaeq0} reduces into
\begin{equation}\label{Theorem  2b  eq2}
\overline{f}^n = cf^{-q},
\end{equation}
and solutions of equation \eqref{Theorem  2b  eq2} are represented as
\begin{equation}\label{Theorem  2b  eq3}
f=c^{\frac{1}{n+q}}\exp[\pi(z)(-q/n)^z],
\end{equation}
where $\pi(z)$ is an arbitrary entire periodic function with period~1; or

\item[(2)] $q$ is even and $q=2q_0$, $q_0\geq1$, and equation \eqref{yanagiharaeq0} reduces into
\begin{equation}\label{Theorem  2b  eq0}
\overline{f}^2 = \frac{P_0(z,f)^2}{Q_0(z,f)^2}(f-1)^{k_1},
\end{equation}
where $k_1$ is an odd integer and $2p_0+k_1< 2q_0$; moreover, we have
\begin{equation}\label{Theorem  2b  eq0a2  fuja}
\begin{split}
Q_0(z,f)=Q_1(z,f)^2(f+1)^{l_0}=\frac{1}{2i}[P_{011}(z,f)^2-P_{012}(z,f)^2(f-1)^{k_1}],
\end{split}
\end{equation}
where $l_0\in \mathbb{N}$ is zero or an odd integer, $Q_1(z,f)$ is a polynomial in $f$ and $P_{011}(z,f)$ and $P_{012}(z,f)$ are two polynomials in $f$ with no common roots such that $P_{011}(z,f)P_{012}(z,f)=P_0(z,f)$; solutions of \eqref{Theorem  2b  eq0} are represented as
\begin{equation*}
f=\frac{1}{2}(\delta^2+\delta^{-2}),
\end{equation*}
and $\delta$ is a function such that
\begin{equation}\label{Theorem  2b  eq0tyrrte010}
\overline{\delta}=\pm i^{1/2}\frac{P_{021}\left(z,\frac{\delta^4+1}{2\delta^{2}}\right)+\theta P_{022}\left(z,\frac{\delta^4+1}{2\delta^{2}}\right)\left(\frac{\delta^2-1}{2^{1/2}\delta}\right)^{k_1}}{P_{021}\left(z,\frac{\delta^4+1}{2\delta^{2}}\right)-\theta P_{022}\left(z,\frac{\delta^4+1}{2\delta^{2}}\right)\left(\frac{\delta^2-1}{2^{1/2}\delta}\right)^{k_1}}, \quad \theta=\pm1,
\end{equation}
when $l_0=0$, where $P_{021}(z,f)$ and $P_{022}(z,f)$ are two polynomials in $f$ with no common roots such that $P_{021}(z,f)P_{022}(z,f)=P_{012}(z,f)$, or such that
\begin{equation}\label{Theorem  2b  eq0tyrrte0100}
\overline{\delta}=\pm (-i)^{1/2}\frac{P_{023}\left(z,\frac{\delta^4+1}{2\delta^{2}}\right)\delta^{t_1}+\theta P_{024}\left(z,\frac{\delta^4+1}{2\delta^{2}}\right)}{P_{023}\left(z,\frac{\delta^4+1}{2\delta^{2}}\right)\delta^{t_1}-\theta P_{024}\left(z,\frac{\delta^4+1}{2\delta^{2}}\right)}, \quad \theta=\pm1,
\end{equation}
when $l_0>0$, where $t_1\in \mathbb{Z}\setminus\{0\}$ is an odd integer, $P_{023}(z,f)$ and $P_{024}(z,f)$ are two polynomials in $f$ with no common roots such that $P_{023}(z,f)P_{024}(z,f)=P_{011}(z,f)$.
\end{enumerate}
\end{theorem}

Equations \eqref{Theorem  2a  eq0} and \eqref{Theorem  2b  eq2}, as well as their solutions \eqref{Theorem  2a  eq1} and \eqref{Theorem  2b  eq3}, are apparently of the same form. We note that in Theorem~\ref{Theorem  2b}, when $n\geq 3$, we only have equation \eqref{Theorem  2b  eq2}. In fact, when $n\geq 3$, if some $\alpha_i$ in \eqref{P} or $\beta_j$ in \eqref{Q} is non-zero, say $\alpha_i\not=0$ for some $i$, then by Lemmas~\ref{basiclemma} and~\ref{keylemma} we have a contradiction to the inequality \eqref{multiplicityinequality} since $\omega\alpha_i$ is a completely ramified rational function of $f$, where $\omega$ is the $n$-th root of~1, and $\infty$ is either a Picard exceptional rational function or a completely ramified rational function of $f$ with multiplicity at least~$n$. Then the reasoning in the proof of Theorem~\ref{Theorem  2a} yields equation \eqref{Theorem  2b  eq2}.

In the autonomous case, the RHS of equation \eqref{Theorem  2b  eq0tyrrte010} or \eqref{Theorem  2b  eq0tyrrte0100} becomes a rational term $R_0(\delta)$ in $\delta$ after multiplying $(2^{1/2}\delta)^{q_0}$ to both the numerator and the denominator and thus $\overline{\delta}=R_0(\delta)$ always has a meromorphic solution for any given $P_0(z,f)$ and $Q_1(z,f)$ such that the relation in \eqref{Theorem  2b  eq0a2  fuja} holds, as mentioned in the introduction.
We note that when $l_0>0$, the polynomials $P_0(z,f)$ and $Q_0(z,f)$ satisfying the relation in \eqref{Theorem  2b  eq0a2  fuja} exist. For example, for the two polynomials $P_0(f)$ and $P_1(f)$ satisfying \eqref{proof4 eqre   1   pu} and \eqref{proof4 eqre   1a0  pu} in the proof of Theorem~\ref{Theorem  1c}, we have $iP_0(f)^2(f-1)-iP_1(f)^2(f+1)=i$. In the simplest case $p_0=k_1=l_0=1$, we have $2i(f-1/2)^2(f+1)-2i(f+1/2)^2(f-1)=i$.

\renewcommand{\proofname}{Proof of Theorem~\ref{Theorem  2b}.}
\begin{proof}
First, we show that under the assumptions of Theorem~\ref{Theorem  2b} the case $p>q$ cannot occur.
Since $n\nmid(p-q)$, then from the proof of Lemma~\ref{basiclemma0}, we see that $\infty$ is a Picard exceptional rational function of $f$ no matter whether or not some $\alpha_i$ in \eqref{P} is zero. Let $\beta$ be any root of $Q(z,f)$. Then from \eqref{yanagiharaeq0} we see that the equation $f-\beta=0$ has at most finitely many roots and so $\beta$ is a Picard exceptional rational function of $f$. By Picard's theorem we see that there is only one such $\beta$. Then by Lemma~\ref{basiclemma} and the inequality \eqref{multiplicityinequality} we conclude that \eqref{yanagiharaeq0} takes the following form:
    \begin{equation}\label{eqnonmero}
    \overline{f}^n=\frac{P_0(z,f)^n}{(f-\beta)^q},
    \end{equation}
where $\beta$ is a rational function. Moreover, since $p>q$, we see that $\beta\not\equiv0$ since otherwise all the roots of $P_0(z,f)$ are Picard exceptional rational functions of $f$, a contradiction to Picard's theorem. Denote $g=(f-\beta)^{1/n}$. Then $g$ is an algebroid function with at most finitely many algebraic branch points and it follows that $f=g^n+\beta$. Then we can rewrite equation \eqref{eqnonmero} as follows
    \begin{equation}\label{eqnonmero ne1}
    \overline{g}^{n}g^{q}=P_0(z,g^n+\beta)-\overline{\beta}g^{q}.
    \end{equation}
Note that $0$ and $\infty$ are both Picard exceptional rational functions of $g$. Let $u_0$ be a function such that
    \begin{equation*}
    P_0(z,u_0^n+\beta)-\overline{\beta}u_0^{q}=0.
    \end{equation*}
Then $u_0$ is an algebraic function. Since $q<p$ and $P_0(z,f)$ has at least one root distinct from $\beta$, we see that the equation above has at least one non-zero root and from \eqref{eqnonmero ne1} we see that for the non-zero $u_0$ we have that $g-u_0=0$ has at most finitely many roots, i.e., $u_0$ is a Picard exceptional rational function of $g$, a contradiction to Picard's theorem. Therefore, the case where $p>q\geq 1$ and $n\nmid (q-p)$ cannot occur. In particular, since $p\not=q$, the above analysis also implies that the case where~0 is a root of $P(z,f)$ of order~$k_0$ such that $n\nmid k_0$ cannot occur; otherwise, by doing a bilinear transformation $f\rightarrow1/f$ to \eqref{yanagiharaeq0}, we always get
    \begin{equation}\label{cannotouu}
    \overline{f}^n=\frac{P_1(z,f)}{Q_1(z,f)},
    \end{equation}
where $P_1(z,f)$ is a polynomial in $f$ of degree $d$ and $Q_1(z,f)$ is a polynomial in $f$ of degree $d-k_0$, which is impossible from previous discussions.

Second, we show that $N_c=1$. Suppose on the contrary that $N_c\geq 2$. Then at least one of $\alpha_i$ and $\beta_j$ is non-zero, say $\alpha_i$. By Lemmas~\ref{basiclemma} and~\ref{keylemma}, $\omega\alpha_i$ is a completely ramified rational function of $f$, where $\omega$ is the $n$-th root of~1. Moreover, since $p<q$ and $n\nmid |p-q|$, then by Lemma~\ref{basiclemma} it follows that $\infty$ is also a completely ramified rational function of $f$ with multiplicity at least~$n$. By the inequality \eqref{multiplicityinequality}, we must have $n=2$. In this case, we see that $N_c$ is an odd integer. If $N_c\geq 3$, then none of $\alpha_i$ and $\beta_j$ is zero; otherwise, $0$ has multiplicity at least~2 and it follows that $\infty$ has multiplicity at least~4, a contradiction to the inequality \eqref{multiplicityinequality}. For convenience, we denote these $\alpha_i$ and $\beta_j$ by $\gamma_j$, $j=1,2,\cdots,k$. By Lemmas~\ref{basiclemma} and~\ref{keylemma}, $\pm\gamma_j$ are all completely ramified rational functions of $f$. Then by Theorem~\ref{completelyrm} we conclude that $\gamma_j^2$ must be equal to each other for all $j$, i.e., $\gamma_1^2=\gamma_2^2=\cdots=\gamma_k^2$ holds for $k\geq 3$, which is impossible. Therefore, we must have $N_c=1$.

Third, we show that~0 cannot be a root of $Q(z,f)$ of order~$l_0<q$ such that $n\nmid l_0$. Note that now we have $0\leq p<q$ under the assumption $p\not=q$. Otherwise, \eqref{yanagiharaeq0} can be written as
    \begin{equation}\label{PPrst0}
    \overline{f}^n=\frac{P(z,f)}{f^{l_0}Q_2(z,f)},
    \end{equation}
where $1\leq l_0\leq q-1$ satisfies $n\nmid l_0$, and $Q_2(z,f)$ is a polynomial in $f$ of degree $q-l_0$. Now we must have $n\mid (q-l_0)$ and $n\mid p$ by the fact that $N_c=1$ and it follows that $q-p=k_0n+l_0$ for some integer $k_0$. Suppose that $f$ has infinitely many zeros and let $z_0\in\mathbb{C}$ be a zero of $f$ with multiplicity $m\in\mathbb{Z}^{+}$. We may choose $z_0$ such that $|z_0|$ is large enough so that none of the coefficients of $P(z,f)$ and $Q(z,f)$ has poles or zeros outside of $\{z\in \mathbb{C}: |z|<|z_0|\}$. By \eqref{PPrst0}, $z_0+1$ is a pole of $f$ of order $l_0m/n$. It follows that $z=z_0+2$ is a zero of $f$ of order $l_0(k_0n+l_0)m/n^2$. Then, by iteration it follows that $z=z_0+2s$, $s\in\mathbb{N}$, is a zero of $f$ of order $(k_0n+l_0)^sl_0^{s}m/n^{2s}$. Since $n^2\nmid (k_0n+l_0)l_0$, then by letting $s\to\infty$, it follows that there is necessarily a branch point of $f$ at $z_0+2s_0$ for some $s_0\in\mathbb{N}$, a contradiction to our assumption that $f$ is meromorphic. Therefore, $f$ has at most finitely many zeros, i.e., $0$ is a Picard exceptional rational function of $f$. Also, from \eqref{PPrst0} we see that $f$ has at most finitely many poles since $p<q$ and then, since $l_0\leq q-1$, it follows that there exists another non-zero $\beta$ such that $\beta$ is a root of $Q(z,f)$ and $f-\beta=0$ has at most finitely many roots, that is to say, $f$ has at least~3 Picard exceptional rational functions, a contradiction to Picard's theorem. Therefore,~$0$ cannot be a root of $Q(z,f)$ of order $l_0<q$ such that $n\nmid l_0$ when assuming $n \nmid (q-p)$.

By combining all the above results together, we see that we only need to consider two cases of \eqref{yanagiharaeq0} under our assumption: (1) the case where $0$ is the only root of $Q(z,f)$; or (2) the case where $0$ is not a root of $Q(z,f)$ of order $0<l_0<q$ such that $n\nmid l_0$. In particular, in the latter case if $0$ is not a root of $Q_0(z,f)$, then from the previous discussions we can assume that $0$ is not a root of $P(z,f)$ of order $k_0$ such that $n\nmid k_0$.

Now, if $0$ is the only root of $Q(z,f)$ and $n\nmid q$, from the above reasoning we must have
\begin{eqnarray*}
\overline{f}^n = \frac{P_0(z,f)^n}{f^{q}}.
\end{eqnarray*}
In this case, since $p<q$ and $n\nmid (q-p)$, then by applying exactly the same analysis as in the previous
discussions on the case where $0$ is a root of $Q(z,f)$ of order $l_0<q$ and $n\nmid l_0$, we can show that $f$ has at most finitely many zeros and poles, i.e., $0$ and $\infty$ are both Picard exceptional rational functions of $f$. By Picard's theorem, we see that $P_0(z,f)$ has no non-zero roots and thus $P_0(z,f)^n=a_p$. Then by doing a linear transformation $f\to \alpha f$ with a suitably chosen rational function $\alpha$, we get
\begin{eqnarray}\label{equa    perf 0}
\overline{f}^n = cf^{-q},
\end{eqnarray}
where $c=\frac{a_p}{\overline{\alpha}^n\alpha^q}$ is a non-zero constant. Also, as in Theorem~\ref{Theorem  2a}, solutions of \eqref{equa    perf 0} are represented as
\begin{equation*}
f=c^{\frac{1}{n+q}}\exp[\pi(z)(-q/n)^z],
\end{equation*}
where $\pi(z)$ is an arbitrary entire periodic function with period~1. Otherwise, we have that $0$ is the only root of $Q(z,f)$ with $n\mid q$, or $0$ is not the only root of $Q(z,f)$. Recalling that we have excluded out the two possibilities that $0$ is a root of $P(z,f)$ of order $k_0$ such that $n\nmid k_0$ and that $0$ is a root of $Q(z,f)$ of order $l_0$ such that $l_0<q$ and $n\nmid l_0$, we see that in either case the only $\alpha_i$ in $P(z,f)$ or $\beta_j$ in $\eqref{Q}$ is non-zero. Therefore, we only need to consider the following two equations:
\begin{eqnarray}
\overline{f}^2 &=& \frac{P_0(z,f)^2}{Q_0(z,f)^2}(f-\alpha_1)^{k_1},
\label{eq  restrict c1}\\
\overline{f}^2 &=& \frac{P_0(z,f)^2}{Q_0(z,f)^2}\frac{1}{(f-\beta_1)^{l_1}}, \label{eq  restrict c3}
\end{eqnarray}
where $\alpha_1\not=0$ and $\beta_1\not=0$, $k_1,l_1\in \mathbb{N}$ are odd integers, and in equation \eqref{eq  restrict c1} we have $2p_0+k_1< 2q_0$ and in equation \eqref{eq  restrict c3} we have $2p_0< 2q_0+l_1$. Below, we discuss the two equations above separately.

\vskip 4pt

\noindent\textbf{Subcase~1:} Equation \eqref{eq  restrict c1}.

\vskip 4pt

From the previous discussions we see that $\pm \alpha_1$ and $\infty$ are all completely ramified rational functions of $f$. In fact, $\infty$ cannot be a Picard exceptional rational function of $f$; otherwise, the roots of $Q(z,f)$ are all Picard's exceptional rational functions of $f$, which is impossible by the inequality \eqref{multiplicityinequality}. Also, from \eqref{eq  restrict c1} we see that all roots of $f-\alpha_1=0$ with at most finitely many exceptions are of even multiplicities. If $Q_0(z,f)$ has a root, say $\beta$, of odd order, then by applying the same analysis as in the proof of Lemma~\ref{basiclemma} and considering the multiplicities of the roots of $f-\beta=0$ together with the fact that $\infty$ is a completely ramified rational function of $f$, we obtain that $\beta$ is a completely ramified rational function of $f$. Moreover, $\beta\not=0$ since otherwise from the proof of Lemma~\ref{basiclemma} we have that $\pm \alpha_1$ are completely ramified rational functions of $f$ with multiplicity~4, a contradiction to the inequality~\eqref{multiplicityinequality}. By Lemma~\ref{keylemma} it follows that $\pm \beta$ are both completely ramified rational functions of $f$. Then by Theorem~\ref{completelyrm} we must have $\beta=-\alpha_1$ and there is only one such $\beta$. We conclude that $Q_0(z,f)$ is of the form $Q_0(z,f)=Q_1(z,f)^2(f+\alpha_1)^{l_0}$ for some polynomial $Q_1(z,f)$ in $f$ and $l_0\in \mathbb{N}$ is $0$ or an odd integer. We may let $\alpha_1=-1$ by doing a linear transformation $f\to \alpha_1f$. We consider
\begin{equation}\label{eq  restrict c2tye000}
\overline{f}^2-1= \frac{P_0(z,f)^2(f-1)^{k_1}-Q_0(z,f)^2}{Q_0(z,f)^2}.
\end{equation}
Note that the numerator of the RHS of \eqref{eq  restrict c2tye000} is a polynomial in $f$ with degree $2q_0$. If it has two distinct roots, say $\gamma_1$ and $\gamma_2$, of odd orders, then by considering the multiplicities of the roots of $f-\gamma_1=0$ and $f-\gamma_2=0$, respectively, together with the fact that $\pm 1$ are both completely ramified rational functions of $f$ and that the roots of $f\pm1=0$ have even multiplicities with at most finitely many exceptions, we obtain that $\gamma_1$ and $\gamma_2$ are both completely ramified rational functions of $f$. By Lemma~\ref{keylemma}, it follows that $\pm \gamma_1$ and $\pm \gamma_2$ are all completely ramified rational functions of $f$; this yields a contradiction to Theorem~\ref{completelyrm} even when $\gamma_1+\gamma_2=0$. Therefore, the numerator of the RHS of \eqref{eq  restrict c2tye000} must be of the form $P_1(z,f)^2$ for some polynomial $P_1(z,f)$ in $f$. Now we have
\begin{equation}\label{eq  restrict c2tye0}
\overline{f}^2= \frac{P_0(z,f)^2(f-1)^{k_1}}{Q_0(z,f)^2},
\end{equation}
and further that
\begin{equation}\label{eq  restrict c2tye1}
\overline{f}^2-1 =\frac{P_1(z,f)^2}{Q_0(z,f)^2}.
\end{equation}
It follows that
\begin{equation*}
\begin{split}
P_0(z,f)^2(f-1)^{k_1}=P_1(z,f)^2+Q_0(z,f)^2=[P_1(z,f)+iQ_0(z,f)][P_1(z,f)-iQ_0(z,f)],
\end{split}
\end{equation*}
and so
\begin{equation}\label{eq  restrict c2tye1 copunpq1}
\begin{split}
P_1(z,f)+iQ_0(z,f)&=P_{01}(z,f),\\
P_1(z,f)-iQ_0(z,f)&=P_{02}(z,f),
\end{split}
\end{equation}
where $P_{01}(z,f)$ and $P_{02}(z,f)$ are two polynomials in $f$ such that $P_{01}(z,f)P_{02}(z,f)=P_0(z,f)^2(f-1)^{k_1}$. Since $P_{01}(z,f)$ and $P_{02}(z,f)$ have no common roots, without loss of generality, we may write
\begin{equation*}
\begin{split}
P_{01}(z,f)&=P_{011}(z,f)^2,\\
P_{02}(z,f)&=P_{012}(z,f)^2(f-1)^{k_1},
\end{split}
\end{equation*}
where $P_{011}(z,f)$ and $P_{012}(z,f)$ are two polynomials in $f$ with no common roots and $P_{011}(z,f)P_{012}(z,f)=P_0(z,f)$. From equations in \eqref{eq  restrict c2tye1 copunpq1} together with previous discussions we have
\begin{equation}\label{eq  restrict c2tye1 copunpq3}
\begin{split}
Q_0(z,f)&=\frac{1}{2i}[P_{011}(z,f)^2-P_{012}(z,f)^2(f-1)^{k_1}]=Q_1(z,f)^2(f+1)^{l_0},\\
P_1(z,f)&=\frac{1}{2}[P_{011}(z,f)^2+P_{012}(z,f)^2(f-1)^{k_1}].
\end{split}
\end{equation}
With the above two expressions for $Q_0(z,f)$ and $P_1(z,f)$, solutions of equation \eqref{eq  restrict c2tye0} can be obtained in the following way. Put
\begin{equation}\label{eq  restrict c2tye3beforeo}
f=\frac{1}{2}(\lambda+\lambda^{-1}).
\end{equation}
Note that the leading coefficient of the polynomial $P_1(z,f)^2$ is $-1$. From \eqref{eq  restrict c2tye1} we see that $\lambda$ is a meromorphic function and it follows that
\begin{equation}\label{eq  restrict c2tye3}
\frac{1}{2}\frac{\overline{\lambda}^2-1}{\overline{\lambda}}= \frac{P_1(z,f)}{Q_0(z,f)}.
\end{equation}
Since $\infty$ is a completely ramified rational function of $f$, then we see that $0$ and $\infty$ are both completely ramified rational functions of $\lambda$. Moreover, since all zeros of $f-1$ have multiplicities at least~2, from \eqref{eq  restrict c2tye0} we see that the leading coefficient of the polynomial $P_0(z,f)^2(f-1)^{k_1}$ is a square of some rational function and it follows that all poles of $f$ have even multiplicities. Put
\begin{equation}\label{eq  restrict c2tye3beforeo0}
\lambda=\delta^2.
\end{equation}
Then $\delta$ is a meromorphic function. By solving
equation \eqref{eq  restrict c2tye3} together with \eqref{eq  restrict c2tye0}, \eqref{eq  restrict c2tye1} and \eqref{eq  restrict c2tye3beforeo0}, we get
\begin{equation*}
\overline{\delta}^2= \frac{P_1(z,f)}{Q_0(z,f)}\pm\frac{P_0(z,f)(f-1)^{k_1/2}}{Q_0(z,f)}= \frac{P_1\left(z,\frac{\delta^4+1}{2\delta^{2}}\right)\pm P_0\left(z,\frac{\delta^4+1}{2\delta^{2}}\right)\left(\frac{\delta^2-1}{2^{1/2}\delta}\right)^{k_1}}{Q_0\left(z,\frac{\delta^4+1}{2\delta^{2}}\right)}.
\end{equation*}
Note that $P_0(z,f)=P_{011}(z,f)P_{012}(z,f)$. By combining the equation above and the equations in \eqref{eq  restrict c2tye1 copunpq3}, we have
\begin{equation}\label{Theorem  2b  eq0a2  fuja   proof5}
\overline{\delta}^2=i\frac{P_{011}\left(z,\frac{\delta^4+1}{2\delta^{2}}\right)+\theta P_{012}\left(z,\frac{\delta^4+1}{2\delta^{2}}\right)\left(\frac{\delta^2-1}{2^{1/2}\delta}\right)^{k_1}}{P_{011}\left(z,\frac{\delta^4+1}{2\delta^{2}}\right)-\theta P_{012}\left(z,\frac{\delta^4+1}{2\delta^{2}}\right)\left(\frac{\delta^2-1}{2^{1/2}\delta}\right)^{k_1}}, \quad \theta=\pm1.
\end{equation}
Denote the degrees of the polynomials $Q_1(z,f)$, $P_{011}(z,f)$ and $P_{012}(z,f)$ by $s_0$, $s_1$, $s_2$, respectively. By comparing the degrees of the polynomials in the first equation of \eqref{eq  restrict c2tye1 copunpq3} on both sides, we see that if $l_0=0$, then $s_0=s_1$ and if $l_0>0$, then $2s_0+l_0=2s_2+k_1$. We discuss these two cases below further.

When $l_0=0$, we have from the first equation in \eqref{eq  restrict c2tye1 copunpq3} that
\begin{equation*}
\begin{split}
[P_{011}(z,f)+(2i)^{1/2}Q_1(z,f)][P_{011}(z,f)-(2i)^{1/2}Q_1(z,f)]=P_{012}(z,f)^2(f-1)^{k_1}.
\end{split}
\end{equation*}
It follows that
\begin{equation}\label{Theorem  2b  eq0a2  fuja asw1}
\begin{split}
P_{011}(z,f)+(2i)^{1/2}Q_1(z,f)&=P_{013}(z,f),\\
P_{011}(z,f)-(2i)^{1/2}Q_1(z,f)&=P_{014}(z,f),
\end{split}
\end{equation}
where $P_{013}(z,f)$ and $P_{014}(z,f)$ are two polynomials in $f$ such that $P_{013}(z,f)P_{014}(z,f)=P_{012}(z,f)^2(f-1)^{k_1}$. Since $P_{013}(z,f)$ and $P_{014}(z,f)$ have no common roots, without loss of generality, we may write
\begin{equation}\label{eq  restrict c2tye1 copunpqaq00}
\begin{split}
P_{013}(z,f)&=P_{021}(z,f)^2,\\
P_{014}(z,f)&=P_{022}(z,f)^2(f-1)^{k_1},
\end{split}
\end{equation}
where $P_{021}(z,f)$ and $P_{022}(z,f)$ are two polynomials in $f$ with no common roots and $P_{021}(z,f)P_{022}(z,f)=P_{012}(z,f)$. Then we have from equations in \eqref{Theorem  2b  eq0a2  fuja asw1} that
\begin{equation}\label{Theorem  2b  eq0a2  fuja asw2}
\begin{split}
P_{011}(z,f)&=\frac{1}{2}[P_{013}(z,f)+P_{014}(z,f)].
\end{split}
\end{equation}
By combining equations in \eqref{eq  restrict c2tye1 copunpqaq00} and \eqref{Theorem  2b  eq0a2  fuja asw2}, we have from equation \eqref{Theorem  2b  eq0a2  fuja   proof5} that
\begin{equation}\label{Theorem  2b  eq0a2  fuja   proof5hh}
\overline{\delta}=\pm i^{1/2}\frac{P_{021}\left(z,\frac{\delta^4+1}{2\delta^{2}}\right)+\theta P_{022}\left(z,\frac{\delta^4+1}{2\delta^{2}}\right)\left(\frac{\delta^2-1}{2^{1/2}\delta}\right)^{k_1}}{P_{021}\left(z,\frac{\delta^4+1}{2\delta^{2}}\right)-\theta P_{022}\left(z,\frac{\delta^4+1}{2\delta^{2}}\right)\left(\frac{\delta^2-1}{2^{1/2}\delta}\right)^{k_1}}, \quad \theta=\pm1.
\end{equation}

When $l_0>0$, we let $h$ and $g$ be such that $h^2+1=f$ and $g^2-1=f$, respectively. Then we have from the first equation in \eqref{eq  restrict c2tye1 copunpq3} that
\begin{equation*}
\begin{split}
[P_{012}(z,f)h^{k_1}+(-2i)^{1/2}Q_1(z,f)g^{l_0}][P_{012}(z,f)h^{k_1}-(-2i)^{1/2}Q_1(z,f)g^{l_0}]=P_{011}(z,f)^2.
\end{split}
\end{equation*}
By \eqref{eq  restrict c2tye3beforeo} and \eqref{eq  restrict c2tye3beforeo0}, we may write $h$ and $g$ as $h=\frac{\delta^2-1}{2^{1/2}\delta}$ and $g=\frac{\delta^2+1}{2^{1/2}\delta}$, respectively, and it follows that
\begin{equation*}
\begin{split}
P_{012}(z,f)h^{k_1}&=P_{012}\left(z,\frac{\delta^4+1}{2\delta^{2}}\right)\left(\frac{\delta^2-1}{2^{1/2}\delta}\right)^{k_1}:=\frac{P_{0121}(z,\delta^2)}{(2^{1/2})^{2s_2+k_1}\delta^{2s_2+k_1}},\\
Q_1(z,f)g^{l_0}&=Q_1\left(z,\frac{\delta^4+1}{2\delta^{2}}\right)\left(\frac{\delta^2+1}{2^{1/2}\delta}\right)^{l_0}:=\frac{Q_{11}(z,\delta^2)}{(2^{1/2})^{2s_0+l_0}\delta^{2s_0+l_0}},
\end{split}
\end{equation*}
where $P_{0121}(z,\delta^2)$ and $Q_{11}(z,\delta^2)$ are two polynomials in $\delta^2$. Here, none of the roots of $P_{0121}(z,\delta^2)$ and $Q_{11}(z,\delta^2)$ in $\delta^2$ is equal to $\pm1$. Note that the leading coefficients of the two polynomials $P_{012}(z,f)^2(f-1)^{k_1}$ and $-2iQ_1(z,f)^2(f+1)^{l_0}$ are equal. Recalling that $2s_0+l_0=2s_2+k_1$, we see that one of the two polynomials $P_{0121}(z,\delta^2)+(-2i)^{1/2}Q_{11}(z,\delta^2)$ and $P_{0121}(z,\delta^2)-(-2i)^{1/2}Q_{11}(z,\delta^2)$ in $\delta^2$ has some zero roots. Consider the pair of equations
\begin{equation*}
\begin{split}
P_{012}(z,f)h^{k_1}+(-2i)^{1/2}Q_1(z,f)g^{l_0}&=0,\\
P_{012}(z,f)h^{k_1}-(-2i)^{1/2}Q_1(z,f)g^{l_0}&=0.
\end{split}
\end{equation*}
Since the two polynomials $P_{012}(z,f)$ and $Q_1(f)$ have no common roots, then together with the relations $h^2+1=f$ and $g^2-1=f$ we see that the above two equations with respect to $f$ cannot have any common solution and thus each root of the polynomial $P_{011}(z,f)^2$ satisfy only one of them. This implies the following two results: First, the two polynomials $P_{0121}(z,\delta^2)+(-2i)^{1/2}Q_{11}(z,\delta^2)$ and $P_{0121}(z,\delta^2)-(-2i)^{1/2}Q_{11}(z,\delta^2)$ in $\delta^2$ have no common non-zero roots; second, for each root $\gamma$ of $P_{011}(z,f)$, the two distinct non-zero solutions of the equation $\delta^4-2\gamma\delta^2+1=0$ with respect to $\delta^2$ must be either both roots of the polynomial $P_{0121}(z,\delta^2)+(-2i)^{1/2}Q_{11}(z,\delta^2)$ in $\delta^2$ or the polynomial $P_{0121}(z,\delta^2)-(-2i)^{1/2}Q_{11}(z,\delta^2)$ in $\delta^2$. Since $2s_0+l_0$ is an odd integer and since $P_{011}(z,f)=P_{011}\left(z,\frac{1}{2}(\delta^2+\delta^{-2})\right)$ is a rational function in $\delta^2$ whose denominator has only zero root, then from the above discussions we see that there must be an odd integer $t_1\in \mathbb{Z}\setminus\{0\}$ such that
\begin{equation}\label{Theorem  2b  eq0a2  fuja asw4}
\begin{split}
P_{012}(z,f)h^{k_1}+(-2i)^{1/2}Q_1(z,f)g^{l_0}&=P_{015}(z,f)\delta^{t_1},\\
P_{012}(z,f)h^{k_1}-(-2i)^{1/2}Q_1(z,f)g^{l_0}&=P_{016}(z,f)\delta^{-t_1},
\end{split}
\end{equation}
where $P_{015}(z,f)$ and $P_{016}(z,f)$ are two polynomials
in $f$ with no common roots and $P_{015}(z,f)P_{016}(z,f)=P_{011}(z,f)^2$. We may write
\begin{equation}\label{Theorem  2b  eq0a2  fuja asw5}
\begin{split}
P_{015}(z,f)&=P_{023}(z,f)^2,\\
P_{016}(z,f)&=P_{024}(z,f)^2,
\end{split}
\end{equation}
where $P_{023}(z,f)$ and $P_{024}(z,f)$ are two polynomials in $f$ with no common roots
and $P_{023}(z,f)P_{024}(z,f)=P_{011}(z,f)$. Then we have from equations in \eqref{Theorem  2b  eq0a2  fuja asw4} that
\begin{equation}\label{Theorem  2b  eq0a2  fuja asw6}
\begin{split}
P_{012}(z,f)h^{k_1}&=\frac{1}{2}[P_{015}(z,f)\delta^{t_1}+P_{016}(z,f)\delta^{-t_1}].
\end{split}
\end{equation}
By combining equations \eqref{Theorem  2b  eq0a2  fuja asw5} and \eqref{Theorem  2b  eq0a2  fuja asw6}, we have from \eqref{Theorem  2b  eq0a2  fuja   proof5} that
\begin{equation}\label{Theorem  2b  eq0a2  fuja asw7}
\overline{\delta}=\pm (-i)^{1/2}\frac{P_{023}\left(z,\frac{\delta^4+1}{2\delta^{2}}\right)\delta^{t_1}+\theta P_{024}\left(z,\frac{\delta^4+1}{2\delta^{2}}\right)}{P_{023}\left(z,\frac{\delta^4+1}{2\delta^{2}}\right)\delta^{t_1}-\theta P_{024}\left(z,\frac{\delta^4+1}{2\delta^{2}}\right)}, \quad \theta=\pm1.
\end{equation}
We conclude that solutions of \eqref{eq  restrict c2tye0} are represented by \eqref{eq  restrict c2tye3beforeo} and \eqref{eq  restrict c2tye3beforeo0}, i.e., $f=\frac{1}{2}(\delta^2+\delta^{-2})$ with $\delta$ being a solution of \eqref{Theorem  2b  eq0a2  fuja   proof5hh} or \eqref{Theorem  2b  eq0a2  fuja asw7}.

\vskip 4pt

\noindent\textbf{Subcase~2:} Equation \eqref{eq  restrict c3}.

\vskip 4pt

Since $\beta_1\not=0$, then by similar arguments as in the previous case we know that $\pm\beta_1$ and $\infty$ are all completely ramified rational functions of $f$. We may let $\beta_1=1$ by doing a linear transformation $f\to\beta_1f$. Further, by considering the multiplicities of the roots of $f-1=0$ together with the fact that $\infty$ has multiplicity at least~$2$ and by Lemma~\ref{keylemma}, we obtain from \eqref{eq  restrict c3} that $\pm 1$ both have multiplicities~$4$ and it follows that $\Theta(\infty,f)+\Theta(1,f)+\Theta(-1,f)=2$. Then by applying the analysis in the proof of Lemma~\ref{keylemma0} to the roots of $P_0(z,f)$ and poles of $f$ and comparing the multiplicities of the zeros on both sides of \eqref{eq  restrict c3} and to the roots of $Q_0(z,f)^2(f-1)^{l_1}$ and comparing the multiplicities of the poles on both sides of \eqref{eq  restrict c3}, respectively, we obtain that $l_1=1$, $p_0=q_0$, $P_0(z,f)$ has simple roots only, $Q_0(z,f)=Q_1(z,f)^2$ for some polynomial $Q_1(z,f)$ in $f$ with simple roots only and none of these roots equals $\pm1$. We consider
\begin{equation}\label{eq  restrict c3put0hg}
\overline{f}^2-1 = \frac{P_0(z,f)^2-Q_0(z,f)^2(f-1)}{Q_0(z,f)^2(f-1)}.
\end{equation}
The numerator of the RHS of \eqref{eq  restrict c3put0hg} is a polynomial in $f$ with degree $2q_0+1$ and thus has at least one root, say $\gamma_1$, of odd order. By applying the same analysis as in the proof of Lemma~\ref{basiclemma} and considering the multiplicities of the roots of $f-\gamma_1=0$ together with the fact that both $\pm 1$ have multiplicities~4, we obtain that $\gamma_1$ is also a completely ramified rational function of $f$. By Lemma~\ref{keylemma} and Theorem~\ref{completelyrm} we must have $\gamma_1=-1$.
Then by Lemma~\ref{keylemma1} we conclude that the order of $\gamma_1$ equals~1. Also, by Theorem~\ref{completelyrm} we see that there is only one such $\gamma_1$ that has odd order. Now we have
\begin{equation}\label{eq  restrict c3put0}
\overline{f}^2= \frac{P_0(z,f)^2}{Q_0(z,f)^2(f-1)},
\end{equation}
and also that
\begin{equation}\label{eq  restrict c3put1}
\overline{f}^2-1 = \frac{P_1(z,f)^2(f+1)}{Q_0(z,f)^2(f-1)}.
\end{equation}
Moreover, by Lemma~\ref{keylemma1} it follows that $P_1(z,f)=P_2(z,f)^2$ for some polynomials $P_2(z,f)$ with simple roots only. We let $(f+1)/(f-1)=g^2$. Then $g$ is an algebroid function with at most finitely many algebraic branch points and it follows that the RHS of \eqref{eq  restrict c3put1} becomes $[P_1(z,f)g/Q_0(z,f)]^2$. Put
\begin{equation}\label{eq  restrict c2tye3beforeo1}
f=\frac{1}{2}(\lambda+\lambda^{-1}).
\end{equation}
With the setting $(f+1)/(f-1)=g^2$ we see from \eqref{eq  restrict c3put1} that $\lambda$ is an algebroid function with at most finitely many algebraic branch points. Recall that $\pm 1$ both have multiplicities~$4$ and $\infty$ has multiplicity~$2$. By a simple multiplicity analysis on the zeros, poles and $\pm1$-points of $\lambda$, we obtain from the definition in \eqref{eq  restrict c2tye3beforeo1} that $0$, $\infty$ and $\pm1$ are all completely ramified rational functions of $\lambda$ with multiplicities~$2$.
By substituting \eqref{eq  restrict c2tye3beforeo1} into \eqref{eq  restrict c3put1} we get
\begin{equation}\label{eq  restrict c3put2}
\frac{1}{2}\frac{\overline{\lambda}^2-1}{\overline{\lambda}}= \frac{P_1(z,f)}{Q_0(z,f)}\frac{\lambda+1}{\lambda-1}.
\end{equation}
We put
\begin{equation}\label{eq  restrict c2tye3beforeo2}
\lambda=\delta^2.
\end{equation}
Then $\delta$ is an algebroid function with at most finitely many algebraic branch points. Moreover, by the definition of $\lambda$ we see that $\pm1$ and $\pm i$ are all completely ramified rational functions of $\delta$ with multiplicity~$2$. By solving equation \eqref{eq  restrict c3put2} together with equations \eqref{eq  restrict c3put0}, \eqref{eq  restrict c3put1} and \eqref{eq  restrict c2tye3beforeo2}, we get
\begin{equation}\label{eq  restrict c3put2frty0be}
\overline{\delta}^2= \frac{P_1(z,f)}{Q_0(z,f)}\frac{\lambda+1}{\lambda-1}\pm\frac{P_0(z,f)}{Q_0(z,f)(f-1)^{1/2}}= \frac{P_1\left(z,\frac{\delta^4+1}{2\delta^2}\right)(\delta^2+1)\pm P_0\left(z,\frac{\delta^4+1}{2\delta^2}\right)(2^{1/2}\delta)}{Q_0\left(z,\frac{\delta^4+1}{2\delta^2}\right)(\delta^2-1)}.
\end{equation}
Recall that $p_0=q_0$ from the beginning of this subcase. Also, from previous discussions we see that the degree of the polynomial $P_1(z,f)$ in \eqref{eq  restrict c3put1}, denoted by $p_1$, satisfies $p_1=p_0$. By multiplying $(2\delta^2)^{p_0}$ to both of the numerator and the denominator of the RHS of \eqref{eq  restrict c3put2frty0be}, we have
\begin{equation}\label{eq  restrict c3put2frty0}
\overline{\delta}^2= \frac{P_{10}(z,\delta^2)(\delta^2+1)\pm P_{00}(z,\delta^2)(2^{1/2}\delta)}{Q_{00}(z,\delta^2)(\delta^2-1)},
\end{equation}
where $P_{10}(z,\delta^2)$, $P_{00}(z,\delta^2)$ and $Q_{00}(z,\delta^2)$ are polynomials in $\delta$ of the same degrees $4q_0$. Note that $\delta$ cannot have any other completely ramified rational function besides $\pm1$ and $\pm i$. Since $\pm1$ are not roots of $P_0(z,f)$, $Q_0(z,f)$ or $P_1(z,f)$, then from the above reasoning we see that $\pm1$ and $\pm i$ are not roots of the polynomials $P_{10}(z,\delta^2)$, $P_{00}(z,\delta^2)$ or $Q_{00}(z,\delta^2)$. By applying Lemma~\ref{keylemma0} to equation \eqref{eq  restrict c3put2frty0}, we conclude that the numerator of the RHS of \eqref{eq  restrict c3put2frty0} must be a square of some polynomial in $\delta$ with simple roots only and none of these roots is equal to $\pm 1$ or $\pm i$. Moreover, since $p<q$ and we have assumed $b_q=1$, we see from the above reasoning that the leading coefficient of the numerator of the RHS of \eqref{eq  restrict c3put2frty0} is $\pm i$. Therefore, we may write
\begin{equation}\label{eq  restrict c3put2frty1}
P_{10}(z,\delta^2)^2(\delta^2+1)\pm P_{00}(z,\delta^2)2^{1/2}\delta=\pm i(\delta-\epsilon_1)^2\cdots (\delta-\epsilon_{2q_{0}+1})^2:=P_{11}(z,\delta),
\end{equation}
where the roots $\epsilon_1$, $\cdots$, $\epsilon_{2q_{0}+1}$ are in general algebraic functions, distinct from each other and from $\pm1$ and $\pm i$. Since $P_{11}(z,\delta)$ is a polynomial in $\delta$ of degree $4q_0+2$, we may denote by $\varrho_{4q_{0}+2}$, $\cdots$, $\varrho_1$ the roots of $P_{11}(z,\delta)$. Then from the equations in \eqref{eq  restrict c3put2frty0be} and \eqref{eq  restrict c3put2frty0} we see that $0$ is not a root of the polynomial $P_{11}(z,\delta)$ and that $\varrho_{4q_{0}+2}^{-1}$, $\cdots$, $\varrho_1^{-1}$ are also roots of $P_{11}(z,\delta)$. However, since $\epsilon_1$, $\cdots$, $\epsilon_{2q_{0}+1}$ are distinct from each other, we see that $\varrho_i=\varrho_i^{-1}$ for at least one $i\in\{4q_0+2,\cdots,1\}$ and thus $\varrho_i=\pm 1$, a contradiction to our previous observations. This implies that equation \eqref{eq  restrict c3put0} cannot have any meromorphic solution, which completes the proof.

\end{proof}

For the case $q\geq 1$ and $2\leq n<d$ and $n\mid |p-q|$ of equation \eqref{yanagiharaeq0}, we prove the following Theorem~\ref{Theorem  2c}. In this theorem, we will do a bilinear transformation to $f$ in \eqref{yanagiharaeq0}; for simplicity, we always assume that the leading coefficient of the resulting $Q_0(z,f)^n$ is~1 and denote the leading coefficient of the resulting $P_0(z,f)^n$ by $A$.

\begin{theorem}\label{Theorem  2c}
Suppose that $q\geq 1$ and $2\leq n<d$ and $n\mid |p-q|$. Let $f$ be a transcendental meromorphic solution of equation \eqref{yanagiharaeq0}. Then there exists a rational or algebraic function $\alpha$ such that by doing a linear transformation $f\rightarrow \alpha f$ or $f\rightarrow 1/(\alpha f)$, \eqref{yanagiharaeq0} becomes \eqref{yanagiharaeq0    9}, \eqref{yanagiharaeq0    5} or \eqref{Theorem  2b  eq0}, or one of the following equations:

\begin{enumerate}
\item [(1)] the first equation is
\begin{equation}\label{Theorem  2c  eq12 first0}
\overline{f}^2= \frac{P_0(z,f)^2}{Q_0(z,f)^2}(f-1)(f-\kappa),
\end{equation}
where $\kappa\not=0,\pm1$, $P_0(z,f)$ and $Q_0(z,f)$ are two polynomials in $f$ with simple roots only and $p_0-q_0\in\{-2,-1,0\}$; moreover,
\begin{equation}\label{Theorem  2c  eq12 first1}
\begin{split}
P_0(z,f)^2(f-1)(f-\kappa)-Q_0(z,f)^2&=P_1(z,f)^2(f+1)(f+\kappa),\\
P_0(z,f)^2(f-1)(f-\kappa)-\overline{\kappa}^2 Q_0(z,f)^2&=P_2(z,f)^2,
\end{split}
\end{equation}
where $P_1(z,f)$ and $P_2(z,f)$ are two polynomials in $f$ with simple roots only and of degrees $p_1$ and $p_2$ such that $2p_1+2\in\{d,d-2\}$ and $2p_2\in\{d,d-2\}$;

\item [(2)] the second equation is
\begin{equation}\label{Theorem  2c  eq12}
\overline{f}^2= \frac{P_0(z,f)^2}{Q_0(z,f)^2}(f-1)^{k_1}(f+1)^{k_2},
\end{equation}
where $k_1,k_2$ are two positive odd integers such that $P_0(z,f)^2(f-1)^{k_1}(f+1)^{k_2}-Q_0(z,f)^2=P_1(z,f)^2$ for some polynomial $P_1(z,f)$ in $f$;

\item [(3)] the third equation is
\begin{equation}\label{Theorem  2c  eq12 third0}
\overline{f}^2= \frac{P_0(z,f)^2}{Q_0(z,f)^2}(f^2-1),
\end{equation}
where $P_0(z,f)$ and $Q_0(z,f)$ are two polynomials in $f$ with simple roots only such that $p_0-q_0\in\{-2,-1,0\}$; moreover,
\begin{equation}\label{Theorem  2c  eq12 third1}
\begin{split}
P_0(z,f)^2(f^2-1)-Q_0(z,f)^2&=P_1(z,f)^2(f^2-\gamma^2),\\
P_0(z,f)^2(f^2-1)-\overline{\gamma}^2 Q_0(z,f)^2&=P_2(z,f)^2,
\end{split}
\end{equation}
where $\gamma\not=0,\pm1$, $P_1(z,f)$ and $P_2(z,f)$ are two polynomials in $f$ with simple roots only and of degrees $p_1$ and $p_2$ such that $2p_1+2\in\{d,d-2\}$ and $2p_2\in\{d,d-2\}$;

\item [(4)] the fourth equation is
\begin{equation}\label{Theorem  2c  eq12 fourth0}
\overline{f}^2= \frac{P_0(z,f)^2(f-\kappa)}{Q_0(z,f)^2(f-1)},
\end{equation}
where $\kappa\not=0,1$, $P_0(z,f)$ and $Q_0(z,f)$ are two polynomials in $f$ with simple roots only and $p_0-q_0\in\{-1,0,1\}$; moreover, when $\kappa=-1$, we have
\begin{equation}\label{Theorem  2c  eq12 fourth1}
\begin{split}
P_0(z,f)^2(f+1)-Q_0(z,f)^2(f-1)&=P_1(z,f)^2(f-\gamma),\\
P_0(z,f)^2(f-1)-\overline{\gamma}^2 Q_0(z,f)^2(f-1)&=P_2(z,f)^2(f+\gamma),
\end{split}
\end{equation}
where $\gamma\not=0,\pm1$, $P_1(z,f)$ and $P_2(z,f)$ are two polynomials in $f$ with simple roots only and of degrees $p_1$ and $p_2$ such that $2p_1+1\in\{d,d-2\}$ and $2p_2+1\in\{d,d-2\}$; or when $\kappa\not=-1$, we have
\begin{equation}\label{Theorem  2c  eq12 fourth2}
\begin{split}
P_0(z,f)^2(f-\kappa)-Q_0(z,f)^2(f-1)&=P_1(z,f)^2(f+\kappa),\\
P_0(z,f)^2(f-\kappa)-\overline{\kappa}^2 Q_0(z,f)^2(f-1)&=P_2(z,f)^2(f+1),
\end{split}
\end{equation}
where $P_1(z,f)$ and $P_2(z,f)$ are two polynomials in $f$ with simple roots only and of degrees $p_1$ and $p_2$ such that $2p_1+1\in\{d,d-2\}$ and $2p_2+1\in\{d,d-2\}$;

\item [(5)] the fifth equation is
\begin{equation}\label{Theorem  2c  eq17}
\overline{f}^3= \frac{P_0(z,f)^3(f-1)}{Q_0(z,f)^3(f-\eta)},
\end{equation}
where $\eta$ is a cubic root of $1$ such that $\eta^2+\eta+1=0$, $P_0(z,f)$ and $Q_0(z,f)$ are two polynomials in $f$ with simple roots only and $p_0-q_0\in\{-1,0,1\}$, and $P_0(z,f)^3(f-1)-Q_0(z,f)^3(f-\eta)=P_1(z,f)^3(f-\eta^2)$ for some polynomial $P_1(z,f)$ in $f$ with simple roots only and of degree $p_1$ such that $3p_1+1\in\{d,d-3\}$;

\item [(6)] the sixth equation is
\begin{equation}\label{Theorem  2c  eq0}
\overline{f}^3= \frac{P_0(z,f)^3}{Q_0(z,f)^3}(f^3-1),
\end{equation}
where $P_0(z,f)$ and $Q_0(z,f)$ are two polynomials in $f$ with simple roots only such that $p_0-q_0\in\{-2,-1,0\}$, and $P_0(z,f)^3(f^3-1)-Q_0(z,f)^3=P_1(z,f)^3$ for some polynomial $P_1(z,f)$ in $f$ with simple roots only and of degree $p_1$ such that $3p_1\in\{d,d-3\}$;

\item [(7)] the seventh equation is
\begin{equation}\label{Theorem  2c  eq12 seventh0}
\overline{f}^2= \frac{P_0(z,f)^2(f^2-\kappa^2)}{Q_0(z,f)^2(f^2-1)},
\end{equation}
where $\kappa\not=0,\pm1$, $P_0(z,f)$ and $Q_0(z,f)$ are two polynomials in $f$ with simple roots only such that $p_0-q_0\in\{-1,0,1\}$; moreover,
\begin{equation}\label{Theorem  2c  eq12 seventh1}
\begin{split}
P_0(z,f)^2(f^2-\kappa^2)-Q_0(z,f)^2(f^2-1)&=P_1(z,f)^2,\\
P_0(z,f)^2(f^2-\kappa^2)-\overline{\kappa}^2Q_0(z,f)^2(f^2-1)&=P_2(z,f)^2,
\end{split}
\end{equation}
where $P_1(z,f)$ and $P_2(z,f)$ are two polynomials in $f$ with simple roots only and of degrees $p_1$ and $p_2$ such that $2p_1\in\{d,d-2\}$ and $2p_2\in\{d,d-2\}$;

\item [(8)] the eighth equation is
\begin{equation}\label{Theorem  2c  eq12 eightth0}
\overline{f}^2= \frac{P_0(z,f)^2(f-\kappa)(f-1)}{Q_0(z,f)^2(f+\kappa)(f+1)},
\end{equation}
where $\kappa\not=0,\pm1$, $P_0(z,f)$ and $Q_0(z,f)$ are two polynomials in $f$ with simple roots only such that $p_0-q_0\in\{-1,0,1\}$; moreover,
\begin{equation}\label{Theorem  2c  eq12 eightth1}
\begin{split}
P_0(z,f)^2(f-\kappa)(f-1)-Q_0(z,f)^2(f+\kappa)(f+1)&=P_1(z,f)^2,\\
P_0(z,f)^2(f-\kappa)(f-1)-\overline{\kappa}^2Q_0(z,f)^2(f+\kappa)(f+1)&=P_2(z,f)^2,
\end{split}
\end{equation}
where $P_1(z,f)$ and $P_2(z,f)$ are two polynomials in $f$ with simple roots only and of degrees $p_1$ and $p_2$ such that $2p_1\in\{d,d-2\}$ and $2p_2\in\{d,d-2\}$;

\item [(9)] the ninth equation is
\begin{equation}\label{Theorem  2c  eq12 ninth0}
\overline{f}^2= \frac{P_0(z,f)^2}{Q_0(z,f)^2}(f^2-\kappa^2)(f^2-1),
\end{equation}
where $\kappa\not=0,\pm1$, $P_0(z,f)$ and $Q_0(z,f)$ are two polynomials in $f$ with simple roots only such that $p_0-q_0\in\{-3,-2,-1\}$; moreover,
\begin{equation}\label{Theorem  2c  eq12 ninth1}
\begin{split}
P_0(z,f)^2(f^2-\kappa^2)(f^2-1)-Q_0(z,f)^2&=P_1(z,f)^2,\\
P_0(z,f)^2(f^2-\kappa^2)(f^2-1)-\overline{\kappa}^2Q_0(z,f)^2&=P_2(z,f)^2,
\end{split}
\end{equation}
where $P_1(z,f)$ and $P_2(z,f)$ are two polynomials in $f$ with simple roots only and of degrees $p_1$ and $p_2$ such that $2p_1\in\{d,d-2\}$ and $2p_2\in\{d,d-2\}$;

\end{enumerate}

\end{theorem}

Before getting into the proof, we make some remarks on how to solve the nine equations in Theorem~\ref{Theorem  2c}. Unlike in the case of equation \eqref{yanagiharaeq0} with $n=d$, it seems impossible to give explicit expressions for the two polynomials $P_0(z,f)$ and $Q_0(z,f)$ in the general case. However, when the degrees $p_0$ and $q_0$ are given, it is possible to do this via some basic computations.

First, we show how to solve the second equation \eqref{Theorem  2c  eq12}. Since $P_0(z,f)^2(f-1)^{k_1}(f+1)^{k_2}-Q_0(z,f)^2=P_1(z,f)^2$ for some polynomial $P_1(z,f)$ in $f$, we have the equation
\begin{equation}\label{eq  restrict a08 e  p eq2s1 t}
\overline{f}^2-1= \frac{P_1(z,f)^2}{Q_0(z,f)^2},
\end{equation}
and may write
\begin{equation}\label{eq  restrict a08 e  p eq2s2 t purt1alpot1xsa1}
\begin{split}
P_1(z,f)+iQ_0(z,f)&=P_{01}(z,f),\\
P_1(z,f)-iQ_0(z,f)&=P_{02}(z,f),
\end{split}
\end{equation}
where $P_{01}(z,f)$ and $P_{02}(z,f)$ are two polynomials in $f$ such that $P_{01}(z,f)P_{02}(z,f)=P_0(z,f)^2(f-1)^{k_1}(f+1)^{k_2}$. Here, $P_{01}(z,f)$ and $P_{02}(z,f)$ have no common roots. Without loss of generality, we may write
\begin{equation}\label{eq  restrict a08 e  p eq2s2 t purt1alpot1xsa1hgtbbbv0}
\begin{split}
P_{01}(z,f)&=P_{011}(z,f)^2(f-1)^{k_1}(f+1)^{k_2},\\
P_{02}(z,f)&=P_{012}(z,f)^2,
\end{split}
\end{equation}
or
\begin{equation}\label{eq  restrict a08 e  p eq2s2 t purt1alpot1xsa1hgtbbbv1}
\begin{split}
P_{01}(z,f)&=P_{011}(z,f)^2(f-1)^{k_1},\\
P_{02}(z,f)&=P_{012}(z,f)^2(f+1)^{k_2},
\end{split}
\end{equation}
where $P_{011}(z,f)$ and $P_{012}(z,f)$ are two polynomials in $f$ such that $P_0(z,f)=P_{011}(z,f)P_{012}(z,f)$. Now, $P_{011}(z,f)$ and $P_{012}(z,f)$ have no common roots. From equations in \eqref{eq  restrict a08 e  p eq2s2 t purt1alpot1xsa1} we get
\begin{equation}\label{eq  restrict a08 e  p eq2s2 t purt1alpot1xsa1tere}
\begin{split}
Q_0(z,f)&=\frac{1}{2i}[P_{01}(z,f)-P_{02}(z,f)],\\
P_1(z,f)&=\frac{1}{2}[P_{01}(z,f)+P_{02}(z,f)].
\end{split}
\end{equation}
Put $f=\frac{1}{2}(\lambda+\lambda^{-1})$. Then from \eqref{eq  restrict a08 e  p eq2s1 t} we see that $\lambda$ is an algebroid function with at most finitely many algebraic branch points and it follows from \eqref{eq  restrict a08 e  p eq2s1 t} that
\begin{equation}\label{eq  restrict a08 e  p eq6}
\frac{1}{2}(\overline{\lambda}-\overline{\lambda}^{-1})=\frac{P_1(z,f)}{Q_0(z,f)}.
\end{equation}
By solving equation \eqref{eq  restrict a08 e  p eq6} together with equation \eqref{eq  restrict a08 e  p eq2s1 t} and equations in \eqref{eq  restrict a08 e  p eq2s2 t purt1alpot1xsa1tere}, we get
\begin{equation*}
\overline{\lambda}=\frac{\frac{1}{2}\left[P_{01}\left(z,\frac{\lambda^2+1}{2\lambda}\right)+P_{02}\left(z,\frac{\lambda^2+1}{2\lambda}\right)\right]\pm P_0\left(z,\frac{\lambda^2+1}{2\lambda}\right)\frac{(\lambda-1)^{k_1}(\lambda+1)^{k_2}}{(2\lambda)^{(k_1+k_2)/2}}}{\frac{1}{2i}\left[P_{01}\left(z,\frac{\lambda^2+1}{2\lambda}\right)-P_{02}\left(z,\frac{\lambda^2+1}{2\lambda}\right)\right]}.
\end{equation*}
Therefore, when equations in \eqref{eq  restrict a08 e  p eq2s2 t purt1alpot1xsa1hgtbbbv0} hold, we have
\begin{equation}\label{eq  restrict a08 e  p eq6aftop0}
\overline{\lambda}=i\frac{P_{011}\left(z,\frac{\lambda^2+1}{2\lambda}\right)\frac{(\lambda-1)^{k_1}(\lambda+1)^{k_2}}{(2\lambda)^{(k_1+k_2)/2}}+\theta P_{012}\left(z,\frac{\lambda^2+1}{2\lambda}\right)}{P_{011}\left(z,\frac{\lambda^2+1}{2\lambda}\right)\frac{(\lambda-1)^{k_1}(\lambda+1)^{k_2}}{(2\lambda)^{(k_1+k_2)/2}}-\theta P_{012}\left(z,\frac{\lambda^2+1}{2\lambda}\right)}, \quad \theta=\pm1,
\end{equation}
or, when equations in \eqref{eq  restrict a08 e  p eq2s2 t purt1alpot1xsa1hgtbbbv1} hold, we have
\begin{equation}\label{eq  restrict a08 e  p eq6aftop1}
\overline{\lambda}=i\frac{P_{011}\left(z,\frac{\lambda^2+1}{2\lambda}\right)\frac{(\lambda-1)^{k_1}}{(2\lambda)^{k_1/2}}+\theta P_{012}\left(z,\frac{\lambda^2+1}{2\lambda}\right)\frac{(\lambda+1)^{k_2}}{(2\lambda)^{k_2/2}}}{P_{011}\left(z,\frac{\lambda^2+1}{2\lambda}\right)\frac{(\lambda-1)^{k_1}}{(2\lambda)^{k_1/2}}-\theta P_{012}\left(z,\frac{\lambda^2+1}{2\lambda}\right)\frac{(\lambda+1)^{k_2}}{(2\lambda)^{k_2/2}}}, \quad \theta=\pm1.
\end{equation}
We conclude that solutions of equation \eqref{eq  restrict a08 e  p eq2s1 t} are represented by $f=\frac{1}{2}(\lambda+\lambda^{-1})$ with $\lambda$ being a solution of equation \eqref{eq  restrict a08 e  p eq6aftop0} or \eqref{eq  restrict a08 e  p eq6aftop1}. In the autonomous case, the RHS of equation \eqref{eq  restrict a08 e  p eq6aftop0} or \eqref{eq  restrict a08 e  p eq6aftop1} becomes a rational term $R_1(\lambda)$ in $\lambda$ after multiplying $(2\lambda)^{d}$ to both of the numerator and the denominator; thus the equation $\overline{\lambda}=R_1(\lambda)$ always has a meromorphic solution for any given $P_{01}(f)$ and $P_{02}(f)$ in \eqref{eq  restrict a08 e  p eq2s2 t purt1alpot1xsa1hgtbbbv0} or in \eqref{eq  restrict a08 e  p eq2s2 t purt1alpot1xsa1hgtbbbv1} and any given $Q_0(f)$ such that \eqref{eq  restrict a08 e  p eq2s2 t purt1alpot1xsa1tere} holds, as mentioned in the introduction.

Second, we show how to determine the polynomials $P_0(z,f)$ and $Q_0(z,f)$ from the seven pairs of polynomial equations in \eqref{Theorem  2c  eq12 first1}, \eqref{Theorem  2c  eq12 third1}, \eqref{Theorem  2c  eq12 fourth1}, \eqref{Theorem  2c  eq12 fourth2}, \eqref{Theorem  2c  eq12 seventh1}, \eqref{Theorem  2c  eq12 eightth1} and \eqref{Theorem  2c  eq12 ninth1}. Take the pair of equations in \eqref{Theorem  2c  eq12 first1} as an example. When $p_0=q_0-1$, we write
\begin{equation}\label{Theorem  2c  eq12 first1 solv1}
\begin{split}
P_0(z,f)^2&=A[(f-a_{p_0})(f-a_{p_0-1})\cdots(f-a_1)]^2,\\
Q_0(z,f)^2&=[(f-b_{q_0})(f-b_{q_0-1})\cdots(f-b_1)]^2,
\end{split}
\end{equation}
where $a_i$ and $b_j$ are in general algebraic functions, distinct from each other. If $A\not=1,\overline{\kappa}^2$, we also write
\begin{equation}\label{Theorem  2c  eq12 first1 solv2}
\begin{split}
P_1(z,f)^2&=(A-1)[(f-c_{p_0})(f-c_{p_0-1})\cdots(f-c_1)]^2,\\
P_2(z,f)^2&=(A-\overline{\kappa}^2)[(f-d_{q_0})(f-d_{q_0-1})\cdots(f-d_1)]^2,
\end{split}
\end{equation}
where $c_k$ and $d_l$ are in general algebraic functions, distinct from each other. By comparing the coefficients on both sides of the two polynomial equations in \eqref{Theorem  2c  eq12 first1} we obtain $4p_0+4$ polynomial equations with respect to the unknowns $A$, $\kappa$, $a_i$, $b_j$, $c_k$ and $d_l$, whose combined number is $4p_0+4$. This implies that $\overline{\kappa}$ and $\kappa$ satisfy a polynomial equation $U(\overline{\kappa},\kappa)=0$ with respect to $\overline{\kappa}$ and $\kappa$. In general, it is difficult to determine whether or not $\kappa$ is a constant, but this can be done once the degrees $p_0$ and $q_0$ are given. The case when $A=1$ or $A=\overline{\kappa}^2$ is similar. It seems that the pair of polynomial equations in \eqref{Theorem  2c  eq12 first1} is solvable only when $p_0=q_0-1$ since we would obtain $4p_0+4$ polynomial equations with respect to $4p_0+3$ unknowns when $p_0=q_0$ and we would obtain $4p_0+8$ polynomial equations with respect to $4p_0+7$ unknowns when $p_0=q_0-2$. The polynomial equations in \eqref{Theorem  2c  eq12 third1}, \eqref{Theorem  2c  eq12 fourth1}, \eqref{Theorem  2c  eq12 fourth2}, \eqref{Theorem  2c  eq17}, \eqref{Theorem  2c  eq12 seventh1}, \eqref{Theorem  2c  eq12 eightth1} and \eqref{Theorem  2c  eq12 ninth1} can be given similar discussions. Also, for the three pairs of polynomial equations in \eqref{Theorem  2c  eq12 fourth1}, \eqref{Theorem  2c  eq12 fourth2} and \eqref{Theorem  2c  eq12 eightth1}, it seems that they are solvable only when $p_0=q_0$. For given small $p_0$ and $q_0$, one may readily obtain some explicit examples for each of the seven pairs of polynomial equations by doing some basic computations. In particular, in all the three cases $p=q$, $p=q+2$ and $p=q-2$ one may obtain some examples for each of the three equations \eqref{Theorem  2c  eq12 third0}, \eqref{Theorem  2c  eq12 eightth0} and \eqref{Theorem  2c  eq12 ninth0}, where $P_0(z,f)$ and $Q_0(z,f)$ are in fact both polynomials in $f^2$. It is not known if this is always true for any given $p_0$ and $q_0$. For the fifth equation \eqref{Theorem  2c  eq17}, we note that the polynomial equation $P_0(z,f)^3(f-1)-Q_0(z,f)^3(f-\eta)=P_1(z,f)^3(f-\eta^2)$ is solvable when $p_0=q_0$ by the same arguments as above. In fact, we may even show that all coefficients of $P_0(z,f)^3(f^3-1)$ and $Q_0(z,f)^3$ are constants since we do not need to deal with the shift of a function.

Third, for the sixth equation \eqref{Theorem  2c  eq0}, below we show that the polynomial equation $P_0(z,f)^3(f^3-1)-Q_0(z,f)^3=P_1(z,f)^3$ is also solvable. We may write
\begin{equation*}
[P_1(z,f)+Q_0(z,f)][P_1(z,f)+\eta Q_0(z,f)][P_1(z,f)+\eta^2Q_0(z,f)]=P_0(z,f)^3(f^3-1).
\end{equation*}
Since any two of the polynomials $P_1(z,f)+Q_0(z,f)$, $P_1(z,f)+\eta Q_0(z,f)$ and $P_1(z,f)+\eta^2Q_0(z,f)$ have no common roots, we may write
    \begin{equation}\label{eq  restrict a01 e  p eq2    juy1}
    \begin{split}
    P_1(z,f)+Q_0(z,f)&=P_{01}(z,f)^3(f-1)^{\theta_{11}}(f-\eta)^{\theta_{12}}(f-\eta^2)^{\theta_{13}},\\
    P_1(z,f)+\eta Q_0(z,f)&=P_{02}(z,f)^3(f-1)^{\theta_{21}}(f-\eta)^{\theta_{22}}(f-\eta^2)^{\theta_{23}},\\
    P_1(z,f)+\eta^2Q_0(z,f)&=P_{03}(z,f)^3(f-1)^{\theta_{31}}(f-\eta)^{\theta_{32}}(f-\eta^2)^{\theta_{33}},
    \end{split}
    \end{equation}
where $\theta_{ij}\in\{0,1\}$ and $\theta_{1j}+\theta_{2j}+\theta_{3j}=1$, and $P_{01}(z,f)$, $P_{02}(z,f)$ and $P_{03}(z,f)$ are three polynomials in $f$ such that $P_{01}(z,f)P_{02}(z,f)P_{03}(z,f)=P_0(z,f)$ and any two of them have no common roots. Denote the degrees of the three polynomials on the RHS of equation \eqref{eq  restrict a01 e  p eq2    juy1} by $p_1$, $p_2$ and $p_3$, respectively. Consider the case where $p_1\leq p_2\leq p_3$. By eliminating $P_1(z,f)$ from the first two equations in \eqref{eq  restrict a01 e  p eq2    juy1} and then from the second and the third equations, respectively, we can obtain two expressions for the polynomial $Q_0(z,f)$. By comparing the degrees of these two polynomials, we see that the cases where $p_1< p_2< p_3$ and $p_1= p_2< p_3$ cannot occur. Therefore, we have $p_1<p_2=p_3$ or $p_1=p_2=p_3$. Moreover, if the case $p_1<p_2=p_3$ occurs, then we must have that the two polynomials $P_1(z,f)$ and $Q_0(z,f)$ have the same degrees and also that the leading coefficients of them have opposite signs; in this case, we must have $p<q$ and we see from equations in \eqref{eq  restrict a01 e  p eq2    juy1} that $p_1+p_2+p_3=p=q-3$. It follows that $p_1-p_2+3p_2=3q_0-3=3p_2-3$, which gives $p_2-p_1=3$. We conclude that the three integers $p_1$, $p_2$ and $p_3$ are equal to each other or one is less by~3 than the other two. Recall that $P_0(z,f)$ and $Q_0(z,f)$ can have simple roots only. From the above discussions, we have only the following two possibilities:
    \begin{equation}\label{eq  restrict a01 e  p eq2    juy2}
    \begin{split}
    P_1(z,f)+Q_0(z,f)&=P_{01}(z,f)^3(f-\eta_1),\\
    P_1(z,f)+\eta Q_0(z,f)&=P_{02}(z,f)^3(f-\eta_2),\\
    P_1(z,f)+\eta^2Q_0(z,f)&=P_{03}(z,f)^3(f-\eta_3),
    \end{split}
    \end{equation}
where $\eta_1$, $\eta_2$ and $\eta_3$ are the three distinct roots of~1, or
    \begin{equation}\label{eq  restrict a01 e  p eq2    juy20}
    \begin{split}
    P_1(z,f)+Q_0(z,f)&=P_{01}(z,f)^3(f^3-1),\\
    P_1(z,f)+\eta Q_0(z,f)&=P_{02}(z,f)^3,\\
    P_1(z,f)+\eta^2Q_0(z,f)&=P_{03}(z,f)^3.
    \end{split}
    \end{equation}
Denote the degrees of the three polynomials $P_{01}(z,f)$, $P_{02}(z,f)$ and $P_{03}(z,f)$ by $s_1$, $s_2$ and $s_3$, respectively, and denote the leading coefficients of the three polynomials $P_{01}(z,f)^3$, $P_{02}(z,f)^3$ and $P_{03}(z,f)^3$ by $a_{11}$, $a_{12}$ and $a_{13}$, respectively.

For the three equations in \eqref{eq  restrict a01 e  p eq2    juy2}, when $p=q+3$, we consider the case when $\eta_2=\eta\eta_1$ and $\eta_3=\eta^2\eta_1$ and also that $P_{01}(z,f)^3=P_{02}(z,\eta f)^3=P_{03}(z,\eta^2 f)^3$, i.e., $P_{02}(z,f)^3=P_{01}(z,\eta^2f)^3$ and $P_{03}(z,f)^3=P_{01}(z,\eta f)^3$. In this case, if $P_1(z,f)$ is of the form $P_1(z,f)=P_{11}(z,f^3)f$ for some polynomial $P_{11}(z,f^3)$ in $f^3$ and $Q_0(z,f)$ is of the form $Q_0(z,f)=Q_{11}(z,f^3)$ for some polynomial $Q_{11}(z,f^3)$ in $f^3$, then by doing the transformation $f\to \eta f$ for the second equation in \eqref{eq  restrict a01 e  p eq2    juy2} and dividing by $\eta$ both sides of the resulting equation and applying the transformation $f\to \eta^2f$ for the third equation in \eqref{eq  restrict a01 e  p eq2    juy2} and dividing by $\eta^2$ both sides of the resulting equation, respectively, we get exactly the first equation in \eqref{eq  restrict a01 e  p eq2    juy2}. With $P_1(z,f)$ and $Q_0(z,f)$ above, the first equation in \eqref{eq  restrict a01 e  p eq2    juy2} is solvable by the same arguments as in the discussions previously.
We note that when $p=q$ and $A=1$, we can deal with \eqref{eq  restrict a01 e  p eq2    juy2} in exactly the same way as above by just changing the positions of $P_1(z,f)$ and $Q_0(z,f)$. On the other hand, for the case when $p=q-3$, we consider equations in \eqref{eq  restrict a01 e  p eq2    juy20}. We have
    \begin{equation}\label{eq  restrict a01 e  p eq2    juy21aswhhgg0jhryt7}
    \begin{split}
    P_{01}(z,f)^3(f^3-1)+\eta P_{02}(z,f)^3+\eta^2P_{03}(z,f)^3=0.
    \end{split}
    \end{equation}
When $p=q-3$, we may suppose that the leading coefficient of $P_1(z,f)$ is $-\eta$. Then if we write
\begin{equation}\label{Theorem  2c  eq12 first1 solv154}
\begin{split}
P_{01}(z,f)^3&=a_{11}[(f^3-a_{s_1})(f^3-a_{s_1-1})\cdots(f^3-a_1)]^3,\\
P_{02}(z,f)^3&=a_{12}[(f^3-b_{s_1})(f^3-b_{s_1-1})\cdots(f^3-b_1)]^3,\\
P_{03}(z,f)^3&=a_{13}[(f^3-c_{s_1})(f^3-c_{s_1-1})\cdots(f^3-c_1)]^3f^3,\\
\end{split}
\end{equation}
where $a_i$, $b_j$ and $c_k$ are in general algebraic functions and $a_{11}=-\eta+1$, $a_{13}=-\eta+\eta^2$ and $a_{11}a_{12}a_{13}=A$. Denoting $g=f^3$, then equation \eqref{eq  restrict a01 e  p eq2    juy21aswhhgg0jhryt7} becomes a polynomial equation with respect to $g$. By comparing the coefficients on both sides of this equation  we obtain $3s_1+1$ polynomial equations with respect to the unknowns $A$, $a_i$, $b_j$ and $c_k$, whose combined number is $3s_1+1$. For example, when $s_1=0$, we get $P_{01}(z,f)^3=1-\eta$, $P_{02}(z,f)^3=\eta^2-1$ and $P_{03}(z,f)^3=(\eta^2-\eta)f^3$. It follows that $P_0(z,f)^3=3(\eta-\eta^2)f^3(f^3-1)$, $Q_0(z,f)=f^3+\eta$ and $P_1(z,f)=-\eta(f^3+\eta^2)$. By looking at the above examples, one may ask if $P_0(z,f)$ and $Q_0(z,f)$ are both polynomials in $f^3$ for any given $p_0$ and $q_0$.

Solutions of the eight equations \eqref{Theorem  2c  eq12 first0}, \eqref{Theorem  2c  eq12 third0}, \eqref{Theorem  2c  eq12 fourth0}, \eqref{Theorem  2c  eq17}, \eqref{Theorem  2c  eq0}, \eqref{Theorem  2c  eq12 seventh0}, \eqref{Theorem  2c  eq12 eightth0} and \eqref{Theorem  2c  eq12 ninth0} in the autonomous case are elliptic functions. We will discuss them further in section~\ref{discussions}. Below we begin to prove Theorem~\ref{Theorem  2c}.

\renewcommand{\proofname}{Proof of Theorem~\ref{Theorem  2c}.}
\begin{proof}
Since $n\mid |p-q|$, then we have $2\leq N_c\leq 4$. Below we consider the three cases where $N_c= 2$, $N_c=3$ and $N_c=4$, respectively.

\vskip 4pt

\noindent\textbf{Case~1:} $N_c= 2$.

\vskip 4pt

Since $n\mid |p-q|$, then we have the following three possibilities:
\begin{eqnarray}
\overline{f}^n &=& \frac{P_0(z,f)^n}{Q_0(z,f)^n}(f-\alpha_1)^{k_1}(f-\alpha_2)^{k_2}, \label{eq  restrict a08 e}\\
\overline{f}^n &=& \frac{P_0(z,f)^n}{Q_0(z,f)^n}\frac{1}{(f-\beta_1)^{l_1}(f-\beta_2)^{l_2}}, \label{eq  restrict a09 e}\\
\overline{f}^n &=& \frac{P_0(z,f)^n}{Q_0(z,f)^n}\frac{(f-\alpha_1)^{k_1}}{(f-\beta_1)^{l_1}},  \label{eq  restrict a1}
\end{eqnarray}
where in \eqref{eq  restrict a08 e} we have $n\mid (k_1+k_2)$, in \eqref{eq  restrict a09 e} we have $n\mid (l_1+l_2)$ and in \eqref{eq  restrict a1} we have $n\mid |k_1-l_1|$. Since $n\mid |p-q|$, we see that $n\nmid k_i$ and $n\nmid l_j$. For convenience, we denote $\alpha_i$ or $\beta_j$ in each of the above equations by $\gamma_1$ and $\gamma_2$. By Lemmas~\ref{basiclemma} and~\ref{keylemma} it follows that if $\gamma_i\not\equiv0$, then $\omega\gamma_i$ is a completely ramified rational function of $f$ with multiplicity at least~$n$, where $\omega$ is the $n$-th root of $1$. Therefore, if one of $\gamma_1$ and $\gamma_2$ is zero, we must have $n=2$; otherwise, say $\gamma_1=0$, if $n\geq 3$, then $0$ and $\omega\gamma_2$ are completely ramified rational functions of $f$ with multiplicity at least~$3$, a contradiction to the inequality \eqref{multiplicityinequality}.

For equation \eqref{eq  restrict a09 e}, we divide the following two cases: (1) $\beta_1$ and $\beta_2$ are both non-zero; or (2) at least one of $\beta_1$ and $\beta_2$ is zero. In the first case, if $p_0\geq 2$ and $P_0(z,f)^n=a_pf^{np_0}$ particularly, then the analysis in the proof of Lemma~\ref{basiclemma0} applies and~$0$ is a Picard exceptional rational function of $f$; also, we have $p\geq q$ for otherwise $\infty$ is also a Picard exceptional rational function of $f$ by analyzing on the poles of $f$ and it follows that $\beta_1$ and $\beta_2$ are both Picard exceptional rational functions of $f$, a contradiction to Picard's theorem. Moreover, by the inequality \eqref{multiplicityinequality}, we must have $n=2$ and $\beta_1+\beta_2=0$ and it follows by Lemma~\ref{keylemma0} that $l_1=l_2=1$. Then by doing a bilinear transformation $f\to 1/f$, we get equation \eqref{proof3 eqre   1}, which leads to equation \eqref{yanagiharaeq0    5} in Theorem~\ref{Theorem  1c}. Otherwise, by doing a bilinear transformation $f\to 1/f$, we get equation \eqref{eq  restrict a08 e}. In the second case, since $n=2$, then by the bilinear transformation $f\to 1/f$, we get equation \eqref{eq  restrict c1}, which leads to equation \eqref{Theorem  2b  eq0} in Theorem~\ref{Theorem  2b}. Therefore, in this section we only need to consider equations \eqref{eq  restrict a08 e} and \eqref{eq  restrict a1}.

Further, for equation \eqref{eq  restrict a08 e}, we may suppose that $\alpha_1$ and $\alpha_2$ are both non-zero; otherwise, say $\alpha_1=0$, by doing the transformation $f\to1/f$ we get equation \eqref{cannotouu} with $P_1(z,f)$ being a polynomial in $f$ of degree $d$ and $Q_1(z,f)$ being a polynomial in $f$ of degree $d-k_1$, which cannot have any meromorphic solution as shown in the proof of Theorem~\ref{Theorem  2b}. For equation \eqref{eq  restrict a1}, when $\alpha_1=0$, by Lemma~\ref{basiclemma0} it follows that $0$ is a Picard exceptional rational function of $f$. Similarly as in the previous paragraph, it follows that $p\geq q$. Moreover, we have $P_0(z,f)^n=a_p$; otherwise, the roots of $P_0(z,f)$ are also Picard's exceptional rational functions, a contradiction to Picard's theorem. Since $n=2$, then by doing a bilinear transformation $f\to 1/f$, we get equation \eqref{proof4 eqre   1}, which leads to equation \eqref{yanagiharaeq0    9} in Theorem~\ref{Theorem  1c}. When $\beta_1=0$, since $n=2$, by doing a bilinear transformation $f\to 1/f$, we get equation \eqref{eq  restrict c3} since $2\mid |k_1-l_1|$, which cannot have any meromorphic solution as shown in the proof of Theorem~\ref{Theorem  2b}. Therefore, in this section we only need to consider equations \eqref{eq  restrict a08 e} and \eqref{eq  restrict a1} for the case where $\alpha_i$ and $\beta_j$ are both non-zero. Under this assumptions, below we consider these two equations separately.

\vskip 4pt

\noindent\textbf{Subcase~1:} Equation \eqref{eq  restrict a08 e} with $\alpha_1\alpha_2\not\equiv0$.

\vskip 4pt

In this case, we claim that $n=2$. Suppose that $n\geq 3$. By Lemma~\ref{keylemma} it follows that $\omega\alpha_1$ and $\omega\alpha_2$ are completely ramified rational functions of $f$ with multiplicity at least~3, where $\omega$ is the $n$-th root of~1. By the inequality \eqref{multiplicityinequality} we must have $n=3$ and also that $\alpha_1^3=\alpha_2^3$. However, since $3\mid(k_1+k_2)$, we have $k_1\nmid 3$ or $k_2\nmid 3$ and in either case we will get a contradiction to Lemma~\ref{keylemma0}. Therefore, $n=2$.

When $\alpha_1+\alpha_2\not=0$, $\pm \alpha_1$ and $\pm\alpha_2$ are all completely ramified rational functions of $f$ with multiplicities~2. By Lemma~\ref{keylemma0} it follows that $P_0(z,f)$ and $Q_0(z,f)$ both have simple roots only and also that $k_1=k_2=1$ and $p-q\in \{-2,0,2\}$. We consider
\begin{equation}\label{eq  restrict a08 e  p eq2s2 t purt1avbf0jjjjy}
\overline{f}^2-\overline{\alpha}_1^2= \frac{P_0(z,f)^2(f-\alpha_1)(f-\alpha_2)-\overline{\alpha}_1^2Q_0(z,f)^2}{Q_0(z,f)^2}.
\end{equation}
By analyzing the multiplicities of poles of $f$ as in the proof Lemma~\ref{keylemma1} we get that the numerator of the RHS of \eqref{eq  restrict a08 e  p eq2s2 t purt1avbf0jjjjy} is a polynomial in $f$ with even degree; moreover, if some root of this polynomial is not equal to $\pm\alpha_1$ or $\pm \alpha_2$, then this root has order two. Therefore, $-\alpha_1$ and $-\alpha_2$ are either both simple roots of the numerator of the RHS of \eqref{eq  restrict a08 e  p eq2s2 t purt1avbf0jjjjy}, or neither of them are. In the first case, we then consider $\overline{f}^2-\overline{\alpha}_2^2$ and conclude by Lemma~\ref{keylemma1} that the polynomial $P_0(z,f)^2(f-\alpha_1)(f-\alpha_2)-\overline{\alpha}_2^2Q_0(z,f)^2$ is a square of some polynomial in $f$. In the latter case, we claim that $-\alpha_1$ and $-\alpha_2$ are both simple roots of the polynomial $P_0(z,f)^2(f-\alpha_1)(f-\alpha_2)-\overline{\alpha}_2^2Q_0(z,f)^2$ when considering $\overline{f}^2-\overline{\alpha}_2^2$. Otherwise, we may suppose that $R(z,-\alpha_1)=\gamma^2$ for some algebraic function $\gamma^2$ which is distinct from $\overline{\alpha}_1^2$ and $\overline{\alpha}_2^2$. Let $z_0\in \mathbb{C}$ be such that $f(z_0)+\alpha_1(z_0)=0$. Then we have $f(z_0+1)^2-\gamma(z_0)^2=0$ and by applying the analysis in the proof of Lemma~\ref{keylemma0} we get that $z_0$ is a root of the equation $f(z_0+1)-\gamma(z_0)=0$ or $f(z_0+1)+\gamma(z_0)=0$ with multiplicity~2 and there are $T(r,f)+o(T(r,f))$ many such points, where $r\to\infty$ outside an exceptional set of finite linear measure. But we then get a contradiction to the inequality \eqref{multiplicityinequality} by computing the quantity $\Theta(\gamma,\overline{f})$ or $\Theta(-\gamma,\overline{f})$ as in the proof of Lemma~\ref{keylemma0} since at least one of these two quantities is strictly positive. Without loss of generality, we may suppose the following two equations:
\begin{equation}\label{eq  restrict a08 e  p eq2s2 t purt1avbf0jjj0}
\overline{f}^2-\overline{\alpha}_1^2= \frac{P_1(z,f)^2(f+\alpha_1)(f+\alpha_2)}{Q_0(z,f)^2},
\end{equation}
and
\begin{equation}\label{eq  restrict a08 e  p eq2s2 t purt1avbf0jjj1}
\overline{f}^2-\overline{\alpha}_2^2= \frac{P_2(z,f)^2}{Q_0(z,f)^2},
\end{equation}
where $P_1(z,f)$ and $P_2(z,f)$ are two polynomials in $f$. By Lemma~\ref{keylemma1}, $P_1(z,f)$ and $P_2(z,f)$ both have simple roots only and none of the roots of $P_1(z,f)$ and $P_2(z,f)$ is equal to $\pm\alpha_1$ or $\pm \alpha_2$. Moreover, when $p=q$, if $a_p=\overline{\alpha}_1^2$ or $a_p=\overline{\alpha}_2^2$, then the degree of the numerator in \eqref{eq  restrict a08 e  p eq2s2 t purt1avbf0jjj0} or in \eqref{eq  restrict a08 e  p eq2s2 t purt1avbf0jjj1} decreases since the terms with the highest degrees in $P_0(z,f)^2$ and $Q_0(z,f)^2$ cancel out when considering $\overline{f}^2-\overline{\alpha}_1^2$ or $\overline{f}^2-\overline{\alpha}_2^2$; by Lemma~\ref{keylemma1} we see that it decreases by $2$.
By doing the transformation $f\to \alpha_1 f$, we get the first equation of Theorem~\ref{Theorem  2c}.

When $\alpha_1+\alpha_2=0$, we may let $\alpha_1=1$ and $\alpha_2=-1$ by doing a linear transformation $f\to \alpha_1f$. We consider
\begin{equation}\label{eq  restrict a08 e  p eq2s2 t purt1avbf0}
\overline{f}^2-1= \frac{P_0(z,f)^2(f-1)^{k_1}(f+1)^{k_2}-Q_0(z,f)^2}{Q_0(z,f)^2}.
\end{equation}
If the numerator of the RHS of \eqref{eq  restrict a08 e  p eq2s2 t purt1avbf0} has at least three distinct roots, say $\gamma_i$, $i=1,2,3$, of odd order, then by applying the analysis in the proof of Lemma~\ref{basiclemma} together with the fact that $\pm1$ are completely ramified rational functions of $f$ and that the roots of $f\pm1=0$ have even multiplicities, we obtain that $\gamma_1,\gamma_2,\gamma_3$ are all completely ramified rational functions of $f$, a contradiction to Theorem~\ref{completelyrm} since $\gamma_i$ are all distinct from $\pm1$. Therefore, the numerator of the RHS of \eqref{eq  restrict a08 e  p eq2s2 t purt1avbf0} can have at most two distinct roots of odd order. Suppose that there is only one such root, say $\gamma_1$. Since $p$ and $q$ are both even integers, then we must have $p=q$ and $A=1$ in which case the terms with the highest degrees in the polynomials $P_0(z,f)^2(f-1)^{k_1}(f+1)^{k_2}$ and $Q_0(z,f)^2$ cancel out when considering \eqref{eq  restrict a08 e  p eq2s2 t purt1avbf0} so that the degree of the polynomial $P_0(z,f)^2(f-1)^{k_1}(f+1)^{k_2}-Q_0(z,f)^2$ decreases to be an odd integer. Then by considering the multiplicity of the poles of $f$ as in the proof of Lemma~\ref{basiclemma}, we see from \eqref{eq  restrict a08 e  p eq2s2 t purt1avbf0} that $\infty$ is also a completely ramified rational function of $f$. If $\gamma_1=0$, then $0$ is a completely ramified rational function of $f$. Then by considering the multiplicities of the roots of $f\pm1=0$ as in the proof of Lemma~\ref{basiclemma} and then by Lemma~\ref{keylemma}, we obtain from \eqref{eq  restrict a08 e} that $\pm 1$ both have multiplicities at least~4, a contradiction to the inequality \eqref{multiplicityinequality}.
On the other hand, if $\gamma_1\not\equiv0$, then by Lemma~\ref{keylemma} it follows that $\pm \gamma_1$ are both completely ramified rational functions of $f$, a contradiction to Theorem~\ref{completelyrm} since $\infty$ is also a completely ramified rational function of $f$. Therefore,
the numerator of the RHS of \eqref{eq  restrict a08 e  p eq2s2 t purt1avbf0} has no roots of odd order, or  has two distinct roots of odd order.

If the numerator of the RHS of \eqref{eq  restrict a08 e  p eq2s2 t purt1avbf0} has no roots of odd order, then we have $P_0(z,f)^2(f-1)^{k_1}(f+1)^{k_2}-1=P_1(z,f)^2$, where $P_1(z,f)$ is a polynomial in $f$. This gives the second equation of Theorem~\ref{Theorem  2c}.

If the numerator of the RHS of \eqref{eq  restrict a08 e  p eq2s2 t purt1avbf0} has two distinct roots of odd order, then we have $P_0(z,f)^2(f-1)^{k_1}(f+1)^{k_2}-Q_0(z,f)^2=P_1(z,f)^2(f-\gamma_1)^{t_1}(f-\gamma_2)^{t_2}$ for some polynomial $P_1(z,f)$ in $f$, and $\gamma_1$ and $\gamma_2$ are distinct from each other. In this case, $\gamma_1$ and $\gamma_2$ are also both completely ramified rational functions of $f$ and by Lemma~\ref{keylemma} it follows that if $\gamma_i\not\equiv0$, then $\pm\gamma_i$ are both completely ramified rational functions of $f$. By Theorem~\ref{completelyrm} we must have $\gamma_1+\gamma_2=0$. Moreover, by Lemmas~\ref{keylemma0} and~\ref{keylemma1}, we have $k_1=k_2=t_1=t_2=1$. Therefore, by denoting $\gamma=\gamma_1$, we have
\begin{equation}\label{eq  restrict a08 e  p eq2s2 t}
\overline{f}^2-1= \frac{P_1(z,f)^2}{Q_0(z,f)^2}(f^2-\gamma^2).
\end{equation}
Now $\pm 1$ and $\pm \gamma$ are all completely ramified rational functions of $f$ with multiplicities~2 and $f$ has no other completely ramified rational functions and has no Picard exceptional rational functions. Thus by applying the analysis on \eqref{eq  restrict a08 e  p eq2s2 t purt1avbf0} to $\overline{f}^2-\overline{\gamma}^2$ we conclude that the polynomial $P_0(z,f)^2-\overline{\gamma}^2Q_0(z,f)^2$ cannot have any root of odd order; otherwise, this root is distinct from $\pm1$ and $\pm \gamma$ and is also a completely ramified rational function of $f$, a contradiction to Theorem~\ref{completelyrm}. Therefore, we have the following
\begin{equation}\label{eq  restrict a08 e  p eq2s2 t purt0a  pp}
\overline{f}^2-\overline{\gamma}^2= \frac{P_2(z,f)^2}{Q_0(z,f)^2},
\end{equation}
where $P_2(z,f)$ is a polynomial in $f$. By Lemmas~\ref{keylemma0} and~\ref{keylemma1}, we have $P_0(z,f)$, $Q_0(z,f)$, $P_1(z,f)$ and $P_2(z,f)$ all have simple roots only and none of the roots of $P_1(z,f)$ and $P_2(z,f)$ is equal to $\pm1$ or $\pm \gamma$, and the degrees $p$ and $q$ satisfy $p-q\in\{-2,0,2\}$.
As for equations \eqref{eq  restrict a08 e  p eq2s2 t purt1avbf0jjj0} and \eqref{eq  restrict a08 e  p eq2s2 t purt1avbf0jjj1}, when $p=q$, if $A=1$ or $A=\overline{\gamma}^2$, the degree of the numerator in \eqref{eq  restrict a08 e  p eq2s2 t} or in \eqref{eq  restrict a08 e  p eq2s2 t purt0a  pp} decreases and by Lemma~\ref{keylemma1} it decreases by $2$. This gives the third equation of Theorem~\ref{Theorem  2c}.

\vskip 4pt

\noindent\textbf{Subcase~2:} Equation \eqref{eq  restrict a1} with $\alpha_1\beta_1\not\equiv0$.

\vskip 4pt

In this case, we discuss the two cases $n=2$ and $n\geq 3$ separately.

When $n=2$, by Lemmas~\ref{basiclemma} and~\ref{keylemma} it follows that $\pm\alpha_1$, as well as $\pm\beta_1$, are completely ramified rational functions of $f$. We may let $\alpha_1=\kappa$ and $\beta_1=1$ by doing a linear transformation $f\to \beta_1f$. We consider
\begin{equation}\label{eq  restrict a1 fgweq0}
\overline{f}^2-1= \frac{P_0(z,f)^2(f-\kappa)^{k_1}-Q_0(z,f)^2(f-1)^{l_1}}{Q_0(z,f)^2(f-1)^{l_1}}.
\end{equation}
Recall that the leading coefficient of the polynomial $P_0(z,f)^2$ is denoted by $A$ and that the polynomial $Q_0(z,f)^2$ is monic. When $p=q$ and $A=1$, the degree of the numerator of the RHS of \eqref{eq  restrict a1 fgweq0} decreases due to the cancellation of the terms with the highest degrees in $P_0(z,f)^2(f-\kappa)^{k_1}$ and $Q_0(z,f)^2(f-1)^{l_1}$. Suppose that the degree of the polynomial $P_0(z,f)^2(f-\kappa)^{k_1}-Q_0(z,f)^2(f-1)^{l_1}$ in $f$ decreases to be an even integer. By considering the multiplicities of the poles of $f$ together with the fact that $\pm1$ are both completely ramified rational functions of $f$ and that the roots of $f\pm 1=0$ have even multiplicities with at most finitely many exceptions, we get that $\infty$ is also a completely ramified rational function of $f$. Further, by considering the multiplicities of the roots of $f-1=0$ as in the proof of Lemma~\ref{keylemma}, we obtain from \eqref{eq  restrict a1} that $\pm 1$ both have multiplicities at least~4 and it follows that $\infty$ is also a completely ramified function of $f$ with multiplicity at least~4, a contradiction to the inequality \eqref{multiplicityinequality}. This implies that the numerator of the RHS of equation \eqref{eq  restrict a1 fgweq0} always has a root, say $\gamma$, of odd order. Then by applying the same analysis as in the proof of Lemma~\ref{basiclemma} together with the fact that $\pm 1$ are completely ramified rational functions of $f$ we obtain that $\gamma$ is a completely ramified rational function of $f$. If $\kappa\not=-1$, then by Lemma~\ref{keylemma}, $\pm 1$ and $\pm \kappa$ are all completely ramified rational functions of $f$ and thus by Theorem~\ref{completelyrm} we have $\gamma=-1$ or $\gamma=-\kappa$ and thus $\gamma\not=0$. If $\kappa=-1$ and $\gamma=0$, then by considering the multiplicities of the roots of $f+1=0$ as in the proof of Lemma~\ref{basiclemma} and then by Lemma~\ref{keylemma}, we obtain from \eqref{eq  restrict a1} that $\pm 1$ both have multiplicities at least~4. But it follows by repeating the analysis after \eqref{eq  restrict a1 fgweq0} that $0$ is a completely ramified rational function of $f$ with multiplicity at least~4, a contradiction to the inequality \eqref{multiplicityinequality}. Therefore, when $\kappa=-1$, we also have $\gamma\not=0$. Now, by Lemma~\ref{keylemma} it follows that $\pm\gamma$ are both completely ramified rational functions of $f$. From the above reasoning, we see that $f$ has four completely ramified rational functions $\pm 1$ and $\pm \kappa$ (or $\pm\gamma$), all of which have multiplicities~2. By Lemma~\ref{keylemma0} we must have $k_1=l_1=1$ and all the roots of $P_0(z,f)$ and $Q_0(z,f)$ are simple and also that $p_0-q_0\in\{-1,0,1\}$. As for equations \eqref{eq  restrict a08 e  p eq2s2 t purt1avbf0jjj0} and \eqref{eq  restrict a08 e  p eq2s2 t purt1avbf0jjj1}, when $p=q$, if the degree of the numerator in \eqref{eq  restrict a1 fgweq0} decreases, then by Lemma~\ref{keylemma1} it decreases by $2$. Now we have
\begin{equation}\label{eq  restrict a1 gf}
\overline{f}^2 = \frac{P_0(z,f)^2}{Q_0(z,f)^2}\frac{(f-\kappa)}{(f-1)}.
\end{equation}
When $\kappa=-1$, from the above discussions we know that $P_0(z,f)^2(f+1)-Q_0(z,f)^2(f-1)$ is of the form $P_1(z,f)^2(f-\gamma)$ for some non-zero algebraic function $\gamma$ and a polynomial $P_1(z,f)$ in $f$. Since $\pm1$ and $\pm \gamma$ are all completely ramified rational functions of $f$, then by Theorem~\ref{completelyrm} and Lemma~\ref{keylemma1} and considering $\overline{f}^2-\overline{\gamma}^2$, we see that $P_0(z,f)^2(f+1)-\overline{\gamma}^2Q_0(z,f)^2(f-1)$ must be of the form $P_2(z,f)^2(f+\gamma)$ for some polynomial $P_2(z,f)$ in $f$; moreover, both $P_1(z,f)$ and $P_2(z,f)$ have simple roots only and none of these roots equals $\pm1$ or $\pm \gamma$. Therefore, we have
\begin{equation}\label{eq  restrict a1 gfhj092}
\overline{f}^2-1= \frac{P_1(z,f)^2}{Q_0(z,f)^2}\frac{(f-\gamma)}{(f-1)},
\end{equation}
and
\begin{equation}\label{eq  restrict a1 gfhj093}
\overline{f}^2-\overline{\gamma}^2= \frac{P_2(z,f)^2}{Q_0(z,f)^2}\frac{(f+\gamma)}{(f-1)}.
\end{equation}
Note that, when $p=q$, if the degree of the numerator in \eqref{eq  restrict a1 gfhj092} or in \eqref{eq  restrict a1 gfhj093} decreases, then by Lemma~\ref{keylemma1} it decreases by $2$. On the other hand, when $\kappa\not=-1$, by Theorem~\ref{completelyrm} and Lemma~\ref{keylemma1} and considering $\overline{f}^2-1$ and $\overline{f}^2-\overline{\kappa}^2$, respectively, we see that $P_0(z,f)^2(f-\kappa)-Q_0(z,f)^2(f-1)$ (and also $P_0(z,f)^2(f-\kappa)-\overline{\kappa}^2Q_0(z,f)^2(f-1)$) must be of the form $P_1(z,f)^2(f+1)$ or $P_2(z,f)^2(f+\kappa)$ for some polynomials $P_1(z,f)$ and $P_2(z,f)$ in $f$ with simple roots only and none of these roots equals $\pm \kappa$ and $\pm 1$. Without loss of generality, we may consider the following two equations:
\begin{equation}\label{eq  restrict a1 gfhj090}
\overline{f}^2-1= \frac{P_1(z,f)^2}{Q_0(z,f)^2}\frac{(f+\kappa)}{(f-1)},
\end{equation}
and
\begin{equation}\label{eq  restrict a1 gfhj091}
\overline{f}^2-\overline{\kappa}^2= \frac{P_2(z,f)^2}{Q_0(z,f)^2}\frac{(f+1)}{(f-1)}.
\end{equation}
Note that, when $p=q$, if the degree of the numerator in \eqref{eq  restrict a1 gfhj090} or in \eqref{eq  restrict a1 gfhj091} decreases, then by Lemma~\ref{keylemma1} it decreases by $2$. This gives the fourth equation of Theorem~\ref{Theorem  2c}.

Now we consider the case when $n\geq 3$. Since $\omega\alpha_1$ and $\omega\beta_1$ are all completely ramified rational functions of $f$ with multiplicities at least~3, where $\omega$ is the $n$-th root of~1, then by the inequality \eqref{multiplicityinequality} we must have $n=3$ and $\alpha_1^3=\beta_1^3$. By Lemma~\ref{keylemma0} we conclude that $k_1=l_1=1$. We fix one $\eta$ such that $\eta^2+\eta+1=0$ and choose without loss of generality that $\beta_1=\eta\alpha_1$. We may let $\alpha_1=1$ by doing a linear transformation $f\to\alpha_1f$. Then we have
\begin{equation}\label{eq  restrict a1     tttt0}
\overline{f}^3= \frac{P_0(z,f)^3(f-1)}{Q_0(z,f)^3(f-\eta)}.
\end{equation}
Also, by Lemma~\ref{keylemma0}, we conclude that $P_0(z,f)$ and $Q_0(z,f)$ can have simple roots only and $p_0-q_0\in\{-1,0,1\}$. We consider
\begin{equation}\label{eq  restrict a1 fg}
\overline{f}^3-1= \frac{P_0(z,f)^3(f-1)-Q_0(z,f)^3(f-\eta)}{Q_0(z,f)^3(f-\eta)}.
\end{equation}
Note that $f$ has three completely ramified rational functions with multiplicities~3. As for equations \eqref{eq  restrict a08 e  p eq2s2 t purt1avbf0jjj0} and \eqref{eq  restrict a08 e  p eq2s2 t purt1avbf0jjj1}, when $p_0=q_0$ and the leading coefficient of the polynomial $P_0(z,f)^3$ satisfies $A=1$, the degree of the numerator of the RHS of \eqref{eq  restrict a1 fg} decreases and by Lemma~\ref{keylemma1} it decreases by $3$. Therefore, the numerator of the RHS of the equation above always has one root of order $l_1$ such that $3\nmid l_1$ and by Lemma~\ref{keylemma1} we must have $l_1=1$ and then by the inequality \eqref{multiplicityinequality} we see that the root must be $\eta^2$. Therefore, the numerator of the RHS of \eqref{eq  restrict a1 fg} is of the form $P_1(z,f)^3(f-\eta^2)$ for a polynomial $P_1(z,f)$ in $f$ with simple roots only, i.e.,
\begin{equation}\label{eq  restrict a1 fgtrytr}
\overline{f}^3-1= \frac{P_1(z,f)^3(f-\eta^2)}{Q_0(z,f)^3(f-\eta)}.
\end{equation}
This gives the fifth equation of Theorem~\ref{Theorem  2c}.

\vskip 4pt

\noindent\textbf{Case~2:} $N_c=3$.

\vskip 4pt

Since $n\mid |p-q|$, we must have $n\geq 3$. For convenience, we denote the three roots by $\gamma_1$, $\gamma_2$ and $\gamma_3$ and their orders by $t_1$, $t_2$ and $t_3$, respectively. Without loss of generality, we may suppose that $\gamma_1\gamma_2\not=0$. Since $n\geq 3$, then by Lemmas~\ref{basiclemma} and~\ref{keylemma} it follows that $\omega\gamma_1$ is a completely ramified rational function of $f$, where $\omega$ is the $n$-th root of~1, and so by the inequality \eqref{multiplicityinequality} we must have $n=3$ or $4$. However, when $n=4$, $t_1$ and $t_2$ must be both even integers; otherwise, $\omega\gamma_1$ (or $\omega\gamma_2$) would have multiplicity at least~$4$, where $\omega$ is the fourth root of~1, which is impossible. But since $n\mid |p-q|$, we see that $t_3$ is also an even integer, a contradiction to our assumption that at least one of $\alpha_i$ and $\beta_j$ in \eqref{P} and \eqref{Q} has no common factors with $n$. Therefore, we must have $n=3$. We see that $\eta\gamma_1$ has multiplicity~$3$ since we must have $(n,t_1)=1$, where $\eta$ is the cubic root of~1. Moreover, by the inequality \eqref{multiplicityinequality} we have none of $\gamma_1$, $\gamma_2$ and $\gamma_3$ is zero and by Lemma~\ref{keylemma0} we also have $t_1=t_2=t_3=1$. By noting that $n\mid |p-q|$, when $n=3$ we have only the following two possibilities:
\begin{eqnarray}
\overline{f}^3 &=& \frac{P_0(z,f)^3}{Q_0(z,f)^3}(f-\alpha_1)(f-\alpha_2)(f-\alpha_3),
\label{eq  restrict a01 e}\\
\overline{f}^3 &=& \frac{P_0(z,f)^3}{Q_0(z,f)^3}\frac{1}{(f-\beta_1)(f-\beta_2)(f-\beta_3)}.
\label{eq  restrict a01 e1}
\end{eqnarray}
For each of the above two equations, by Lemma~\ref{keylemma0} we have that all the roots of $P_0(z,f)$ and $Q_0(z,f)$ are simple and also that $p-q\in\{-3,0,3\}$. Now, for equation \eqref{eq  restrict a01 e1}, if $0$ is the only root of $P_0(z,f)$, then $p=3$ and it follows that $q=6$ under our assumption that $n<d$; in this case we have $p_0=q_0=1$ and that $Q_0(z,f)$ has a non-zero root. Since none of $\alpha_i$ and $\beta_j$ is zero, then by doing a bilinear transformation $f\to 1/f$, both of the above two cases of equation \eqref{eq  restrict a01 e1} become \eqref{eq  restrict a01 e}. Thus we only need to consider equation \eqref{eq  restrict a01 e}.

Since $\eta\alpha_i$ has multiplicity~$3$, where $\eta$ is the cubic root of~1, then by the inequality \eqref{multiplicityinequality} we must have $\alpha_1^3=\alpha_2^3=\alpha_3^3$. We may let $\alpha_1=1$ by doing a linear transformation $f\to \alpha_1f$. We consider
\begin{equation}\label{eq  restrict a01 e  p eq1}
\overline{f}^3-1=\frac{P_0(z,f)^3(f^3-1)-Q_0(z,f)^3}{Q_0(z,f)^3}.
\end{equation}
Let $\eta$ a fixed cubic root of~1 such that $\eta^2+\eta+1=0$. Since $1$, $\eta$ and $\eta^2$ all have multiplicities~$3$, then by Lemma~\ref{keylemma1} we conclude that the numerator of the RHS of equation \eqref{eq  restrict a01 e  p eq1} is of the form $P_1(z,f)^3$ for some polynomial $P_1(z,f)$ in $f$ with simple roots only and these roots are distinct from $1$, $\eta$ and $\eta^2$. As for equations \eqref{eq  restrict a08 e  p eq2s2 t purt1avbf0jjj0} and \eqref{eq  restrict a08 e  p eq2s2 t purt1avbf0jjj1} when $p=q$ and $A=1$, the degree of the numerator in \eqref{eq  restrict a01 e  p eq1} decreases and by Lemma~\ref{keylemma1} it decreases by $3$.
This gives the sixth equation of Theorem~\ref{Theorem  2c}.

\vskip 4pt

\noindent\textbf{Case~3:} $N_c= 4$.

\vskip 4pt

In this case, by Lemma~\ref{basiclemma} we know that $\alpha_i$ and $\beta_j$ are all completely ramified rational functions of $f$. Then by the inequality \eqref{multiplicityinequality} we must have $n=2$. By noting that $2\mid |p-q|$ and Lemma~\ref{keylemma0} we have the following possibilities:
\begin{eqnarray}
\overline{f}^2 &=& \frac{P_0(z,f)^2}{Q_0(z,f)^2}\frac{(f-\alpha_1)(f-\alpha_2)}{(f-\beta_1)(f-\beta_2)},
\label{eq  restrict a06 e}\\
\overline{f}^2 &=& \frac{P_0(z,f)^2}{Q_0(z,f)^2}(f-\alpha_1)(f-\alpha_2)(f-\alpha_3)(f-\alpha_4),
 \label{eq  restrict a06 e0}\\
\overline{f}^2 &=& \frac{P_0(z,f)^2}{Q_0(z,f)^2}\frac{1}{(f-\beta_1)(f-\beta_2)(f-\beta_3)(f-\beta_4)}, \label{eq  restrict a003 nwe0}\\
\overline{f}^2 &=& \frac{P_0(z,f)^2}{Q_0(z,f)^2}\frac{(f-\alpha_1)(f-\alpha_2)(f-\alpha_3)}{(f-\beta_1)},
 \label{eq  restrict a07 e}\\
\overline{f}^2 &=& \frac{P_0(z,f)^2}{Q_0(z,f)^2}\frac{(f-\alpha_1)}{(f-\beta_1)(f-\beta_2)(f-\beta_3)}. \label{eq  restrict a005 nwe1}
\end{eqnarray}
For convenience, we denote the four roots $\alpha_i$ and $\beta_j$ in each of the above equations by $\gamma_1$, $\gamma_2$, $\gamma_3$ and $\gamma_4$. If $\gamma_i\not\equiv0$ for some $i$, then by Lemmas~\ref{basiclemma} and~\ref{keylemma} it follows that $\pm \gamma_i$ are both completely ramified rational functions of $f$ with multiplicities~2. This implies that none of $\gamma_1$, $\gamma_2$, $\gamma_3$ and $\gamma_4$ is zero for otherwise $f$ would have at least five completely ramified rational functions, a contradiction to Theorem~\ref{completelyrm}. Moreover, by the inequality \eqref{multiplicityinequality} we must have $\gamma_1^2=\gamma_2^2$ and $\gamma_3^2=\gamma_4^2$, apart from permutations. Also, by Lemma~\ref{keylemma0} we know that in each of the above equations all the roots of $P_0(z,f)$ and $Q_0(z,f)$ are simple and the degrees of $P(z,f)$ and $Q(z,f)$ satisfy $p-q\in\{-2,0,2\}$. In particular, for equation \eqref{eq  restrict a003 nwe0}, we see that if $0$ is the only root of $P_0(z,f)$ then we must have $p=2$ and $q=4$ under the the assumption that $n<d$. Therefore, by doing a linear transformation $f\to1/f$, equations \eqref{eq  restrict a003 nwe0} and \eqref{eq  restrict a005 nwe1} become \eqref{eq  restrict a06 e0} and \eqref{eq  restrict a07 e}, respectively. From the above discussions, we conclude that we only need to consider the three equations \eqref{eq  restrict a06 e}, \eqref{eq  restrict a06 e0} and \eqref{eq  restrict a07 e}.

Further, equation \eqref{eq  restrict a07 e} cannot have any meromorphic solution, as is shown below. From the previous discussions, we may suppose $\alpha_1+\alpha_2=0$ and $\alpha_3+\beta_1=0$. We consider
\begin{equation}\label{eq  restrict a06 e  p eq1  jhgde}
\overline{f}^2-\overline{\alpha}_1^2= \frac{P_0(z,f)^2(f^2-\alpha_1^2)(f+\beta_1)-\overline{\alpha}_1^2Q_0(z,f)^2(f-\beta_1)}{Q_0(z,f)^2(f-\beta_1)}.
\end{equation}
Since $\alpha_1$, $\alpha_2$, $\alpha_3$ and $\beta_1$ are four completely ramified rational functions of $f$, then by Lemma~\ref{keylemma1} we conclude that the numerator of the RHS of equation \eqref{eq  restrict a06 e  p eq1  jhgde} is of the form $P_1(z,f)^2$ for some polynomial $P_1(z,f)$ in $f$ with simple roots only and none of these roots is equal to~$\alpha_1$, $\alpha_2$, $\alpha_3$ or $\beta_1$. Note that $p-q\in\{-2,0,2\}$. Since the degrees $p$ and $q$ are both odd integers, this is possible only when $p=q$ and the leading coefficient $a_p$ of the numerator $P(z,f)$ satisfies $a_p=\overline{\alpha}_1^2$ so that the terms with the highest degree in the two polynomials $P_0(z,f)^2(f^2-\alpha_1^2)(f+\beta_1)$ and $\overline{\alpha}_1^2Q_0(z,f)^2(f-\beta_1)$ cancel out. It follows by these arguments that $a_p=\overline{\alpha}_1^2=\overline{\alpha}_2^2=\overline{\alpha}_3^2=\overline{\beta}_2^2$, which is impossible. Therefore, we only need to consider the two equations \eqref{eq  restrict a06 e} and \eqref{eq  restrict a06 e0}. Below we discuss them, respectively.

\vskip 4pt

\noindent\textbf{Subcase~1:} Equation \eqref{eq  restrict a06 e}.

\vskip 4pt

From the previous discussions, we have two cases to consider: (1) $\alpha_1+\alpha_2=0$ and $\beta_1+\beta_2=0$; or (2) $\alpha_1+\beta_1=0$ and $\alpha_2+\beta_2=0$.

When $\alpha_1+\alpha_2=0$ and $\beta_1+\beta_2=0$, we may let $\alpha_1=\kappa$ and $\beta_1=1$ by doing a linear transformation $f\to \beta_1 f$. We consider
\begin{equation}\label{eq  restrict a06 e  p eq1  e}
\overline{f}^2-1= \frac{P_0(z,f)^2(f^2-\kappa^2)-\overline{\kappa}^2Q_0(z,f)^2(f^2-1)}{Q_0(z,f)^2(f^2-1)}.
\end{equation}
Since~$\pm1$ and $\pm \kappa$ are four completely ramified rational functions of $f$, then by Lemma~\ref{keylemma1} we conclude that the numerator of the RHS of equation \eqref{eq  restrict a06 e  p eq1  e} is of the form $P_1(z,f)^2$ for some polynomial $P_1(z,f)$ in $f$ with simple roots only and none of these roots is $\pm 1$ or $\pm \kappa$. Similarly, by considering $\overline{f}^2-1$, we also have $P_0(z,f)^2(f^2-\kappa^2)-\overline{\kappa}^2Q_0(z,f)^2(f^2-1)=P_2(z,f)^2$ for some polynomial $P_2(z,f)$ in $f$ with simple roots only and none of these roots is $\pm 1$ or $\pm \kappa$. Now we have
\begin{equation}\label{eq  restrict a06 e  p eq1  e0tr0nv0}
\overline{f}^2-1= \frac{P_1(z,f)^2}{Q_0(z,f)^2(f^2-1)},
\end{equation}
and
\begin{equation}\label{eq  restrict a06 e  p eq1  e0tr0nv1}
\overline{f}^2-\overline{\kappa}^2= \frac{P_2(z,f)^2}{Q_0(z,f)^2(f^2-1)}.
\end{equation}
As for equations \eqref{eq  restrict a08 e  p eq2s2 t purt1avbf0jjj0} and \eqref{eq  restrict a08 e  p eq2s2 t purt1avbf0jjj1}, when $p=q$, if $A=1$ or $A=\overline{\kappa}^2$, the degree of the numerator in \eqref{eq  restrict a06 e  p eq1  e0tr0nv0} or in \eqref{eq  restrict a06 e  p eq1  e0tr0nv1} decreases and by Lemma~\ref{keylemma1} it decreases by $2$. This gives the seventh equation of Theorem~\ref{Theorem  2c}.

When $\alpha_1+\beta_1=0$ and $\alpha_2+\beta_2=0$, we may let $\alpha_1=\kappa$ and $\alpha_2=1$ by doing a linear transformation $f\to \alpha_2f$.
By considering $\overline{f}^2-1$ and $\overline{f}^2-\overline{\kappa}^2$ similarly as in the first case, respectively, we have
\begin{equation}\label{eq  restrict a06 e  p eq1  e0tr0jhg ue0}
\overline{f}^2-1= \frac{P_1(z,f)^2}{Q_0(z,f)^2(f+\kappa)(f+1)},
\end{equation}
and
\begin{equation}\label{eq  restrict a06 e  p eq1  e0tr0jhg ue1}
\overline{f}^2-\overline{\kappa}^2= \frac{P_2(z,f)^2}{Q_0(z,f)^2(f+\kappa)(f+1)},
\end{equation}
where $P_1(z,f)$ and $P_2(z,f)$ are two polynomials in $f$ with simple roots only and none of these roots is equal to $\pm \kappa$ or $\pm 1$. As for equations \eqref{eq  restrict a08 e  p eq2s2 t purt1avbf0jjj0} and \eqref{eq  restrict a08 e  p eq2s2 t purt1avbf0jjj1}, when $p=q$, if $A=1$ or $A=\overline{\kappa}^2$, the degree of the numerator in \eqref{eq  restrict a06 e  p eq1  e0tr0jhg ue0} or in \eqref{eq  restrict a06 e  p eq1  e0tr0jhg ue1} decreases and by Lemma~\ref{keylemma1} it decreases by $2$. This gives the eighth equation of Theorem~\ref{Theorem  2c}.

\vskip 4pt

\noindent\textbf{Subcase~2:} Equation \eqref{eq  restrict a06 e0}.

\vskip 4pt

We have $\alpha_1+\alpha_2=0$ and $\alpha_3+\alpha_4=0$. We may suppose that $\alpha_1=1$ and $\alpha_3=\kappa$ by doing a linear transformation $f\to \alpha_1f$. Now we have
\begin{equation}\label{eq  restrict a06 e  p eq1   a0}
\overline{f}^2 = \frac{P_0(z,f)^2}{Q_0(z,f)^2}(f^2-1)(f^2-\kappa^2),
\end{equation}
and, further, by applying the analysis after equation \eqref{eq  restrict a06 e  p eq1  e} to $\overline{f}^2-1$ and $\overline{f}^2-\overline{\kappa}^2$, respectively, we have
\begin{equation}\label{eq  restrict a06 e  p eq1   a1}
\overline{f}^2-1= \frac{P_1(z,f)^2}{Q_0(z,f)^2},
\end{equation}
and
\begin{equation}\label{eq  restrict a06 e  p eq1   a2}
\overline{f}^2-\overline{\kappa}^2= \frac{P_2(z,f)^2}{Q_0(z,f)^2},
\end{equation}
where $P_1(z,f)$ and $P_2(z,f)$ are two polynomials in $f$ with simple roots only and none of these roots is equal to $\pm \kappa$ or $\pm 1$. As for equations \eqref{eq  restrict a08 e  p eq2s2 t purt1avbf0jjj0} and \eqref{eq  restrict a08 e  p eq2s2 t purt1avbf0jjj1}, when $p=q$, if $A=1$ or $A=\overline{\kappa}^2$, the degree of the numerator in \eqref{eq  restrict a06 e  p eq1   a1} or in \eqref{eq  restrict a06 e  p eq1   a2} decreases and by Lemma~\ref{keylemma1} it decreases by $2$. This gives the ninth equation of Theorem~\ref{Theorem  2c} and also completes the proof.

\end{proof}

\section{Discussion}\label{discussions}

In sections~\ref{Proof 1} and \ref{Proof 2}, we gave a classification of equation \eqref{yanagiharaeq0} under the assumptions that equation \eqref{yanagiharaeq0} has a transcendental meromorphic solution and the degree of $R(z,f)$ in $f$ satisfies $d\not=n$. These results together with the main theorem in our previous paper~\cite{Korhonenzhang2020}, where the case $d=n$ of equation \eqref{yanagiharaeq0} was considered, provide a complete difference analogue of Steinmetz' generalization of Malmquist's theorem. The classification in sections~\ref{Proof 1} and \ref{Proof 2} is according to the number $N_c$ of the roots $\alpha_i$ in \eqref{P} and $\beta_j$ in \eqref{Q} and whether some of these roots is zero. We did this by mainly using five lemmas, i.e., Lemmas~\ref{basiclemma}--\ref{keylemma1} in section~\ref{Proof 0}. From their proofs, we see that with some simple adjustments these lemmas also apply to the case $d=n$ of equation \eqref{yanagiharaeq0}. In~\cite{ZhangKorhonen:2022}, we have shown how to simplify the proof of the main theorem in~\cite{Korhonenzhang2020}.

We have shown that if equation \eqref{yanagiharaeq0} with $d\not=n$ has a transcendental meromorphic solution, then \eqref{yanagiharaeq0} reduces into one in a list of 17 equations. In the beginning of section~\ref{Proof 0}, we point out that equation \eqref{yanagiharaeq0} may reduce into \eqref{yanagiharaeq0 fd} in some special cases. In section~\ref{Proof 1}, we consider the case where $q=0$; from the results in Theorems~\ref{Theorem  1a} and~\ref{Theorem  1c}, we see that the polynomial term $P(z,f)$ takes particular form and the solutions $f$ are expressed in terms of exponential type functions explicitly. In section~\ref{Proof 2}, we consider the case where $q\geq 1$. In this case, if $n>d$ or $3\leq n<d$, then from Theorems~\ref{Theorem  2a} and~\ref{Theorem  2b} we see that solutions of \eqref{yanagiharaeq0} are also expressed in terms of exponential type functions. But for the case $n=2$ and $n<d$, equation \eqref{yanagiharaeq0}, as well as its solutions, becomes much more complicated. When $q\geq 1$, $n=2$ and $n\nmid |p-q|$, the polynomials $P_0(z,f)$ and $Q_0(z,f)$ are determined and the solutions are clear, as seen in Theorem~\ref{Theorem  2b}. When $q\geq 1$, $n=2$ or $n=3$ and $n\mid |p-q|$, we obtain nine equations in Theorem~\ref{Theorem  2c}. Below, we discuss the eight equations \eqref{Theorem  2c  eq12 first0}, \eqref{Theorem  2c  eq12 third0}, \eqref{Theorem  2c  eq12 fourth0}, \eqref{Theorem  2c  eq17}, \eqref{Theorem  2c  eq0}, \eqref{Theorem  2c  eq12 seventh0}, \eqref{Theorem  2c  eq12 eightth0} and \eqref{Theorem  2c  eq12 ninth0} in the autonomous case.

Solutions to the equations \eqref{Theorem  2c  eq17} and \eqref{Theorem  2c  eq0} are Weierstrass elliptic functions, composed with entire functions. Below we show their relations with the Fermat type equation $h(z)^3+g(z)^3=1$; see \cite{Baker1966,Gross1966erratum,Gross1966}. All meromorphic solutions of the Fermat type equation $h^3+g^3=1$ can be represented as: $h=H(\varphi)$, $g=\eta G(\varphi)=\eta H(-\varphi)=H(-\eta^2\varphi)$, where $\varphi=\varphi(z)$ is an entire function and $\eta$ is a cubic root of~1, and
\begin{equation}\label{Weierstrassellptic0}
H(z)=\frac{1+\wp'(z)/\sqrt{3}}{2\wp(z)}, \quad G(z)=\frac{1-\wp'(z)/\sqrt{3}}{2\wp(z)}
\end{equation}
is a pair of solutions of the Fermat equation $H^3+G^3=1$ with $\wp(z)$ being the particular Weierstrass elliptic function such that $\wp'(z)^2=4\wp(z)^3-1$. For equation \eqref{Theorem  2c  eq0}, we have
\begin{equation*}
\overline{f}^3 +\left[-\frac{P_1(f)}{Q_0(f)}\right]^3=1.
\end{equation*}
Therefore, we have $\overline{f}=H(\phi_1)$ and $P_1(f)/Q_0(f)=-\eta G(\phi_1)$, where $\phi_1=\phi_1(z)$ is an entire function, and $H(z)$ and $G(z)$ are defined as in equation \eqref{Weierstrassellptic0}. Moreover, there exist two constants $A_1\not=0$ and $B_1$ dependent on the coefficients of $P(f)$ and $Q(f)$ such that $\overline{\phi}_1=A_1\phi_1+B_1$. On the other hand, equation \eqref{Theorem  2c  eq17} can also be transformed into the Fermat type equation in the following way: Recall from the proof that we have equation \eqref{eq  restrict a1 fgtrytr} and $1,\eta_1,\eta_1^2$, where $\eta_1$ is a cubic root of~1 such that $\eta_1^2+\eta_1+1=0$, are completely ramified values of $f$ with multiplicities~3. We let $(f-\eta_1^2)/(f-\eta_1)=g^3$. Then $g$ is a meromorphic function and it follows that $f=(\eta_1 g^3-\eta_1^2)/(g^3-1)$. By substituting this equation into \eqref{eq  restrict a1 fgtrytr}, we get
\begin{equation*}
\overline{f}^3+\left[-\frac{P_{01}(g^3)g}{Q_{01}(g^3)}\right]^3=1,
\end{equation*}
where $P_{01}(g^3)$ and $Q_{01}(g^3)$ are two polynomials in $g$ with no common roots. Then we have $\overline{f}=H(\phi_1)$ and $P_{01}(g^3)g/Q_{01}(g^3)=-\eta G(\phi_2)$, where $\phi_2=\phi_2(z)$ is an entire function, and $H(z)$ and $G(z)$ are defined as in \eqref{Weierstrassellptic0}. Moreover, there exist two constants $A_2\not=0$ and $B_2$ dependent on the coefficients of $P_0(f)$ and $Q_0(f)$ such that $\overline{\phi}_2=A_2\phi_2+B_2$.

Solutions to the six equations \eqref{Theorem  2c  eq12 first0}, \eqref{Theorem  2c  eq12 third0},
\eqref{Theorem  2c  eq12 fourth0}, \eqref{Theorem  2c  eq12 seventh0}, \eqref{Theorem  2c  eq12 eightth0} and \eqref{Theorem  2c  eq12 ninth0} are Jacobian elliptic functions, composed with entire functions. Below we show their relations with the symmetric biquadratic equation of the form $x^2y^2-(x+y)+c^2=0$; see \cite[p.~471]{Baxter1982}). For equations \eqref{Theorem  2c  eq12 seventh0}, \eqref{Theorem  2c  eq12 eightth0}, \eqref{Theorem  2c  eq12 ninth0}, from the two equations \eqref{eq  restrict a06 e  p eq1  e0tr0nv0} and \eqref{eq  restrict a06 e  p eq1  e0tr0nv1}, or from the two equations \eqref{eq  restrict a06 e  p eq1  e0tr0jhg ue0} and \eqref{eq  restrict a06 e  p eq1  e0tr0jhg ue1}, or from the two equations \eqref{eq  restrict a06 e  p eq1   a1} and \eqref{eq  restrict a06 e  p eq1   a2}, we get an equation of the following form:
    \begin{equation*}
   \frac{\overline{f}^2-1}{\overline{f}^2-\kappa^2}=\frac{P_1(f)^2}{P_2(f)^2},
   \end{equation*}
where $P_1(f)$ and $P_2(f)$ are two polynomials in $f$ with simple roots only and with no common roots.
Denoting $R_1(f)=P_1(f)/P_2(f)$, it follows that
   \begin{equation}\label{Theorem  2c  eq579toge1b}
   \overline{f}^2R_1(f)^2-[\overline{f}^2+R_1(f)^2]+\kappa^2=0,
   \end{equation}
which is a symmetric biquadratic equation in $\overline{f}$ and $R_1$. The equation above can be solved as $\overline{f}=k_1^{1/2}\text{sn}(\varphi_1(z)\pm \tau_1)$ and $R_1(f)=k_1^{1/2}\text{sn}(\varphi_1(z))$, where $k_1$ and $\tau_1$ are two parameters dependent on the constant $\kappa^2$, $\text{sn}(\varphi_1)$ is the Jacobian elliptic function with modulus $k_1$ and $\varphi_1=\varphi_1(z)$ is an entire function. Then there are two constants $C_1\not=0$ and $D_1$ such that $\overline{\varphi}_1=C_1\varphi_1+D_1$. On the other hand, equations \eqref{Theorem  2c  eq12 first0}, \eqref{Theorem  2c  eq12 third0} and \eqref{Theorem  2c  eq12 fourth0} can also be transformed into symmetric biquadratic equations similar to \eqref{Theorem  2c  eq579toge1b} in the following way: For equation \eqref{Theorem  2c  eq12 first0}, we let $(f+\kappa)/(f+1)=g^2$; for equation \eqref{Theorem  2c  eq12 third0} we let $(f+\gamma)/(f-\gamma)=g^2$; for equation \eqref{Theorem  2c  eq12 third0} we let $(f+\gamma)/(f-\gamma)=g^2$ when $\kappa=-1$ and let $(f+\kappa)/(f+1)=g^2$ when $\kappa\not=-1$. From the proof in Theorem~\ref{Theorem  2c}, we see that each $g$ in the above expressions is a meromorphic function. By writing $f$ in terms of $g^2$, then from the two equations \eqref{eq  restrict a08 e  p eq2s2 t purt1avbf0jjj0} and \eqref{eq  restrict a08 e  p eq2s2 t purt1avbf0jjj1}, or from the two equations \eqref{eq  restrict a08 e  p eq2s2 t} and \eqref{eq  restrict a08 e  p eq2s2 t purt0a  pp}, or from the two equations \eqref{eq  restrict a1 gfhj092} and \eqref{eq  restrict a1 gfhj093}, or from the two equations \eqref{eq  restrict a1 gfhj090} and \eqref{eq  restrict a1 gfhj091}, we get an equation of the following form:
    \begin{equation*}
   \frac{\overline{f}^2-\gamma^2}{\overline{f}^2-1}=\frac{P_{01}(g)^2}{P_{02}(g)^2},
   \end{equation*}
where $P_{01}(g)$ and $P_{02}(g)$ are two polynomials in $g$ with simple roots only and with no common roots. Denoting $R_2(g)=P_{01}(g)/P_{02}(g)$, it follows that
   \begin{equation*}
   \overline{f}^2R_2(g)^2-[\overline{f}^2+R_2(g)^2]+\gamma^2=0,
   \end{equation*}
which is a symmetric biquadratic equation in $\overline{f}$ and $R_2$. Then we have $\overline{f}=k_2^{1/2}\text{sn}(\varphi_2(z)\pm \tau_2)$ and $R_2(g)=k_2^{1/2}\text{sn}(\varphi_2(z))$, where $k_2$ and $\tau_2$ are two parameters dependent on the constant $\gamma^2$, $\text{sn}(\varphi_2)$ is the Jacobian elliptic function with modulus $k_2$ and $\varphi_2=\varphi_2(z)$ is an entire function. Moreover, there are two constants $C_2\not=0$ and $D_2$ such that $\overline{\varphi}_2=C_2\varphi_2+D_2$.

\end{document}